\documentclass[a4paper,fleqn,leqno]{article}

\usepackage[utf8]{inputenc}
\usepackage[english]{babel}
\usepackage{latexsym,amssymb,amsthm,color}
\usepackage[centertags]{amsmath}
\usepackage{mathtools}
\usepackage{graphicx}
\usepackage{bbm}
\usepackage{framed}
\usepackage{pdflscape}
\usepackage{hyperref}
\usepackage[normalem]{ulem}
\newcommand{\Rand}[1]{\marginpar{#1}}
\renewcommand{\Rand}[1]{}

\newcommand{\be}[1]{\Rand{\vspace{0,6cm}\tt #1}\begin{equation}\label{#1}}
\newcommand{\beL}[1]{\Rand{\vspace{0,6cm}\tt #1}\begin{lemma}\label{#1}}
\newcommand{\belC}[2]{\Rand{\vspace{0,6cm}\tt #1}\begin{lemma}[#2]\label{#1}}
\newcommand{\beP}[1]{\Rand{\vspace{0,6cm}\tt #1}\begin{proposition}\label{#1}}
\newcommand{\bePC}[2]{\Rand{\vspace{0,6cm}\tt #1}\begin{proposition}[#2]\label{#1}}
\newcommand{\beD}[1]{\Rand{\vspace{0,6cm}\tt #1}\begin{definition}\label{#1}}
\newcommand{\beT}[1]{\Rand{\vspace{0,6cm}\tt #1}\begin{theorem}\label{#1}}
\newcommand{\beC}[1]{\Rand{\vspace{0,6cm}\tt #1}\begin{corollary}\label{#1}}
\newcommand{\ep}{\end{proposition}}
\newcommand{\bea}[1]{\Rand{\vspace{0,7cm}\tt #1\vspace{-0,7cm}}\begin{eqnarray}\label{#1}}

\newcommand{\wt}{\widetilde}
\def\CA{\mathcal{A}}
\def\CB{\mathcal{B}}
\def\CD{\mathcal{D}}

\def\CM{\mathcal{M}}

\def\CH{\mathcal{H}}
\def\CL{\mathcal{L}}

\def\CU{\mathcal{U}}

\def\CX{\mathcal{X}}

\def\A{\mathbb{A}}
\def\G{\mathbb{G}}
\def\K{\mathbb{K}}
\def\N{\mathbb{N}}
\def\M{\mathbb{M}}
\def\R{\mathbb{R}}
\def\U{\mathbb{U}}
\def\Z{\mathbb{Z}}

\def\mfw{\mathfrak{w}}
\def\mfy{\mathfrak{y}}
\def\mfh{\mathfrak{h}}
\def\mfu{\mathfrak{u}}
\def\mfA{\mathfrak{A}}

\def\mfH{\mathfrak{H}}
\def\mfW{\mathfrak{W}}
\def\mfL{\mathfrak{L}}
\def\mfU{\mathfrak{U}}

\def\mfD{\mathfrak{D}}
\def\mfY{\mathfrak{Y}}
\def\mfX{\mathfrak{X}}
\def\mfZ{\mathfrak{Z}}
\def\mfN{\mathfrak{N}}

\DeclareMathOperator{\supp}{supp}

\newcommand{\la}{\longrightarrow }
\newcommand{\La}{\Longrightarrow }
\newcommand{\Ltoo}{{_{\D \longrightarrow \atop L \to \infty}}}
\newcommand{\Lto}{{_{\D \Longrightarrow \atop L \to \infty}}}
\newcommand{\kto}{{_{\D \Longrightarrow \atop k \to \infty}}}

\newcommand{\tto}{{_{\D \Longrightarrow \atop t \to \infty}}}

\newcommand{\nto}{{_{\D \Longrightarrow \atop n \to \infty}}}
\newcommand{\dto}{{_{\D \to \atop d \to \infty}}}

\newcommand{\D}{\displaystyle}
\newcommand{\wh}{\widehat}

\newcommand{\dr}{\underline{\ur}}
\newcommand{\uur}{\underline{\ur}}

\newcommand{\uc}{\underline{c}}
\newcommand{\ulam}{\underline{\lambda}}
\newcommand{\uLam}{\underline{\Lambda}}

\newcommand{\ur}{\underline{r}}

\newcommand{\uv}{\underline{v}}

\newcommand{\smallu}{\mathpzc{u}}
\DeclareMathAlphabet{\mathpzc}{OT1}{pzc}{m}{it}

\newcommand{\dd}{\mathrm{d}}
\newcommand{\ee}{\end{equation}}
\newcommand{\eea}{\end{eqnarray}}
\newcommand{\bean}{\begin{eqnarray*}}
\newcommand{\eean}{\end{eqnarray*}}

\newcommand{\ve}{\varepsilon}

\newcounter{secnum}
\numberwithin{equation}{section}

\makeatletter \@addtoreset{equation}{section} \makeatother
\makeatletter \@addtoreset{enunciato}{section} \makeatother
\newcounter{enunciato}[section]

\newtheorem{ittheorem}{Theorem}
\newtheorem{itlemma}{Lemma}
\newtheorem{itproposition}{Proposition}
\newtheorem{itdefinition}{Definition}
\newtheorem{itremark}{Remark}
\newtheorem{itexample}{Example}
\newtheorem{itclaim}{Claim}
\newtheorem{itcorollary}{Corollary}

\newenvironment{theorem}{\addtocounter{enunciato}{1}
\begin{ittheorem}}{\end{ittheorem}}
\newenvironment{lemma}{\addtocounter{enunciato}{1}
\begin{itlemma}}{\end{itlemma}}
\newenvironment{proposition}{\addtocounter{enunciato}{1}
\begin{itproposition}}{\end{itproposition}}
\newenvironment{definition}{\addtocounter{enunciato}{1}
\begin{itdefinition}}{\end{itdefinition}}
\newenvironment{remark}{\addtocounter{enunciato}{1}
\begin{itremark}}{\end{itremark}}

\newenvironment{corollary}{\addtocounter{enunciato}{1}
\begin{itcorollary}}{\end{itcorollary}}

\usepackage{stackengine}

\begin{document}

\title{Continuum space hierarchical Cannings process:\\ 
II.\ Genealogies and ancestral web}

\author{A. Greven$^1$, F. den Hollander$^2$}

\date{\today}

\maketitle

\begin{abstract}
The present paper is the second in a series of three papers in which we consider a \emph{spatial multi-type population model} arising as the many-individuals-per-site limit of individual-based models. The evolution mechanisms in this model are: migration, Fleming-Viot resampling, $\Lambda_0$-Cannings resampling, and $(\Lambda_k)_{k \in \N}$-Cannings block-resampling. This evolution leads to a collection of \emph{measure-valued} processes indexed by a continuum geographic space $\Omega_N^\infty$, called the \emph{continuum hierarchical group}, obtained as the limit of a discrete hierarchical group $\Omega_N^L$ as $L\to\infty$, where $\Omega_N^L = \oplus_{k \geq -L} Z_N$ with $Z_N$ the cyclical group of order $N$. On this geographic space, each component is a probability measure on a type space $\K$. 

In the first paper \cite{gh1}, we \emph{constructed} the measure-valued process on $\Omega_N^\infty$ describing the type frequencies in the population and identified its main properties. We found a \emph{dichotomy} between this process being \emph{stochastic} (= locally monotype) or \emph{deterministic} (= locally multitype). In the present paper we turn to the \emph{genealogies} of the population on $\Omega_N^\infty$. Our goal is to construct and analyse a dynamics analogous to the genealogical Cannings and Fleming-Viot processes on the discrete hierarchical group $\Omega_N$. We show that the existence of a stochastic continuum-space evolution of the type-frequency process allows for the construction of the $\U^{\Omega_N^\infty \times \K}$-valued \emph{stochastic marked genealogy process}. This construction is based on an approximation by similar processes on $\Omega_N^L$ as $L \to \infty$. In the third paper \cite{gh3}, we look at the \emph{hierarchical mean-field limit} of the processes studied in \cite{gh1} and the present paper, which is a measure-valued process on $\Omega_\infty^\infty$ obtained via approximation by similar processes on $\Omega_N$ as $N\to\infty$. In this limit a number of simplifications occur that allow for the identification of several \emph{universality classes}.

The Cannings process on $\Omega_N^\infty$ exhibits a \emph{dichotomy} also for the \emph{genealogies}, with two regimes. In regime~(1), the population descends, viewed locally and for time-expanding sets of sites, from a single ancestor. This gives a non-trivial stochastic limit dynamics on $\U^{\Omega_N^\infty}$ as well as on the projected type frequencies. In regime~(2), the population of descendants decomposes into individuals descending from different ancestors observed locally, and the evolution of type frequencies is deterministic. In this regime the genealogical structure is degenerate: the family decomposition is dust. Regime~(1) occurs for point-recurrent migration, and for Harris-recurrent migration when block-resampling is present. Regime~(2) occurs for transient migration. For recurrent migration we obtain, for each positive time, a decomposition of $\Omega_N^\infty$ into random balls where the population is mono-type.

A new element in the paper is the construction of the \emph{ancestral web}. This is a forward-evolving population of ancestral paths (called lineages in population genetics) in ancestral-geographic-type space. It replaces populations of individuals by \emph{populations of paths}, and embeds two objects: the genealogical process on $\U$ and the measure-valued process on $\Omega_N^\infty \times \K$. The latter takes values in the space of measures on ancestral paths marked with types and locations. This embedding gives more information on the evolution of the population. Functional and strong duality play a key role in the analysis of the ancestral web. The dual is constructed with the help of the L\'evy-web, the analogue of the Brownian web. 

In \cite{gh1} we carried out a continuum-space rescaling limit (of time, space and rates) on $\Omega_N$, obtaining a special process on $\Omega_N^\infty$: the \emph{stable process}. In the present paper this is done also at the level of genealogical processes. The resulting class of processes identifies the universality classes for the genealogy processes of the Cannings process with block-resampling, via stable versions of these processes (based on stable migration).
\end{abstract}


\medskip\noindent
{\bf Keywords:} Stochastic processes of evolving genealogies, continuum-space rescaling limit, Cannings process with block-resampling, $\U^V$-valued processes, ancestral web, L\'evy-web, functional duality, strong graphical duality.

\medskip\noindent
{\bf AMS Subject Classification:} 60K35.

\medskip\noindent
{\footnotesize $^{1}$ 
Department Mathematik, Universit\"at Erlangen-N\"urnberg, Cauerstrasse 11, D-91058 Erlangen, Germany\\
greven@mi.uni-erlangen.de}\\
{\footnotesize $^{2}$ 
Mathematisch Instituut, Universiteit Leiden, Einsteinweg 55, 2333CC Leiden, The Netherlands\\
denholla@math.leidenuniv.nl}

\small
\tableofcontents
\normalsize


\section{Introduction}
\label{s.intro}


\subsection{Background and goals}
\label{ss.backgoal}

{\bf Basic population models.}
In \cite{GHKK14} we introduced \emph{infinite population} per site versions of \emph{Fleming-Viot} and \emph{Cannings} processes on a specific countable geographic space, the hierarchical group $\Omega_N$, as measure-valued (precisely: $\CM_{\rm loc}(\Omega_N \times \K)$-valued) processes. Types are drawn from a compact set $\K$, while the geographic space is modeled by the hierarchical group 
\begin{equation}
\Omega_N =\bigoplus\limits_{\N_0} Z_N, 
\end{equation}
with $Z_N$ the cyclical group of order $N \in \N$. This results in a system indexed by $\Omega_N$ of interacting measure-valued processes (probability measures on $\K$, interacting via migration), evolving \emph{neutrally with resampling}. The resampling is based on \emph{Fleming-Viot} or Cannings resampling at single sites, and additionally on Cannings \emph{block-resampling} with \emph{reshuffling} acting in \emph{balls} in $\Omega_N^\infty$ of all sizes. This was introduced in \cite{GHKK14} and is analogous to models on $\R^d$ with resampling of the population in random balls studied in~\cite{EVY19,Etheridge2024}. The motivation is to have resampling with \emph{local interaction} between different sites, not only migration between sites.

More specifically, in the individual-based model the process evolves via migration and change of generation, based on Fleming-Viot resampling, $\Lambda$-Cannings resampling at the level of sites (colonies), and \emph{block-resampling} with \emph{uniformly distributed reshuffling} in \emph{$k$-balls} for all ball sizes. The infinite-population limit dynamics is characterised by: an effective migration kernel $a(\cdot,\cdot)$, a volatility constant $d \in [0,\infty)$ of the Fleming-Viot diffusion, $\Lambda \in \CM([0,1])$ for the Cannings single-site resampling, and a sequence of Cannings resampling measures $(\Lambda_k)_{k \in \N}$ on $[0,1]$ for block-resampling.

Since the theory of the $\Omega_N$-model extends to the countable $\Omega_N^L$-model (obtained by splitting sites in $\Omega_N$ successively $L$ times into $N$ finer sites), we were able in~\cite{gh1} to construct the \emph{continuum-space} version, the $\CM_{\rm loc}(\Omega_N^\infty \times \K)$-valued process, with $\Omega_N^\infty$ an \emph{uncountable} and \emph{Polish} space. This is the continuum-space version of the discrete geographic space model on $\Omega_N$. The key finding of that paper was an \emph{approximation procedure} by processes defined on successively refined countable geographic spaces $\Omega_N^L$. It showed that single-site Cannings resampling produces an extra Fleming-Viot single-site effect, while only Cannings \emph{block-resampling} with \emph{reshuffling} remains as a jump component in the processes on $\Omega_N^\infty$. The motivation for looking at the \emph{continuum-space limit} and at the \emph{hierarchical} group as geographic space rather than $\R^d$ was outlined in \cite{gh1} for the {\em measure-valued} process, and applies here as well.


\medskip\noindent{\bf Stochastic genealogical processes.}
The \emph{marked genealogy valued} process on $\Omega_N$ without block-resampling is constructed in~\cite{gklimw} on the state space $\U^{\Omega_N \times \K}$ introduced in~\cite{DGP11}. With the population alive at time $t$ one can associate their \emph{genealogy} (the genealogical distances between pairs of individuals), which evolves as a \emph{stochastic process} in $t$. This is called the marked ultrametric measure space valued dynamics, or more precisely: the $\U^{\Omega^{\infty}_N}$-valued Cannings process $\mfU$. In the present paper we construct this process of \emph{evolving genealogies} on $\Omega_N^\infty$. By introducing marks (locations and types of individuals) we then get the \emph{marked} genealogies in $\U^{\Omega_N^\infty \times \K}$.

The results on the measure-valued version are helpful for deriving qualitative properties of the genealogies. Conversely, the $\U^{\Omega_N^\infty}$-valued process helps understand aspects of the measure-valued process, such as dichotomies for the long-time behaviour on $\Omega_N$ and the formation of clusters of monotype regions. The continuum limit of the {\em genealogical} process as a marked metric measure space-valued process exists only if we have a {\em stochastic} continuum limit for the {\em measure-valued} process on $\wt \CM(\Omega_N \times \K)$, a subset of $\CM_{\rm loc} (\Omega_N \times \K)$ to be introduced.

\medskip\noindent
{\bf The ancestral web.}
As a \emph{new} element for continuum-space models, we consider a richer version of the genealogy process. We replace populations of individuals by populations of whole \emph{paths} in ancestral-type-geographic-space. Such paths are called lineages in population genetics, and they give the analogue of the random graph representing the information of the time-space process of the underlying population evolution. We introduce two key objects. First, the \emph{ancestral web} of~\cite{gklimw,Lambda}: a \emph{population of evolving paths} (the \emph{lineages}) with values in the \emph{ancestral space}. The paths connect \emph{ancestors} at time $0$ to individuals alive at time $t$ (the \emph{descendants}). This view is formalised as the \emph{ancestral coupling}. The paths are marked by type and location. The population is described by a measure on \emph{$(\Omega_N^\infty \times \K)$-marked paths} in the ancestral space. This is the \emph{forward} picture. Second, the \emph{backward} picture, evolving from a time horizon $t$ backwards to time $0$: the ancestral path reversed gives the \emph{path of descent}, studied via the dual process, which is in \emph{strong graphical duality} with the original forward process.

The process $\mfU$ is the backbone of the evolution of genealogies viewed as populations of paths in the ancestral space of \emph{all} individuals \emph{ever alive up to time $t$}, including all \emph{fossils}. This object allows tracing the complete \emph{time-space footprint} of the population and its evolution in $T$, with $T$ the \emph{time horizon} of the evolution of the genealogy. The ancestral space $\mfX$ is the space of all \emph{individuals ever lived}, modelled as a \emph{rooted} or more generally \emph{grounded} \emph{tree-like metric fibred measure space}. Here ``fibred'' refers to the time-$t$ slice of the population alive, equipped with a path of measures on the ``alive at time $t$'' fibres. This enrichment is called the \emph{ancestral coupling}, our second key object. It embeds the genealogies of individuals alive at different times $t$, as well as the ancestral web (see~\cite{gklimw}). To obtain the full picture, the ancestral paths are marked with \emph{types} and \emph{locations}.

Both the ancestral web and the ancestral coupling are constructed by means of a \emph{dual} object: the \emph{L\'evy-web}, a population of \emph{backward paths} (the paths of descent) generated by a \emph{dual process}. This dual process is driven by a \emph{coalescing} and \emph{block-coalescing} system of L\'evy-processes on $\Omega_N^\infty$ with \emph{reshuffling}. It is constructed so that both forward and backward processes sit in a single \emph{marked, rooted, metric and fibred measure tree} $\mfX$. The geodesics of this tree are the objects of primary interest and form a ``continuum'' graph. This embedding is referred to as the \emph{string graphical representation} of the genealogy.

For all these objects we identify \emph{universality classes} of their large space-time scale behaviour by considering \emph{time-space-rate rescaling} limits, which represent the possible \emph{universality classes}.


\medskip\noindent
{\bf Goals.}
The main goals of this paper are to address the following three sets of \emph{themes}.

\medskip\noindent
{\bf \emph{{\rm (1)} Genealogy processes:}}
\begin{itemize}
\item 
\emph{Construction} of the process $\mfU=(\CU_t)_{t \geq 0}$ of \emph{genealogies} of the population alive at the current time $t$, evolving with the Fleming-Viot and Cannings process with block-resampling on the continuum space $\Omega_N^\infty$. The construction uses countable geographic space \emph{approximation} of $\Omega_N^\infty$. It is the limit $L \to \infty$ of the genealogy process starting from $\Omega_N$ and then on analogous \emph{countable} spaces of successive $L$-refinements of sites, denoted $\Omega_N^L$.
\item 
Characterisation of the obtained genealogy process $\mfU$ by a \emph{well-posed martingale problem}.
\item 
Criteria for when the continuum-space dynamics exists as a \emph{$\U^{\Omega^{\infty}_N \times \K}$-valued} process that is truly \emph{stochastic} (not a deterministic flow). In particular: the \emph{existence} of continuum-space $\U^{\Omega_N^\infty \times \K}$-valued processes is equivalent to the existence of a \emph{truly stochastic} measure-valued continuum-space process for the projection on $\Omega_N^\infty \times \K$.
\item 
Derivation of characteristic \emph{geometric properties} of genealogies: for each time $t$, a decomposition of the time-$t$ population for $s \in [0,t]$ into depth-$s$ subfamilies in $2s$-balls. For $t>0$: a decomposition of the geographic space $\Omega_N^\infty$ into collections (the founders) of unions of disjoint balls, each union being covered by exactly one \emph{depth-$h$ subfamily} of the population with \emph{random $h \in (0,t)$}. Description of the \emph{time-space} structure of these decompositions in relation to the time structure of the family decompositions in the tree $\mfX$ of all individuals ever alive, leading to the next theme (ancestral web and coupling).
\end{itemize}

\medskip\noindent
{\bf \emph{{\rm (2)} Ancestral web and random graph representation:}}
\begin{itemize}
\item Construction of the \emph{ancestral web} $\mfW=(\mfW_T)_{T \geq 0}$ for the above genealogy process $\mfU$. This is a forward-evolving population of ancestral paths (\emph{lineages}), constructed via its \emph{dual} object, the \emph{L\'evy-web} $\wh \mfW$, which is a population of \emph{dual backward paths} (the \emph{paths of descent}). The dual is identified as the genealogy of a spatial coalescent starting with a path everywhere in $\R \times \Omega_N^\infty$. The process $(\mfU_t)_{t \in [0,T]}$ is embedded in $\mfW_T$.
\item \emph{Graphical construction} of the genealogy process $\mfU$ on $\U^{\Omega^{\infty}_N}$, the ancestral web $\mfW$, and the dual \emph{L\'evy-web} $\wh \mfW$, all embedded in one \emph{continuum random graph}. Both are populations of paths: forward in $\mfW$ and backward in $\wh \mfW$. The paths are embedded in a single \emph{grounded tree-like marked fibred metric measure space}, the \emph{ancestral space} $\mfX^\dagger$. This space represents the population of \emph{all individuals ever alive} with their ancestral relations given by the lineages. Together with the dual paths of descent they form a \emph{continuum graph}. We show that in this space $\mfW$ and $\wh \mfW$ are in \emph{strong graphical duality}.
\end{itemize}

\medskip\noindent
{\bf \emph{{\rm (3)} Universality and continuum space rescaling limit:}}
\begin{itemize}
\item Carrying out on $\Omega_N$ the \emph{continuum space scaling limit} towards $\Omega_N^\infty$ from~\cite{gh1}, now for the \emph{genealogy} process and the \emph{ancestral web/coupling}. This yields the \emph{universality classes} for the \emph{large time-space scale} behaviour for the genealogy models, represented by a $(\R^+)^2$-indexed set of prototype models based on the \emph{stable processes} as migration and their stable ancestral webs and couplings.
\end{itemize}

The key tool of analysis is a \emph{functional duality relation} based on spatial coalescent processes and stronger versions of this duality: the \emph{strong dual representations} and the \emph{graphical} representations of process and dual in one marked \emph{random graph}.

\medskip\noindent
{\bf Past work.}
 For spaces of continuum cardinality, the analysis of evolving genealogies in $\R^1$ \cite{Etheridge2024,EVY19} and for the Fleming-Viot model in~\cite{GSW} has been carried out previously. That analysis was based on the Brownian web, with Brownian motion as migration. However, a characterisation via a well-posed martingale problem was still missing there. Jump migration processes on $\R^1$ seem intractable. We therefore work on $\Omega_N^\infty$ as geographic space, where progress is possible since the underlying migration is a jump process. In~\cite{gklimw} genealogies were studied on geographic spaces that are countable Abelian groups for Cannings models, but without block-resampling. The Cannings version of the historical process for branching processes was introduced in~\cite{D93,DP91}.

For the Cannings process on countable geographic spaces like $\Omega_N$, we can associate, even in the continuum-frequency limit, the process $\mfU$ of $V$-marked genealogies as Markovian dynamics on a Polish space called $\U^V$, introduced in~\cite{DGP11}. Here the points (the individuals of the population) are marked with types and locations, i.e., $V= \K \times \Omega_N^\infty$. The genealogies are equivalence classes of $V$-marked ultrametric measure spaces. The metric is the genealogical distance of two individuals, which is twice the time back to the most recent common ancestor.

The ancestral web is based on a graphical construction of the individual-based model. It considers the time-space process of the evolution, coded in a random graph, adapted to the continuum-frequency setting. The key tool is duality for all the processes. The graphical construction uses a dual graph generated by a dual process called the L\'evy-web, which we define in a main part of this paper. We construct the ancestral web in the continuum-frequency and continuum-space model on $\Omega_N^\infty$, in the same spirit as the genealogical version of the Brownian web on $\R^1$ studied in~\cite{GSW,NRS05}. The result is an evolving (in time) random graph of forward lineages (ancestral paths) and backward marked paths of descent, embedded in what we call the ancestral space. This is a (random) grounded tree-like metric fibred measure space, where ``fibred'' means equipped with a collection of time-indexed resampling measures on the ancestral paths alive at a specific time. This concept was introduced in~\cite{gklimw,Lambda} for countable discrete geographic spaces. We embed the forward and the backward picture in this ancestral space as forward and backward geodesics. This yields a new object and a strong graphical duality of the forward and backward time-space processes of location-type marked genealogies.


\subsection{Basic models on the hierarchical groups $\Omega_N$, $\Omega_N^L$, $\Omega_N^\infty$: notation and dynamics}
\label{s.basic}

In this section we recall the \emph{geographic} spaces $\Omega_N,\Omega_N^L$ ($L \in \N$) and $\Omega_N^\infty$, the basic \emph{evolution mechanism} of the Cannings model with block-resampling of~\cite{GHKK14} and~\cite{gklimw}, and give the mechanisms of the evolution and its dual. We recall basic properties of the $\CM(\Omega_N^\infty \times \K)$-valued process from~\cite{gh1} \emph{whenever they are used}, with precise references and further explanation. We refer to \cite[Section 1.3]{gh1} for full definitions and notations. 

\medskip\noindent
{\bf Geographic spaces and the Cannings process.}
The \emph{geographic spaces} are the following topological groups:
\begin{align}
&\Omega_N= \oplus_{k \in \N_0} Z_N,  \text{ with $Z_N$ the cyclical group $\{0,1,\ldots,N-1\}$,}\\
&\Omega_N^L=\oplus_{k \geq -L} Z_N, \,k \in \Z,L \in \N,\\
&\Omega_N^\infty = \text{ closure of $Z_N^{\Z,+}$},
\end{align}
with $Z_N^{\Z,+}$ the set of sequences with values in $Z_N^\Z$ being $0$ above some finite index, equipped with the ultrametric of the genealogical distance. The genealogical distance equals $k$ for $k \geq 0$ and $\mathrm{e}^{k-1}$ for $k<0$, where $k$ is the index where two sequences first agree for all entries above. See Section 1.3 in~\cite{gh1} for more details.

On $\Omega_N$ and $\Omega_N^L$ we introduce a $((c_k,\lambda_k,\Lambda_k)_{k\in\Z}, (\Lambda_0,d))$-\emph{Cannings process}, extended to $\Omega_N^L$ by taking $k \in \{-L,-L+1,\ldots\}$, with:
\begin{itemize}
\item \emph{migration} characterised by jump rates $c_k$ for size-\emph{$k$} jumps in $\Omega_N^\infty$,
\item \emph{block-resampling} with \emph{block reshuffling} at rates $\lambda_k$ with $\Lambda_k$-block-resampling for $k \in \N$ and single-site $\Lambda_0$-Cannings resampling respectively rate-$d$ Moran resampling, with $k$ running from $-L,-L+1$ to infinity for $\Omega_N^L$.
\end{itemize}
The processes for $k \in \Z$ with $k \geq -L$ are the processes on $\Omega_N^L$. The full set with $k \in \Z$ is the process we want to define on $\Omega_N^\infty$, i.e., with parameters:
\begin{equation}\label{e302}
\left((c_k,\lambda_k,\Lambda_k)_{k \in \Z},\left(\Lambda_{-\infty},d\right)\right),
\end{equation}
where for $L=\infty$ the last entry becomes \emph{Fleming-Viot} resampling with rate $d+|\Lambda_{-\infty}|$. The concept of \emph{Moran resampling} (often called \emph{Fleming-Viot resampling}), \emph{$\Lambda$-resampling}, and their \emph{$(\Lambda_k)_k$-block-resampling} versions are explained in detail in Section 2.1 of~\cite{gh1}. We also look at the process where $d + |\Lambda_{-\infty}| \to \infty$, which yields a different class of models.

\medskip\noindent
{\bf Duality.}
A key role is played by \emph{duality}. The dual process is driven by a \emph{spatial coalescent}: the $((c_k)_{k \in \Z}, d, \Lambda,(\Lambda_k)_{k \in \Z})$-coalescent, with Kingman coalescence at rate $d$ and Cannings coalescence with $\Lambda=0$ (respectively for finite $L$: rate $d$ and $\Lambda_{-L}$-Cannings coalescence at each site), and furthermore $(\Lambda_k)_{k \in \Z}$ $k$-ball coalescence with uniform $k$-ball reshuffling, with parameters given by~\eqref{e302}. The backward-time \emph{coalescence transitions} are dual to the \emph{resampling transitions} in the forward process.

The parameters allow defining the effective migration via a \emph{Levy-process}: observing a single individual and its location, one obtains an \emph{effective} Levy-process moving by migration at rates $(c_k)_{k \in \Z}$ and reshuffling at rates $(\lambda_k)_{k \in \Z}$, independently.


\subsection{Review of modelling genealogies}
\label{ss.gene}

The study of genealogies of populations started in the 1980s. For extensive references we refer to~\cite{DG23}. We give an intuitive picture of the \emph{evolving} genealogies as stochastic processes and their mathematical description with state spaces of marked metric measure spaces.


\subsubsection{Background on genealogies}
\label{sss.background}

We now explain the description of {\em genealogies} used here. Return to the \emph{individual-based} Cannings model. In the resampling transition we replace a number of individuals marked ``$1$'' (the \emph{ancestors}) by the children of a single individual (the \emph{descendants}). We view this as a \emph{birth} event: one {\em ancestor} dies and is replaced by offspring called his {\em descendants}, and others are replaced by his descendants and simply die. This allows defining a genealogical distance between two individuals as twice the time back to the {\em most recent common ancestor}. The population therefore has a \emph{marked genealogical} structure, since the concepts of \emph{ancestors} and \emph{descendants} induce an order relation in the ancestral space $\mfX$, and since individuals come with \emph{locations} and \emph{types} in the individual-based version. This structure is inherited in \emph{infinite}-population models, which are the focus of this paper.

To obtain a Polish space, take equivalence classes of {\em marked ultrametric measure spaces}, a concept from \cite{GPW09} and \cite{DGP11}. One defines a genealogical distance (twice the time back to the most recent common ancestor of two or more individuals), assigns a \emph{mark} to every individual with its (type, location) information, and considers as state space the set of {\em equivalence classes} of these {\em $V$-marked ultrametric measure spaces} induced by measure- and mark-preserving isometries. Here $V=\K \times \Omega_N$ or $V=\K \times \Omega_N^\infty$, where $\K$ is the \emph{compact} type space.

The basic object is the \emph{state space} of genealogies, denoted
\begin{equation}\label{e295}
\text{$\U_1^V$, \quad if $(V,r_V)$ is the metric \emph{mark} space.}
\end{equation}
We use this for the Cannings model with block-resampling and its continuum-frequency and subsequent continuum-space limit.

In the processes considered in this paper, the stochastic process of marked genealogies arises as the many-individuals-per-site limit of the individual-based models (see~\cite{GPWmp13},~\cite{gklimw}). We refer to this as the {\em genealogical spatial Cannings process on $\Omega_N^\infty$}, denoted $\mfU$. Genealogical processes can be characterised by well-posed martingale problems if the geographic space is $\Omega_N$ (see \cite{GPWmp13} and \cite{GSW} for the Fleming-Viot process, and \cite{gklimw} for the Cannings process; we only need to extend the latter to incorporate \emph{block} resampling). We now have to achieve this on the uncountable geographic space $\Omega_N^\infty$ instead of the discrete countable $\Omega_N$.

We recall below briefly which facts are needed on marked metric measure spaces, before passing to giving the dynamics of the $\U^{\Omega_N^\infty \times \K}$-valued process with locations in $\Omega^{\infty}_N$. We state in the next section that the continuum-space process is well defined.


\subsubsection{Marked metric measure spaces}
\label{sss.mmmsp}

We review the basic concepts from~\cite{GPW09,DGP11} and~\cite{GSW} needed for studying \emph{evolving} genealogies. For a more detailed survey, see \cite{DG23}.

\medskip\noindent
{\bf Marked ultrametric measure spaces.}
Consider quadruples
\be{ag37c}
(U,d, \mu,\kappa),
\ee
where $(U,d)$ is an ultrametric measure space and $\mu$ is a \emph{finite measure} on $(U, \CB(U))$, with $\CB$ denoting the Borel-$\sigma$ algebra on a topological space. Finally $\kappa$ is a map
\be{ag37d}
\kappa \colon U \la V, \quad
V \mbox{ is a Polish space, $\kappa$ is }
(\CB(U), \CB(V))\mbox{-measurable.}
\ee
Consider isomorphy classes of the quadruplet, i.e., equivalence classes with respect to maps $\varphi$,
\be{ag37e}
\varphi \colon \supp(U) \la \supp(U^\prime),
\ee
satisfying
\be{ag37f}
\varphi \mbox{ is a bijection, \quad $\varphi$ is isometric, } \quad
\mu^\prime = \varphi_\ast \mu, \quad \kappa^\prime (u^\prime) = \kappa(\varphi^{-1} (u^\prime)) \mbox{ for all } u^\prime \in U^\prime.
\ee
The set of all isomorphy classes of quadruplets with $U \subseteq \R^\N$ is denoted
\be{ag37g}
\wt \U^V.
\ee
The space $\wt \U^V$ gives rise to a {\em Polish space} (see \cite{DGP11}), denoted
\be{ag37h}
\U^V, \text{write $\U$ if there are no marks,}
\ee
which is obtained by considering instead of mark \emph{maps} $\kappa$ now \emph{transition probabilities $\kappa$} from $(U, \CB(U))$ to $(V,\CB(V))$, and requiring in \eqref{ag37f} that
\be{ag37i}
\mu^\prime \otimes \kappa^\prime = (\varphi_\ast \mu) \otimes \kappa (\varphi^{-1}(\cdot), \cdot).
\ee

We can associate with an element $(U,d,\mu,\kappa)$ a \emph{sampling measure} on $U \times V$, namely $\nu=\mu \otimes \kappa$.

In our situation, we choose $U$ as the set of \emph{individuals}, $d(\cdot,\cdot)$ as the \emph{genealogical distance} (twice the time back to the most recent common ancestor), $\mu$ as the \emph{sampling measure} on the set of individuals $U$, $V = \K \times \Omega_{N}$ (or $\K \times \Omega^{\infty}_N$), and $\kappa$ as the \emph{mark function} assigning (type, location) to every individual in $U$. 

\medskip\noindent
{\bf Topology.}
Consider the case where $\mu$ is a probability measure. The {\em topology on} $\U^V$ is defined by defining \emph{convergence of sequences} based on the idea that all \emph{sampled finite subspaces} from representatives of the sequence converge in law. We denote by 
\be{nudef}
\nu^{\smallu} \in \mathcal{M}_1\bigl(\R_+^{\binom{\N}{2}} \times V^\N\bigr)
\ee 
the (distance, mark)-matrix distribution of $\CU$, defining with the notation $\smallu =((U,r,\mu),\kappa)$:
\be{AG40b}
\begin{aligned}
& R\colon U \times V \la \R^{\N \choose 2}_+ \times V^\N
\\
& R\big((x_i, v_i)_{1 \leq i})\big)
= \big(r(x_i,x_j)_{1 \leq i < j}, (v_i)_{1 \leq i}\big),
\end{aligned}
\ee
\be{AG40c}
\nu^\smallu = R_\ast (\mu^{\otimes \N}), \text{the \emph{distance-mark matrix distribution}.}
\ee
A sequence
\be{AG40d}
(\smallu_n)_{n \in \N}, \quad \smallu_n \in \U^V \mbox{ {\em converges} to }
\smallu \in \U^V \mbox{ if and only if } \nu^{\smallu_n} \La \nu^\smallu \mbox{ as } n \to \infty
\ee
in the \emph{weak topology of measures} on $(\R^+)^{m \choose 2} \times V^m$ for all $m \in \N$. For general finite $\mu$, require convergence of total masses as well.

An important class of continuous bounded functions on $\U^V$ are the \emph{polynomials}: the algebra of functions that are linear combinations of \emph{monomials} $\Phi$ of degree $n$, given by
\begin{equation}
\label{e444}
\begin{aligned}
\Phi (\mfu) = \Phi^{n,\phi} (\mfu)
&= \int_{(U \times V)^n} \nu \big(\underline{(u,v)}\big) \wt \phi (\underline{u,v})\\ 
&= \int_{\R^{m \choose 2} \times V^n} \nu^\mfu 
\left(\dd \underline{\ur},\dd \uv\right) \phi \big(\underline{\ur},\uv\big),
\qquad \phi \in C_b \big(\R^{m \choose 2} \times V^n\big),\, u \in \N_0.
\end{aligned}
\end{equation}
The function $\Phi$ is interpreted as the sample mean of $\phi$ for an $n$-tuple of points, where the natural distances and marks of an $n$-sample are observed. This set of functions determines $\mfu$, and $\{E[\Phi^{n,\phi}], \phi \in C_b(\R_+^{m \choose n} \times V)\}$ determines a $\U^V$-valued random variable $\mfU$, i.e., $\CL[\mfU]$ uniquely. It suffices to use $\phi=\varphi\,\chi$, with $\varphi \in C_b(\R^{m \choose 2}), \chi \in C^c_b (V^n)$, $c$ indicating compact support.

Since we consider the mark space $\Omega^{\infty}_N \times \K$ with ``infinitely'' many components (this applies even to $\Omega_N$), $V$ is not compact. We therefore modify the concepts above slightly and allow \emph{$\sigma$-finite} sampling measures $\mu$ or $\nu$, following \cite{GSW}. This is done by \emph{localising} the states to the population with geographic marks in {\em bounded subsets} $G_i$ with $G_i \uparrow V$, simply requiring convergence for each localisation to a set $G_i$.

Finally, we need the space of equivalence classes of \emph{$V$-marked, grounded tree-like metric measure spaces} with a finite maximal distance from the ground, respectively the fibred version, denoted by
\begin{equation}
\label{e437}
\M^{\downarrow,\rm tree,V},\, \M_{\rm fib}^{\downarrow,\rm tree,V}.
\end{equation}
The topology is the \emph{Gromov weak topology} (\cite{GPWmp13}), generalising~\eqref{AG40c} and ultrametrics above to metrics. This object allows modelling the population of individuals alive at some time, in which the ultrametric level sets above ground describe the genealogy of individuals at a specific time.

In this space we find our marked ancestral space of all individuals ever alive, where the points at distance $t$ above the ground are the genealogies of the population alive at time $t$ (elements of $U_t$). Here ``grounded'' allows considering \emph{forests} of \emph{tree-like} spaces (also called \emph{grounded} tree-like spaces): instead of a \emph{root} we have a set of \emph{founders}. The metric satisfies the \emph{four-point condition}, and the set of points at distance $t$ from the founders is an ultrametric space. The measure is the occupation measure of $(\mu_t)_{t \in [0,T]}$, respectively $(\nu_t)_{t \geq 0}$ in the marked case. This gives~\eqref{e437} a Polish space. We will later introduce a subset of this space with a stronger topology.

\medskip\noindent
{\bf Outline.}
\label{out}
Section~\ref{s.resgen} formulates the main results of the paper, consisting of 14 theorems. This is followed by preparations for the proof of these theorems: Section~\ref{s.3} develops the main lemmas, Section~\ref{ss.pexisten} proves these lemmas, and constructs and analyses the dual \emph{L\'evy-web} as the main tool of proof. The proofs of the first five theorems, concerning the genealogy process, are in Section~\ref{s.proofsXXX}, the proofs of the last nine theorems, concerning the ancestral web, are in Section~\ref{s.proofs}.


\section{Main theorems on evolving genealogies}
\label{s.resgen}

We proceed with spelling out the \emph{three themes} mentioned in Section~\ref{ss.backgoal}. We state and explain 14 theorems in~\ref{ss.contsgen},~\ref{ss.levyweb} and~\ref{ss.csl}, suppressing technicalities. In Section~\ref{s.3} we formalise and detail certain relevant concepts, and formulate two key propositions on L\'evy-webs that are needed for the proofs of the theorems. We prove these propositions in Section~\ref{ss.pexisten}.

\medskip\noindent
{\bf Theme (1): The genealogical process.}
The first theme is to establish that the genealogical process of our model on continuum space $\Omega_N^\infty$ is well defined as the limit of $\Omega_N^L$-truncations. We show that a unique stochastic process of evolving genealogies $\mfU=(\mfU_t)_{t \geq 0}$ can be constructed via a \emph{martingale problem}, provided we are in the {\em stochastic} regime for the measure-valued model on $\Omega_N^\infty \times \K$ treated in~\cite{gh1}. In the deterministic regime the evolution is deterministic and admits only a \emph{trivial genealogy}: no two individuals have a common ancestor, and we have a population of independently migrating ancestral paths. The process is obtained as the $L \to \infty$ limit of genealogies evolving as Markov processes on \emph{countable discrete} geographic spaces $\Omega_N^L$. We state the path properties of this process, its duality relation, its martingale problem characterisation, and its long-time behaviour.

\medskip\noindent
{\bf Theme (2): The ancestral web.}
The second theme is to construct an \emph{ancestral web} $\mfW=(W_T)_{T \geq 0}$ in which the genealogies of the currently-alive populations (described by $\mfU$) are embedded for all times $t \in [0,T]$, together with \emph{all} individuals ever alive before time $T$, as $T$ evolves. We obtain $\mfW$ by constructing its \emph{dual web} $\wh \mfW$. In our model this dual is the so-called \emph{L\'evy-web}, which represents $\mfW$ at every fixed time $t>0$. Moreover, both $\mfW$ and $\wh \mfW$ admit a \emph{strong graphical representation}: they can be embedded as populations of \emph{forward} ancestral paths (lineages) and \emph{backward} paths of descent, respectively, so that both types of paths are oppositely \emph{oriented geodesics} connecting time-$0$-level and time-$t$-level individuals. Both paths live in \emph{a single} space, the \emph{ancestral space $\mfX$}, which is a \emph{marked grounded fibred metric measure space}. This continuum random graph describes the complete \emph{time-space picture} of the forward \emph{and} the dual backward processes.

\medskip\noindent
{\bf Theme (3): Universality.}
The third theme is to establish \emph{universality} properties by using the \emph{continuum-space rescaling limit} (with rescaling of space, time, mass, and resampling rates) of the $(\Omega_N \times \K)$-marked discrete-space genealogy. We show that we obtain special subclasses of our constructed processes from theme (2). These provide the prototypes for the \emph{universality classes} of large time-space scale behaviour for models on $\Omega_N \times \K$ or on $\Omega_N^\infty \times \K$.


\subsection{Results 1: construction, duality and martingale characterisation of the genealogy process $\mfU$} 
\label{ss.contsgen}

We next state the basic results of the paper, organised into four steps:
\begin{itemize}
\item[(i)] 
Construct a process on the geographic space $\Omega_N^\infty$ via \emph{approximation} by $L$-truncation with models on the discrete geographic space $\Omega_N^L$.
\item[(ii)] Formulate the \emph{duality}.
\item[(iii)] Characterise the process by a \emph{well-posed martingale problem}.
\item[(iv)] Specify the test functions on which we define \emph{operators} for the martingale problem.
\end{itemize}

We have three universality classes for the evolution:
\begin{itemize}
\item 
L\'evy-mass transport combined with a Fleming-Viot-type continuum-space process, based on the semigroup of \emph{point-recurrent} L\'evy-processes.
\item 
Pure Cannings block-resampling, with Harris-recurrent reshuffling mass transport.
\item 
A mixture of the above mechanisms with tuned parameters.
\end{itemize}
For the latter two, in Section~\ref{sss.project} we construct a stochastic process. In Section~\ref{sss.charac} we introduce duality and give the long-time behaviour. In Sections~\ref{sss.dyngen} and~\ref{sss.mart} we give the martingale problem with its ingredients.


\subsubsection{Construction of the genealogy process $\mfU$ as a truncation limit}
\label{sss.project}

We construct the $\U^V$-valued process in the stochastic regime as the limit of processes on countable geographic spaces approximating $\Omega_N^\infty$, namely $\Omega_N^L$ (the $L$-refinements of $\Omega_N$). We use \emph{truncated} versions: first truncate in elements of $\Omega_N^\infty$ the levels less than $-L$ for some $L \in \N$, and then let $L \to \infty$. The truncated models on $\Omega_N^L$ are obtained by truncating the parameter sequence at $-L$. This yields a countable geographic space and processes with a well-known duality. We use this \emph{duality} to verify the approximation.

\medskip\noindent
{\bf Defining $\mfU$ via truncation.}
We begin by giving a detailed construction of the finite-dimensional distributions of the processes. We \emph{define a process $\mfU$}, the continuum-space Cannings process with block-resampling, via a \emph{truncation limit}: we turn off migration, reshuffling, and block-resampling in small balls beyond level $-L$, and shift the site resampling to level $-L$. At that level we take Fleming-Viot resampling with rate $d\,N^L$, which keeps enough volatility in $1$-balls. We subsequently give the construction via duality, following the strategy in~\cite{gh1}. These processes correspond to the continuum-space $\Omega_N^\infty$-version of our parameter sequences~\eqref{e302}.

More precisely, we construct processes $\{\mfU^L\colon\, L \in \N\}$ on $\Omega^{\infty}_N$, which are $\U^{\Omega_N^L}$-valued systems embedded in $\Omega_N^\infty$ as geographic space, and define $\mfU$ as the $L \to \infty$ limit (see Theorem~\ref{T.wptvcont} below). The parameters are chosen by truncating $(c_k)_{k \in \Z},(\lambda_k)_{k \in \Z},(\Lambda_k)_{k\in \Z}$ for $k<-L$, and using single-site Fleming-Viot resampling at rates $dN^L$, respectively, Cannings resampling at rates $N^L \Lambda_{-\infty}$.

The basic idea is to use the duality for genealogies of systems for \emph{finite $L$}, which are Cannings processes with block-resampling on a \emph{countable} Abelian group. Genealogies for such processes were constructed in~\cite{gklimw} in the case \emph{without} block-resampling. We show convergence of these processes to obtain a $\U^{\Omega_N^\infty \times \K}$-valued stochastic process $\mfU^\infty$, later characterised as the solution of a martingale problem. We argue that block-resampling can be added and that the results in~\cite{gklimw} on existence and uniqueness of the martingale problem carry over. On $\Omega_N$ the new mechanism for finite $L$ generates a Markov jump process and poses no principal problems. Details are provided in the proof section (the measure-valued case is in~\cite{GHKK14}).

\medskip\noindent
{\bf Stochastic vs.\ deterministic regimes.}
We next give the conditions under which a dynamics on $\U_1^{\Omega^{\infty}_N \times \K}$ can be constructed uniquely. From~\cite{gh1} we know for which choice of parameters the projection of $\mu_t$, $t \geq 0$ onto $\K$ is \emph{stochastic}. We show that in that case also $\mfU_t$, $t \geq 0$, \emph{exists} as a $\U^V$-valued process and is truly stochastic. Assuming block-resampling, we distinguish between \emph{regimes $A$ and $B$} depending on whether
\be{ag55a}
\sum^\infty_{k=0} \frac{1}{c_k} \sum_{\ell=0}^k \lambda_\ell \,\left\{\begin{array}{ll}
= \infty, &\text{regime A},\\
<\infty, &\text{regime B }.
\end{array}
\right.
\ee

\begin{theorem}[Continuum-space process $\mfU$ as $L$-truncation limit]
\label{T.wptvcont}
$\mbox{}$\\
Consider initial states $\CU \in \U^{\Omega^{\infty}_N \times \K}_{\mathrm{ reg }} \subseteq \U^{\Omega^{\infty}_N \times \K}$, where the resampling measure projected on $\Omega^{\infty}_N$ is $\mfh$.
\begin{itemize}
\item[(a)] In the stochastic regime $A$, characterised by
\be{ag6}
\sum_{k=0}^\infty \frac{1}{c_k} \Big( \sum_{\ell=0}^k \lambda_\ell\Big) = \infty,
\ee
the following hold for $\mfU_0 \in \U^{\Omega^{\infty}_N}_{\mathrm{reg}}$:
\begin{itemize}
\item[(i)] 
The process $\mfU^\infty$ is uniquely defined by
\begin{equation}
\label{e604}
\CL \left[\left(\mfU^\infty_t\right)_{t \geq 0}\right] = \lim_{L \to \infty} \CL \left[\left(\mfU_t^{L}\right)_{t \geq 0}\right], 
\end{equation}
in particular, the limit exists in the state space $\U^{\Omega_N^\infty \times \K}$.
\item[(ii)] 
The process $\mfU^\infty$ is a strong Markov and Feller process on $\U^{\Omega_N^\infty \times \K}_{\rm reg}$. In particular, $proj_{\Omega_N^\infty}(\nu)$ $= \mfh$, where $\nu$ is the sampling measure of $\mfU$, and the process has c\`adl\`ag paths, where jumps occur only at times of block-resampling events (in monotype balls no jump is induced).
\item[(iii)] 
In the stochastic regime, $\mfU^\infty$ has for times $t>0$ populations with a random depth-$h$ family decomposition for each $h \in (0,t]$, i.e., a decomposition of $U^\infty_t$ into disjoint open and closed $2h$-balls in $U_t$. In particular, non-trivial evaluation functionals of distances have positive variances. Each $2t$-family covers a random number (possibly countably many) of balls in $\Omega_N^\infty$ with random position and random radius. The decomposition induces a disjoint decomposition of $\Omega_N^\infty$ into disjoint unions of random balls, each covered (up to $\mfh$-null sets) by one of the $2t$-subfamilies, each of one type.
\end{itemize}
\item[(b)]
\begin{itemize} 
\item[(i)] 
In the deterministic regime B, no limiting $\U^{\Omega_N^\infty \times \K}$-valued genealogy process exists.
\item[(ii)] 
In the process on $\Omega_N^L$, the probability for two individuals to have a common ancestor tends to zero as $L \to \infty$, for all times $t>0$. Hence there is dust formation as $L \to \infty$ for the decomposition into depth-$t$ subfamilies. In formulas, for all $\xi  \in \Omega^{\infty}_N$,
\be{ag56}
\lim_{L \to \infty} P\left\{ \int_{U_t \times U_t} \mu_t^L \left(\dd \iota \right) \mu_t^L \left(\dd \iota^\prime \right)
1 \left(r_t^L(\iota, \iota^\prime) > 2t \right) \right\} \geq 1-\wh \delta \quad \forall\, \wh \delta >0,t >0,
\ee
and the population consists of individuals without any common ancestor for any pair of individuals. Hence no limit exists in $\U^{\Omega_N^\infty \times \K}$, and the population is fragmented into a dust of single-ancestor families performing effective L\'evy-processes as $L \to \infty$.
\item[(iii)] 
The measure-valued dynamics of the marks $(\wh X^L_t)_{t \geq 0}$ in the limit as $L \to \infty$ is deterministic and equal to $\mathfrak{h} \otimes \theta$ for $t>0$~\cite{gh1}.
\end{itemize}
\end{itemize}
\end{theorem}

In the sequel we drop the superscript $\infty$ on $\mfU$. We will see later that the limit process can alternatively be characterised as follows, making stochastic analysis applicable.

\begin{corollary}[Characterization $\mfU$]
The law of the process $\mfU$ from~\eqref{e604} is the unique solution of the $(L_\U^{\Omega_N^\infty}, \Pi,\mfU_0)$ martingale problem from Theorem~\ref{T.mp}. Its finite-dimensional distributions are given by~\eqref{e642} or also~\eqref{e526} in Theorem~\ref{th.588}: part (a) by duality, and part (b) by strong duality, where the state is represented by an entrance law of the genealogy process of a spatial coalescent with block-resampling from~\eqref{e302} and below.
\end{corollary}

\medskip\noindent
{\bf Excursion on hotspots.}
The \emph{hotspots} are points in geographic space where volatility develops during the evolution at time $t$. This is a very small set in $\Omega_N^\infty$. We identify the basic set for the state of the process $\mfU$ and the set of points at time $t$ that create \emph{volatility} in the subsequent evolution from time $t$ onwards. Fix a time $t$ and consider the locations in the geographic space $\Omega_N^\infty$ of individuals in $2h$-balls in the \emph{ancestral space} that survive at time $t+h$ with $h>0$. For that purpose, consider the points
\be{ag70}
E^t_h \subseteq \Omega^{\infty}_N, \quad \iota \in E^t_h \mbox{ is the time-$t$ ancestor of an individual at time } t+h,
\ee
which form a countable set. The set of all \emph{founders} $E^t$ at time $t$, i.e., the set of ancestors at time $t$ of the population at time $t+h$ for some $h>0$, is
\be{ag71}
E^t= \bigcup\limits_{h>0} E^t_h.
\ee
Note that (with $\mathfrak{h}$ the Haar measure on $\Omega^{\infty}_N$)
\be{ag73}
\mfh (E^t) =0.
\ee
The set $E^t$ is a \emph{countable} $\mfh$-null set, but it is the \emph{complete source} of \emph{all} subsequent \emph{stochastic} change. Namely, the points of $E^t$ are exactly those that subsequently create surviving populations, thereby generating the \emph{volatility}.

The possible states at time $t$ are systems of balls in $\Omega_N^\infty$ completely covered by descendants of a founder at an earlier time. In the backward view, points in $E_h^t$ are where a partition element covering some ball at time $t+h$ passes through at time $t$.

We can decompose the descendants of a single founder at time $0$ at time $t+h$ into depth-$h$ subfamilies (with a founder at time $t$), arising from the subfamily splitting of the ancestral line to the founder at time $0$ at times $s \in [0,t]$. The founder $\iota$ at time $0$ populates at time $t$ the $K^\ast_t(\iota)$-ball in $\Omega_N^\infty$ around the position of $\iota$ completely. At time $t+h$ it covers completely the set
\be{ag72}
\bigcup\limits_{\iota^\prime \in A_{t,t+h} (\iota)} B_{K^\ast_{t+h}(\iota^\prime)}(\iota^\prime),
\ee
where $A_{s,t}(\iota)$ are the ancestors of $\iota$ from time $t$ at times $s \in [0,t)$.


\subsubsection{Duality on $\Omega_N^L$ and $\Omega_N^\infty$}
\label{sss.charac}

To prove our results, we make heavy use of the \emph{duality relation} and the fact that it is easy to show convergence of our truncated processes as $L \to \infty$, by showing \emph{convergence of their dual processes} to a limit process.

\medskip\noindent
{\bf The dual process.}
We first introduce the \emph{dual processes} and the \emph{duality functions}. The
\begin{equation}\label{e508}
\text{\emph{dual process $\mfY$} will be a $\U^{\Omega_N^\infty \times \K}$-valued process, }
\end{equation}
driven by a \emph{spatial delayed Kingman-coalescent}, $\Lambda$-coalescence, and (for finite $L$) $(\Lambda_k)_{k \in \Z}$-block coalescence with reshuffling. See the end of Section~\ref{s.basic}, in particular,~\eqref{e302}. An ultrametric is extracted from the partition-valued process at backward time $t$ as twice the time for two elements to coalesce, or $2t$ if they do not coalesce. This turns a finite coalescent and its elements into a finite ultrametric space. By assigning the time-$t$ location and type as marks, we obtain a finite element of $\U^V$ for the duality.

The duality function $H$ subtle. In particular,
\begin{equation}
\label{e512}
\text{$H$ is built from specific \emph{polynomials} on $\U^{\Omega_N^\infty \times \K}$, namely, \emph{monomials}.}
\end{equation}
These polynomials are introduced in Section~\ref{sss.mart}, see~\eqref{pol1} and~\eqref{s17}. The idea, following~\cite{gh1}, is to build on these polynomials by using multiple averaging with the Haar measure $\mfh^{\otimes n}$ over the locations of the finite coalescents, specifying the truncated dual functions associated with the $L$-truncated process on $\Omega_N^L$. This reflects that the state arises as the $L$-truncation limit of the $\Omega_N^L$-models.

\begin{theorem}[Duality and strong duality]
\label{th.588}
In regime A (recall \eqref{ag6}), consider the process $\mfU$ starting in $\mfU_0 \in \U_{\rm reg}^{\Omega_N^\infty \times \K}$.
\begin{itemize}
\item[(a)] 
There exists a duality function $H\colon\,\U^V \times \U^V \to \R$ and a dual process $\mfY$ taking values in $\U^V$ with $V=\Omega_N^\infty \times \K$ such that, for the dual evolving backward from every time horizon $t>0$ starting in a subset of $\U^V$ (consisting of those elements with finitely many elements, i.e., dual individuals),
\begin{equation}
\label{e642}
E_{\mfU_0} \left[H\left(\mfU_t,\mfY_0\right)\right] =E_{\mfY_0}\left[H \left(\mfU_0,\mfY_t\right)\right].
\end{equation}
An example is given in~\eqref{e729}.
\item[(b)] 
For every time horizon $t>0$ there exists an entrance law $\wh \mfY$ for the above dual process $\mfY$, starting at backward time $0$ in $\wh \mfY_0$ with countably many dual particles (one at each site of $\Omega_N^{\infty,\rm disc}$), which is a specific countable dense subgroup (see~\eqref{e1141}). Equivalently, this amounts to starting the dual $\mfY$ effectively with ``sampling'' measure $\mfh$ on $\Omega_N^\infty$.
\end{itemize}
For $\mfU_0$ being a single root, a strong dual representation holds:
\begin{equation}
\label{e526}
\CL_{\mfU_0}[\mfU_t]=\CL_{\mfY_0}[\mfY_t].
\end{equation}
The general case for $\mfU_0$ is obtained through the concatenation of $\wh \mfY_t$ with $\mfU_0$ as
\begin{equation}
\label{e592}
\mfU_t=\mfU_0 \vdash \mfY_t.
\end{equation}
Here, $\mfU_1 \vdash \mfU_2$ means glueing $\mfU_2$ to $\mfU_1$ as follows: the $2t$-balls are glued to individuals sampled from $U_1$ i.i.d.\ infinitely often according to $\mu_1$ (the ancestors of the founding fathers; see \eqref{ag70}--\eqref{ag72}), and the distance in the glued object is $r_1+r_2$ on the glued ultrametric space.
\end{theorem}

\begin{remark}[Duality in regime B]
\label{r.caseb}
{\rm In regime B (deterministic), $\mfY$ becomes a system of independent L\'evy-processes, one starting from each point.}
\end{remark}

The entrance law can also be viewed as starting with one particle in each site of $\Omega_N^\infty$, arising as the limit of starting with one individual in every point of $\Omega_N^{\infty,\rm disc}$ and with the $L$-truncated system embedded in $\Omega_N^\infty$ as entrance law in the $L \to \infty$ limit. This is carried out in~\cite{gh1}.

\medskip\noindent
{\bf Polynomials and monomials.}
To obtain the duality function we need the following objects. We introduce on $\U^V$ an \emph{algebra of polynomials} as test functions and build the duality functions on monomials characterised by parameters that identify an initial state of a finite spatial coalescent. A monomial $\Phi$ on $\U^V$ of order $n \in \N$ is an average of the evaluation of an $n$-sample from the set $U$. See~\eqref{pol1} for details.

Duality is not only key for the construction of the process, it also allows us to obtain the long-time behaviour.

\begin{theorem}[Longtime behaviour]
\label{T.ltb}
In regime A the following hold. 
\begin{itemize}
\item[(a)] 
Consider the process $\mfU=(\mfU(t))_{t \geq 0}$ started in a translation-invariant ergodic initial state. Then
\be{ag55}
\CL [\mfU (t)] \tto \CL [\mfU(\infty)],
\ee
where $\CU(\infty)$ is a random variable with values in $\U^V$ having a translation-invariant law, and the corresponding measure-valued state is the unique shift-ergodic equilibrium state with the mean type intensity $\theta$ of $X(0)$. 
\item[(b)] 
If $\mfU$ (the $\U^{\Omega_N^\infty}$-valued process) is projected onto the $\U$-valued component, then for all initial states a unique random genealogy $\mfU_{\rm gen}(\infty)$ is obtained. This random variable is the infinite-time limit of the genealogy of the dual population of the $\U^{\Omega_N^\infty}$-valued $((c_k)_{k \in \Z}, d^\ast,(\Lambda_k)_{k \in \Z})$-coalescent started with one particle at each site of $\Omega_N^{\infty,\rm disc}$.
\item[(c)] 
If $\CU$ in equilibrium is projected onto the mark component in $V$ and the projection of the sampling measure onto $\K$, respectively, $\Omega_N^\infty$ is taken, then the outcome is
\begin{equation}
\label{e513}
\int \delta_{\underline{k}} \theta (\dd k), \quad \theta \in \CM_1(\K), \text{ respectively, $\mfh$,}
\end{equation}
where $\underline{k}$ is the mark function for the type constant equal to $k$.
\end{itemize}
\end{theorem}

There is a related model, obtained by letting $d \to \infty$, which yields as dual the process with instantaneous coalescence upon hitting.

\begin{corollary}[Infinite volatility limit]
\label{cor.volat}
For the process $\\mfU^{[d]}$, $d>0$, in regime A,
\begin{equation}
\label{e676}
\CL \left[\mfU^{[d]}\right] \dto  \CL \left[\mfU^\sim\right],
\end{equation}
where $\mfU^{[\infty]}$ satisfies Theorems~\ref{T.wptvcont}--\ref{T.ltb}.
\end{corollary}

\noindent
This process also arises when studying one-ancestor clusters and their expansion in time and geographic space. However, the analysis of this process poses serious challenges and the subsequent results do not apply to it.


\subsubsection{Characterisation of the genealogy process $\mfU$ via a martingale problem}
\label{sss.dyngen}

Next, we generalise the classical $\U^{G \times \K}$-valued Cannings dynamics for \emph{countable} Abelian groups studied in \cite{gklimw} to \emph{include block-resampling}, and use this to pass from $\Omega_N^L$ to $\Omega_N^\infty$. Even {\em without} types, the dynamics is given by a \emph{well-posed martingale problem} on $\U$, but requires quite a bit of preparation.

We need to introduce \emph{test functions}: certain polynomials from $\Pi$ appearing in our martingale problem. We have to define the operators expressing the evolution mechanisms, and arrive at an operator $\mfL^{\Omega_N^\infty}_{\rm gen}$ defined on these test functions. We first state the key results, and afterwards specify the ingredients in detail. In a further theorem we state the form of these objects.

The characterisation of the genealogical process via a martingale problem reveals the following scenario when passing from $\Omega_N$ to the spatial continuum $\Omega_N^\infty$. Recall the distinction of~\eqref{ag55a}.

\begin{theorem}[Martingale problem characterisation: genealogical Cannings process]
\label{T.mp}
$\mbox{}$\\
In regime $A$, the $(\Pi, \mfL^{\Omega^{\infty}_N}_{\mathrm{gen}}, \delta_{\CU_0})$-martingale problem is well-posed for every $\mfU_0 \in \U^V$. The ingredients $\Pi$ and $\mfL^{\Omega_N^\infty}_{\rm gen}$ are explicitly identified in Theorem~\ref{th.sss.mart} below and in the sequel. The resulting solution, the law of the process $(\mfU_t)_{t \geq 0}$, is in law exactly the genealogical Cannings process with block-resampling $\mfU$ on $\Omega_N^\infty$ from Theorem~\ref{T.wptvcont}.
\end{theorem}

In regime B, there exists \emph{no solution} for the dynamics of $\mfU$ in the given state space $\U^V$ (see~\eqref{ag56} and the sequel). The measure-valued projection of such a process on the mark space $V$ would have to be \emph{deterministic} and equal $\mathfrak{h} \otimes \theta$ for $t>0$. The genealogy process forms dust as $L \to \infty$, consisting only of single-ancestor lines migrating in space.

In the transient migration/reshuffling case, the system stabilises in states that fragment locally into families that $\mfh$-almost all have different founding fathers. This is referred to as \emph{coexistence} or \emph{diversity}. Since in the continuum limit a geographic ball represents many sites in the original discrete model, the local population has a genealogy that \emph{fragments into dust}.

In all cases we obtain the {\em measure-valued Cannings process} $\wh X^\infty$ with values in $\CM(\Omega_N^\infty \times \K)$ as the {\em projection} of the {\em sampling measure} in $\mathfrak{U}$, or as a process of historical paths. However, $\wh X^\infty$ does not correspond to a non-degenerate genealogy with locations and types, i.e., with $V=\Omega_N^\infty \times \K$.


\subsubsection{Martingale problem ingredients: identification of operators and test functions}
\label{sss.mart}

The process $(\mfU_t)_{t \geq 0}$ solves the well-posed $(\mfL_\U^{\Omega_N^\infty}, \Pi, U_0)$-martingale problem. We now identify the ingredients of the operator and the domain of test functions.

\begin{theorem}[Martingale problem characterisation of $\mfU^\infty$: Identification of the ingredients]
\label{th.sss.mart}
$\mbox{}$\\
The operator of the martingale problem is the sum of four operators describing the action of migration/reshuffling, aging via metric growth, resampling acting on distances, and change of type frequencies via resampling:
\be{pol2}
\mfL^{\Omega^{\infty}_N}_\U = \mfA^{\mathrm{mig}} +\mfA^{\rm grow} +  \mfA^{\rm res} + \mfA^{\rm res, \K}
\ee
(see~\eqref{pol0}--\eqref{e:generator}). The domain of $\mfL^{\Omega_N^\infty}_{\U^V}$ are the polynomials $\Pi$ (see below).
\end{theorem}

\begin{definition}[Polynomials]
\label{def.poly}
Define $\Pi$ as the subalgebra of $C_b(\U^V, \R)$ spanned by the monomials $\Phi$,
\begin{equation}
\label{pol1}
\begin{aligned}
\Phi(\mfu) &= \int_{(U \times V)^n} \nu^\mfu \left(\dd({\underline{u,v}}\right) \phi \left((\underline{u,v})\right)
= \int_{\R^m \times V} \nu^\smallu (\dd\dr, \dd \underline{y})\, \varphi(\dr) \chi(\underline{y}),\\ 
&\quad \smallu \in \U^V, \,m= {n \choose 2}, \,n \in \N,
\end{aligned}
\end{equation}
for functions $\chi \in C_b^c ((\Omega^{\infty}_N \times \K)^n, \R)$ and $\varphi \in C_b((\R^+)^{n \choose 2},\R)$ for some $n \in \N$. Take as domain all such $\Phi$ with $\phi=\varphi\,\chi$ and $\chi$ of the form
\be{ag60}
\chi(\underline{y}) = g(\underline{\xi}) f(\underline{v}), \quad \mbox{ with } \underline{y} = (\underline{\xi},\underline{v}) \in (\Omega^{\infty}_N \times \K)^n,\,g \in C^c_b \left(\left(\Omega_N^\infty\right)^n\right),\, f \in C_b \left(\K^n\right). 
\ee
\end{definition}

We next build up the objects, functions, and operators in the above theorem in complete detail. We specify the ingredients of the martingale problem in \emph{two steps}. Since the projection of the process $\mfU$ onto the genealogy component is again a Markov process, as is the evolution of the projection on the marks, we can first study these dynamics separately. The latter is already treated in~\cite{gh1}, so we first look at the case \emph{without marks} and then extend to the marked case.

\medskip\noindent
{\bf Step 1: the pure genealogy evolution.}
$\mbox{}$

\medskip\noindent
\textbf{(i) (a) Test functions.}
Let $\Phi \colon \mathbb{U}_1 \to \R$ be a polynomial of degree~$n \in \N$, i.e., $n$ is the minimal number such that there exists a test function $\varphi$ depending on $(r_{i,j})_{1 \leq i < j \leq n}$ only, such that (recall~\eqref{AG40c} for the distance matrix distribution $\nu^{\mathcal{U}}$)
\be{eq:polynomial}
\begin{aligned}
\Phi(\smallu)=\int_{(U \times V)^n} \nu \left(\dd(\underline{u,v})\right)=
\int_{\R^m} \nu^{\smallu}(\dd\uuline{r})\,\varphi\big(\uuline{r}\big),
\quad
\smallu \in \mathbb{U}_1, \; \Pi^{\rm gen} \text{   is the generated algebra},
\end{aligned}
\ee
which we call the \emph{polynomials}, denoted by $\Pi^{\rm gen}$ to distinguish from the \emph{general case $\Pi$} with marks. We generalise to $\sigma$-finite sampling measures $\nu$ on $U \times (\Omega_N^\infty \times \K)$ in $\Omega_N^\infty$, finite on compact sets (following~\cite{GSW}) by working with test functions that have compact support on $\Omega_N^\infty$.

\medskip\noindent
\textbf{(b) Duality function 1.}
To get a \emph{duality function} $H= \U^V \times \U^V \to \R$ we set, for $\mfu=(U \times V,r,\nu), \mathfrak{y}=(\Pi,r^\Pi,\mathfrak{n})$ with $\Pi$ the \emph{partition},
\begin{equation}
\label{e729}
H(\mfu,\mathfrak{y}) = \int_{\R^m} \nu^\mfu(\dd\uur)\, \varphi(\uur^\Pi), \qquad \uur^\Pi := \uur \left(\Pi(1),\ldots,\Pi(n)\right).
\end{equation}

\medskip\noindent
\textbf{(ii) Operators.}
First, we define the generator of the classical \emph{Fleming-Viot resampling}, i.e., the generator of the \emph{$\U_1$-valued Fleming-Viot process}, as
\be{pol0}
\mfA^{\mathrm{res}} \Phi (\smallu) = d \sum^\infty_{i,j=1}
\int_{\R^m} \nu^\smallu (\dd \dr) (\varphi(R_{\{i,j\}} (\dr)) - \varphi(\dr)),
\ee
where $R_{\{i,j\}} (\dr)$ sets $r_{i,j}$ equal to zero and $r_{l,j}$ equal to $r_{l,i}$ for $l \in \N$.

\medskip\noindent
{\bf Cannings resampling.}
Next, we define further ingredients of the generator to include the \emph{Cannings resampling}. We need random $\{0,1\}$-sequences $w=(w_k)_{k \in \N}$ allowing us to mark the individuals killed and replaced by the next generation. For $x\in[0;1]$, define the measure $p^x$ on $\{0,1\}^\N$ as
\be{AG21b}
\big(x \delta_1 +(1-x)\delta_0\big)^{\otimes \N}.
\ee

Given $I\subseteq\mathbb{N}$ and $\omega\in\{0,1\}^\N$, the \emph{replacement operator} $R_I \colon \R_+^{\N\choose 2}\to \R_+^{\N\choose 2}$ sends a distance matrix of a sequence of points to the distance matrix of a new sequence obtained by identifying all points with an index $k\in I$:
\be{e:replacement}
\begin{aligned}
\big(R_{I}(\uuline{r})\big)_{i,j}
&:=
\begin{cases}
r_{i,j}, & i,j\not\in I,\\
r_{\min I,j}, & i\in I,j\not\in I,\\
r_{\min I,i}, & j\in I,i\not\in I,\\ 
0, & i,j\in I.
\end{cases}
\end{aligned}
\end{equation}

Write for $\omega\in\{0,1\}^{\mathbb{N}}$, $I_\omega:=\{i\in\mathbb{N}\colon\omega_i=1\}$. The generator for the \emph{unmarked} process is the linear operator on ${\mathcal C}_b(\U)$, denoted $\mfA_{\rm gen} \colon \Pi^{\rm gen} \to C_\text{b}(\U))$, given by
\begin{equation}
\label{e:generator}
\begin{aligned}
&\mfA_{\rm gen} \Phi(\smallu)\\
&:= \big\langle\nu^\smallu, (\nabla\varphi,2)\big\rangle
+
\int_{\R^m} \nu^{\smallu}(\mathrm{d}\uuline{r})
\int_0^1\frac{\Lambda(\mathrm{d}x)}{x^2}
\int_{\{0,1\}^\N}p^{(x)}(\mathrm{d}\omega)
\Big(\varphi\big(R_{I_\omega}\uuline{r}\big)-\phi\big(\uuline{r}\big)\Big) + \mfA^{\mathrm{res}} \Phi(\smallu),\\
&= (\mfA^{\mathrm{growth}} + \mfA^{\mathrm{Cann}} +\mfA^{\mathrm{res}}) \Phi (\smallu),
\end{aligned}
\end{equation}
where $(\boldsymbol{\cdot},\boldsymbol{\cdot})$ is the inner product on $\R^{\N \choose 2}$. The first term describes the \emph{growth} of genealogical distances due to {\em aging}. The second term describes the transition upon {\em change of generation according to Cannings resampling}. The third term describes the change of generation via classical \emph{Fleming-Viot} {\em resampling}. The latter absorbs, in the continuum-space approximation by discrete spaces, the single-site Cannings resampling, as proved in~\cite{gh1}.

\medskip\noindent
{\bf Step 2: extension to marked case.}
We now turn to the \emph{marked} case with state space $\U^V$. We first incorporate \emph{types} and then consider the \emph{locations} with the migration dynamics of classical migration and migration by reshuffling at resampling times.

\medskip\noindent
\textbf{(i) (a) Test functions.}
We use the set $\Pi$ from Theorem~\ref{th.sss.mart}, which means complementing in~\eqref{pol1}
\begin{equation}
\label{s17}
\phi \text{ by the further factor } \chi =g\, f \text{ with } g \text{ for locations and } f \text{ for types}
\end{equation}
(recall~\eqref{pol1}, \eqref{ag60}). We have to provide formulas for $\mfA^{\rm mig}$ and $\mfA^{\rm Cann,b}_{(L,\xi)}$.

\medskip\noindent
\textbf{(b) Duality function 2.}
The \emph{duality function} $H:\U^V \times \U^V \to \R$ now includes marks. For $\mfu=(U \times V,r,\nu),\mathfrak{y}=(\Pi,\uur^\Pi,\underline{y},\mathfrak{n})$,
\begin{align}
\label{e792}
H\left(\mfu,\mathfrak{y}\right) &= \int_{\R^m \times V} 
\nu^\mfu \left(\dd\uur,\dd\underline{y}\right) \varphi \left(\uur^\Pi\right)\chi(\underline{y}^\Pi),\\
\uur^\Pi &:=\uur(\Pi(1),\ldots,\Pi(n)), \nonumber\\
\underline{y}^\Pi &:=\left(y \left(\Pi(1)\right), \ldots, y \left(\Pi(n)\right)\right). \nonumber
\end{align}

\medskip\noindent
\textbf{(ii) Operators.}
We verify that the formal operators are well defined on the domain. This was already largely handled in~\cite{gh1}. The linear operator from~\eqref{e:generator} must be \emph{extended} from $\U$ to $\U^V$ with $V=\Omega^{\infty}_N \times \K$, i.e., the operators extend from the domain $\Pi^{\rm gen}$ to $\Pi$. The extended operators act on the \emph{test function $\varphi$ only} and leave the factor $g$ of $\chi$ untouched: they affect only the factor $f$ of $\chi$ by resampling, i.e., by merging variables.

We specify the operator for the \emph{dynamics of the marks} from $V$, giving the various \emph{resampling} operators acting on these marks, followed by migration and reshuffling.

\medskip\noindent
\textbf{(a) Resampling operators.}

\medskip\noindent
\textbf{(i)} Consider first the diffusion operator: the on-site action. The on-site Cannings mechanism on $\Omega_N^L$ becomes a further Fleming-Viot term on $\Omega^{\infty}_N$ with rate $|\Lambda_{-\infty}|$. This Fleming-Viot term creates a complication when we consider the action on \emph{types} rather than genealogies: the compensator is based on a \emph{time-space functional of the paths} up to the present, rather than a function of the current state. It is the \emph{collision time functional} of the different types. This non-standard situation was treated and resolved in~\cite{gh1}. For genealogy processes, the current state already encodes the history of two individuals through their genealogical distance, so this problem does not arise again here. We obtain for the genealogy part the simple Fleming-Viot operator in~\eqref{pol0}, acting only on $f$ and leaving $g$ and $\varphi$ unchanged. The term on $g$ is the classical Fleming-Viot operator, with the \emph{different-type-collision time-functional}. For formulas on a countable geographic space we refer to ~\cite[Section 2]{gklimw}, and for formulas on $\Omega_N^\infty$ to ~\cite[Section 2]{gh1}.

\medskip\noindent
\textbf{(ii)} The simpler case is the \emph{jump operator} for the Cannings \emph{block-resampling part}. This is a jump transition driven by jumps at \emph{random times}, applying in the infinite-population limit a deterministic map to the present state of the population in a specific ball. The locally finite jump rate gives the relevant operator, called $\mfA^{\rm Cann,b}$. A \emph{non-standard} situation arises because \emph{countably} many operators are added. This was resolved in~\cite{gh1} by using instantaneous entrance into locally monotype states, and hence by locally finite reshuffling rates. We write the \emph{block-resampling} operator as
\begin{equation}\label{e804}
\mfA^{\rm Cann,b} = \sum_{L\in\Z,\,\xi\in\Omega_N^L} \lambda_L \mfA^{\rm Cann,b}_{(L,\xi)}.
\end{equation}
Here, $\mfA_{(L,\xi)}^{\rm Cann,b}$ is the block-resampling (without reshuffling) in the $L$-ball around $\xi$. Writing $1_A$ for the restriction operator of a measure or function to points with $\Omega_N^\infty$-component in $A \in \Omega_N^\infty$, we have
\begin{equation}
\label{e808}
\begin{aligned}
&\mfA^{\rm Cann,b}_{(L,\xi)} (\Phi)= \Phi^\ast,\\  
&\Phi^\ast(\nu)=\Phi(\nu^\ast),\\
&\nu^\ast = \int_{U \times B_L(\xi) \times \K} \nu(\dd u^\prime,\dd \xi^\prime,\dd a^\prime)
\,\left[1_{\mathcal{C} B_L(\xi)}\; \nu + (1-\alpha) \nu 1_{B_L(\xi)}
+ \alpha \left(1_{B_L (\xi)}(\delta_{(u^\prime,a^\prime,\xi^\prime)}-\nu\right) \right].
\end{aligned}
\end{equation}
Here, $\alpha$ depends on $L$ through $\Lambda_L$:
\begin{equation}
\label{e821}
\alpha= \int_0^1 x \Lambda_L (\dd x).
\end{equation}

\medskip\noindent
\textbf{(b) Migration and reshuffling.}
We continue with the \emph{migration} and reshuffling operators. The \emph{migration} operator $\mfA^{\mathrm{mig}}$ is added next: its action on the factor $g$ of $\chi$ must be specified, while it acts as the identity on $f$. The migration corresponds to a \emph{deterministic linear mass flow} component: a \emph{first-order operator} describing a mass flow of distances of pairs from different sites, drawn according to the sampling measure $\nu$ of $\mfU$. For monomials $\Phi$ with $\phi=\varphi\,\chi, \chi= fg$,
\begin{equation}\label{e812}
\mfA^{\rm mig} \Phi^{n,\phi}=\Phi^{n,(\varphi, \chi^\ast)}, \text{ with } \chi^\ast = f \, g^\ast,
\end{equation}
where $g^\ast$ is given as
\begin{equation}
\label{e816}
g^\ast=\mfA^{\rm mig,n}(g)
\end{equation}
with
\begin{equation}
\label{e819}
\mfA^{\rm mig,n}(g)=\sum_{i=1}^n \int_{(\Omega_N^\infty)^n} \left(a(\iota,\iota_i)-\delta_{\iota,\iota_i}\right) \mfh^{\otimes n}(\dd\underline{\iota})g \left(\underline{\iota}^{\iota,\iota_i}\right),
\end{equation}
where $\underline{\iota}^{\iota,\iota_i}:=(\iota_1,\iota_2,\ldots,\iota_{i-1},\iota,\iota_{i+1},\ldots,\iota_n)$ and $a(\cdot,\cdot)$ is the transition density of \emph{L\'evy-migration} on $\Omega_N^\infty \times \Omega_N^\infty$.

\medskip\noindent
{\bf Reshuffling operator.}
In addition to $\mfA^{\rm Cann}$ acting on the genealogy, we need an operator describing the \emph{simultaneous resampling-reshuffling} action on the \emph{marks}. The resampling and the reshuffling appear \emph{simultaneously}, with the reshuffling as an independent increment component on $\Omega_N^\infty$. This additional reshuffling operator acts on the \emph{locations} and hence only on $g$, independently of $u$ and $a$, as
\begin{equation}\label{e824}
\mfA^{\rm Cann,b,\Omega_N^\infty}.
\end{equation}
It describes the deterministic map replacing the state of the location in $B_L(\xi)$: it removes fraction $\alpha$ (from~\eqref{e808}) and adds $\alpha \mfh$ on the ball. This term is part of~\eqref{e808}, \emph{activated at the time of the resampling transition} in the ball.

\medskip\noindent
{\bf Coupling vs.\ independence.}
We note that block-resampling also includes simultaneous \emph{reshuffling} from relocations of individuals in the ball. This makes the structure of the population different from models with Fleming-Viot resampling and migration. The choice of actual positions is \emph{independent} of type and genealogy. For functions of singleton samples, and even for functions of pairs of lineages, this difference is not strongly visible. It is \emph{very significant}, however, for $m$-tuples with $m \geq 3$: in the reshuffling dynamics, different locations are \emph{coupled} L\'evy-processes. Nevertheless, in \emph{large} populations with rapid jumps, the behaviour in the limit is again qualitatively \emph{determined by the pairs}.


\subsection{Results 2: ancestral web $\mfW$, dual L\'evy-web $\wh\mfW$, joint graphical embedding in ancestral space $\mfX^\dagger$ }
\label{ss.levyweb}

In~\cite{gh1} we saw that the measure-valued process leads to a decomposition of $\Omega_N^\infty$ into monotype components, which are unions of disjoint balls. We want to better understand the time-space picture of these decompositions and their \emph{evolution in time}. The idea is to monitor the decomposition of the \emph{population of lineages} at time $t$ into descendants of the \emph{founders} at time $0$. More generally, the dependence structure in time arises from the decomposition into subfamilies with most recent common ancestors (MRCA's) at times $s < t$. This requires passing from populations of \emph{individuals} to populations of \emph{paths} of these individuals. We take a basic idea from~\cite{D93,DP91} on the historical process tracing individual paths in geographic space, and enrich it by combining it with our genealogical process $\mfU$ to processes of populations of paths in geographic-type \emph{and} ancestral space, possibly with added \emph{types}.

The present section makes these ideas rigorous by using a \emph{new concept} from~\cite{Lambda,gklimw}, appearing in two forms: the \emph{ancestral web} and the \emph{ancestral coupling}. We construct an embedding of these objects together with the process $\mfU$ in \emph{one} large space, the \emph{ancestral space} $\mfX$, respectively, the fibred measure space $\mfX^\dagger$.

By embedding representatives of all genealogies $\mfU_t$ of individuals alive at time $t \geq 0$ into the genealogy of \emph{all} individuals \emph{alive at some time} $t \geq 0$, in one common \emph{ancestral space} $\mfX$, we obtain the fibred measure version $\mfX^\dagger$. The equivalence class of the \emph{grounded measure} (respectively, \emph{fibred measures}) \emph{$\R$-tree of all individuals ever alive} then contains, for the process $\mfU$, the equivalence classes of \emph{marked ultrametric measure spaces} at \emph{each time $t$}, but \emph{not} connecting the representatives at different $t$ into a single state.

In the space $\mfX$ the \emph{geodesics} from the founders at time $0$ (or any time $s$) to the population at time $t$ ($>0$, respectively, $>s$) are the \emph{lineages} (ancestral paths). Oriented the other way they are the \emph{paths of descent}. We introduce the \emph{$(\Omega_N^\infty \times \K)$-marked} version of ancestral paths and paths of descent as \emph{marked oriented geodesics} in an \emph{equivalence class} of \emph{grounded tree-like metric fibred measure spaces}.


\subsubsection{Concept and construction of ancestral web, coupling as stochastic process}
\label{sss.levyweb}

These objects are constructed for our population model on $\Omega_N^\infty$, based on the fact that we have them on \emph{countable} geographic spaces, in particular, on $\Omega_N$, respectively, on $\Omega_N^L$. The technique is to again work with the $L$-truncation of the dynamics living on countable geographic spaces and pass to the limit $L \to \infty$. The processes are \emph{explicitly constructable} in \emph{individual-based} models on \emph{finite} geographic space.

On $\Omega_N^\infty$ we have so far introduced the process $\mfU$, constructed via duality. We now construct the corresponding
\begin{equation}
\label{e742}
\text{\emph{ancestral web} $\mfW = (\mfW_T)_{T\geq 0}$},
\end{equation}
as a population of ancestral \emph{paths} with time index $t \in [0,T]$, evolving up to \emph{time horizon $T$}. This is again constructed via the \emph{dual} process, the \emph{L\'evy-web} $\wh \mfW$ of the \emph{paths of descent} going backwards, driven in backward time by spatial coalescents.

The final step is to
\begin{equation}
\label{e746}
\begin{aligned}
&\text{embed $\mfW$ and $\wh \mfW$ as paths in one \emph{ancestral space} $\mfX$,}\\
&\text{respectively, $\mfX^\dagger$ when fibred with respect to time,}
\end{aligned}
\end{equation}
to obtain a \emph{graphical construction} of the whole genealogical structure, forward \emph{and} backward \emph{simultaneously}. In this way, $\mfW$ and $\wh \mfW$ are \emph{graph and dual graph} in the time-space continuum graph $\mfX$, viewed as a \emph{net of geodesics}.

\medskip\noindent
{\bf The forward-backward pair.}
We construct, for the given \emph{genealogy process} $(\mfU_t)_{t \in [0,T]}$, 
\begin{equation}
\text{a \emph{stochastic process} $\mfW=(\mfW_T)_{T \geq 0}$ with \emph{dual process} $\wh \mfW^T=(\wh \mfW_t^T)_{t \in [0,T]}$}
\end{equation} 
at time horizon $T>0$, given by the \emph{L\'evy-web} $\wh \mfW=(\wh \mfW^T)_{T >0}$. The process $(\wh \mfW_t^T)_{t \in [0,T]}$ is the \emph{dual (backward) process} for time horizon $T$ of the forward process $\mfW=(\mfW_T)_{T \geq 0}$ of populations of marked lineages evolving as a Markov process. Later we embed both in \emph{one} space $\mfX^\dagger$, so that $\wh \mfW$ arises from $\mfW$ by \emph{reorienting} the \emph{forward} geodesics from $0$ to $t$ to ones running \emph{backward} from $t$ to $0$.

The construction is based on the information obtained in~\cite{gh1} about the measure-valued process $\wh X^\infty$, and in~\cite{gklimw,Lambda} about the process $\mfU$ on countable geographic spaces. Furthermore, $\wh \mfW$ allows us to \emph{construct} $\mfW$ by using the technique of \emph{constructing Markov semigroups via duality} developed in~\cite{DGP23}. We obtain the genealogy and embed the whole process $\mfU$ into the ancestral web $\mfW$ for the dynamics on $\Omega_N^L$. We subsequently use this to get, by $L$-truncation and $L \to \infty$, a limit process satisfying a limiting duality relation, from which we construct the Markov process for the \emph{uncountable} geographic space $\Omega_N^\infty$.

We obtain these various processes as unique solutions to corresponding \emph{martingale problems}, and establish the \emph{Markov property} and the \emph{Feller property}. We also give duality relations, including \emph{strong duality} and even sharper forms: a \emph{graphical} and \emph{strong} duality in $\mfX^\dagger$.

\medskip\noindent
{\bf Main ideas and concepts for individual-based models on finite geographic spaces.}
We start from finite space and population size models, where the new concepts are visible \emph{without technicalities} and we can argue intuitively.

In the theory of individual-based models, stochastic processes can often be explicitly constructed by \emph{random graphs}, the so-called graphical construction. This is the case in \emph{individual-based} models on \emph{finite} geographic spaces in population genetics, such as the Cannings model with block-resampling that we have here. The state space $E$ is of the form $\{1,\ldots,n\}^M$ with $n$ the number of individuals carrying a mark from $M=\K \times \G$, with $\K$ a compact \emph{type} space and $\G$ a finite \emph{geographic} space. One takes a \emph{time-space graph}, a cylinder of \emph{finite} graphs
\begin{equation}
\label{e778}
\mathfrak{X} \subseteq [0,\infty) \times E,
\end{equation}
and constructs explicitly the evolution of the population. The construction starts with \emph{Poisson point processes} generating \emph{death points} and \emph{arrows} specifying \emph{birth events} of new descendants. This gives an \emph{explicit representative} for the states of the genealogy process $\CU$ corresponding to the population alive at time $t$, evolving in $t$ driven by this \emph{random graph} based on the ancestral space $\mfX$.

Looking \emph{forward}, we see the evolving population of \emph{lineages} leading from ancestor to descendant, marked with type and location. Looking \emph{backward}, we see the \emph{path of descent}, going successively back to the ancestors at earlier times, decorated with marks (type and location). By assigning weight one to every lineage alive at time $t$, we obtain a \emph{measure on lineages} representing the relevant information on the population.

In summary, we have a \emph{population of marked lineages} that \emph{evolves forward} as a stochastic process $\mfW$ in time, typically a \emph{strong Markov and Feller process} specified by a \emph{well-posed martingale problem} in $D([0,\infty), \A^\ast)$ with values in $\M^{V,\rm tree}$, in fact, in a subspace $\A^\ast$. This is important when we let the population size per site go to infinity to obtain a \emph{limiting measure-valued state} as description of the set of lineages.

Equipping the set of all individuals alive at some time with the \emph{genealogical distance} and taking the closure, we obtain the \emph{ancestral space $X$} of the population. From every point alive at time $t$ there is a unique ancestor at time $0$, and hence a unique ancestral path (a \emph{geodesic}), viewed forward representing a lineage. All these lineages carry a measure $\mu^\ast_t$ giving weight one to each individual alive at time $t$.

Thus, we have created from the \emph{random graph $X$} two objects:
\begin{itemize}
\item 
\emph{Pre-ancestral web} $W=(X,r,(\mu^\ast_t)_{t \geq 0})$: the \emph{population of lineages ever alive}.
\item
 \emph{Pre-ancestral coupling} $X^\dagger=(X,r,(\mu_t)_{t \geq 0})$: via the \emph{ancestral space}, the genealogical \emph{marked $V$-fibred measure $\R$-tree} $X^\dagger$ of \emph{all} individuals ever alive.
\end{itemize}
The equivalence classes $\mfW$ and $\mfX^\dagger$ are the final objects (without ``pre-''). Here, $\mu_t$ is the projection of $\mu^\ast_t$ onto the endpoint of the path at time $t$. The states of our stochastic processes are defined as the \emph{equivalence classes}
\begin{equation}
\label{e789}
\mfW = \left[X,r,(\mu^\ast_t)_{t \geq 0}\right], \qquad \mfX^\dagger=\left[X,r,(\mu_t)_{t \geq 0}\right],
\end{equation}
where equivalence is with respect to measure-preserving isometries. For these objects we construct
\begin{equation}
\label{e875}
\text{\emph{Polish state spaces} $\A^\ast$, respectively, $\A$ in Section~\ref{sss.formaldef}.}
\end{equation}

Next, we pass to \emph{$V$-marked} objects: we replace $X$ by $X \times V$, $\mu^\ast$ by measures $\nu^\ast$ on $V$-marked ancestral paths with marks in $V$, and $\A^\ast,\A$ by $\A^{\ast,V}$ respectively $\A^V$.


\subsubsection{Generalisation to continuum frequency models}
\label{sss.resogen}

The objects in~\eqref{e789} can also be defined in the \emph{infinite-population-per-site limit} and for \emph{infinite} geographic spaces like $\Omega_N$ via classical approximation techniques. Such processes have been constructed on \emph{countable} geographic spaces (see~\cite{gklimw}). To obtain them on $\Omega_N^\infty$, we follow the same line of argument as for the measure-valued process $\wh X^\infty$ and the genealogy process $\mfU^{\infty}$ on $\Omega_N^\infty$, using \emph{$L$-truncation} combined with the \emph{dual representations} of $\mfU^L$ available on $\Omega_N^L$.

The state space now describes populations of paths rather than individual points, which creates technical difficulties that are known from the study of historical processes~\cite{DP91}. Therefore, even though the totality of results is parallel to the previous theorems, we do not proceed entirely in parallel in order to minimise technicalities.

\begin{theorem}[$\mfW$ as limiting $\mfW^L$]
\label{th.mfw}
Consider $\mfW^L,\mfX^{L,\dagger}$ on $\Omega_N^L$ arising from the $L$-truncation of the objects on $\Omega_N^\infty$. Then a new process on $\Omega_N^\infty$ is obtained, namely, the limit $L \to \infty$ of  
\begin{equation}
\label{e908}
\text{$\CL[\mfW^L], \CL[\mfX^{L,\dagger}]$ exists in $\A^{\ast,V}$, respectively, $\A^V$.}
\end{equation}
The limit is denoted by $\mfW$, respectively, $\mfX^\dagger$ (dropping the $\infty$ index). These processes are strong Markov and Feller processes with c\`adl\`ag paths, and with continuous paths when there is no block-resampling.
\end{theorem}

Again we can determine the behaviour of the process $\mfW$ as time tends to $\infty$.

\begin{theorem}[Longtime behaviour $\CL(\mfW_T)_{T \geq 0}$]
\label{th.858}
The process $(\mfW_T)_{T \geq 0}$ converges as $T \to \infty$ to an equilibrium random variable:
\begin{equation}
\label{e861}
\CL [\mfW_T] \Rightarrow \CL [\mfW_\infty].
\end{equation}
The variable $W_\infty$ can be expressed, via strong duality, as the long-time limit $T \to \infty$ of the dual limit state $\wh \mfW_\infty$ of the L\'evy-web $\wh \mfW$, which exists almost surely (see Theorem~\ref{P.Levy} below). 
\end{theorem}

\medskip\noindent
{\bf Martingale problem for $\mfW$.}
We next give a forward characterisation of the \emph{ancestral web process $\mfW=(W_T)_{T \geq 0}$} as the \emph{solution} to a \emph{martingale problem}. The process $\mfW_T$ describes the ancestral paths of all individuals alive before time horizon $T$. For this, consider ancestral webs over time intervals $[0,T]$ and vary $T$, to obtain a stochastic process $(\mfW_T)_{T \geq 0}$. The evolution of $\mfW$ is driven by $\mfU$, and the operators from~\eqref{pol0} and the rest must be lifted from $\U^V$ to operators on $\A^{\ast,V}$ or $\A^V$.

The test functions are \emph{polynomials} denoted $\Pi^{\rm web}$, defined on measures on \emph{paths} instead of points in $U \times V$. The polynomials represent drawing $n$-samples of paths from the sampling measure $\mu^\ast_t$ (or $\nu_t^\ast$ in the marked case) at $m$ time points $t\leq T$, from the \emph{population of paths}, and afterwards applying specific evaluations. These evaluations observe the path at finitely many time points, giving position in ancestral space, type, and location. The formal definitions are in Definition~\ref{def.991} below for test functions and Definition~\ref{def.995} below for operators.

\begin{theorem}[Martingale problem ancestral web]
\label{th.843}
The process $(\mfW_T)_{T \geq 0}$ is the solution of a well-posed $(\mfL^{\rm web},\Pi^{\rm web},\CL[\mfW_0])$-martingale problem on $D([0,\infty),\A^{\ast,V})$.
\end{theorem}

We return to the rigorous and formal concepts for \emph{webs} and \emph{ancestral couplings} in Section~\ref{ss.formal}.


\subsubsection{The classical and strong graphical duality of $\mfW$ and the dual web $\wh \mfW$}
\label{sss.dualW}

Our goal is to give a \emph{graphicalal representation} of both the \emph{forward processes} $\mfU$ and the ancestral web $\mfW$ on $\Omega^{\infty}_N$, embedded in one space evolving as process $\mfX$ in the time-horizon $T$, and to relate it to a backward process $\wh \mfY$, respectively, $\wh \mfW$. Both backward processes are driven by the evolution of the populations of paths created by a \emph{$V$-marked genealogy-valued spatial coalescent}. We embed both forward and backward processes into $\mfX$, to obtain a \emph{strong graphical duality}. Here, we lift the duality $\mfU,\wh \mfY$ already stated to the level of paths.

\medskip\noindent
{\bf Strong duality between $\mfW$ and $\wh \mfW$.}
We first present the \emph{strong} duality between $\mfW,\mfU$ and $\wh \mfW,\wh \mfU$, and in the next section the construction of the graph $\mfX$ in which both $\mfW$ and $\wh \mfW$ as forward and backward \emph{oriented graphs} can be embedded. The genealogies forward and backward are adapted as lineages, respectively, paths of descent, yielding a strong \emph{graphical} duality.

We need to pass to $\U^{\Omega_N^\infty \times \K}$-valued states involving measure-valued states: measures on $U \times (\Omega_N^\infty \times \K)$ from mark-marginals in $\CM_{\rm loc}(\Omega_N^\infty \times \K)$, lifted to measures $\mu^\ast_t$ on \emph{paths} in ancestral space with marks in $\Omega_N^\infty \times \K$. The basic idea of the \emph{historical process} $\mfH$ of Dawson and Perkins~\cite{DP91} is that it takes values in measures on \emph{paths} with values in $\Omega_N^\infty \times \K$, associated with our measure-valued process $\wh X^\infty$. The ancestral web $\mfW$ enriches $\mfH$ with \emph{genealogical information} and is constructed with marks in $\Omega_N^\infty \times \K$.

The \emph{historical process} $\mfH$ turns out to be a \emph{Markovian functional} of the ancestral web $\mfW$. By definition, $\mfW$ also \emph{includes} the process $\mfU$ of genealogies of individuals alive at time $t$ and the \emph{fossil process} of genealogies of \emph{all individuals ever alive up to time $t>0$} as a functional. With the knowledge of the path of $W$ embedded in $X$ we can construct the path $U$, and with it $\wh W$ and the object $\wh U$. It therefore suffices to obtain
\begin{equation}
\label{e925}
\mfW,\,\wh \mfW,\,\mfX^\dagger.
\end{equation}

\medskip\noindent
{\bf The L\'evy-web as dual.}
To construct $\mfW$, we use \emph{duality} and therefore first construct the \emph{dual process $\wh \mfW$}, the \emph{L\'evy-web}. To obtain this process  on $\Omega_N^\infty$ we use the duality for a process on a \emph{countable} group $\Omega_N^L$ embedded into $\Omega_N^\infty$, denoted $\wh \mfW^L$. This process is known and given for the backward paths by the genealogies of spatial coalescents. We want to \emph{lift} the duality of $\mfU$ with $\mfY$ from Theorem~\ref{th.588} to a duality of $\mfW$ and the dual random graphs given via backward paths, called $\wh \mfW$. Here we lift the duality function $H$ to $H^{\rm web}$, a function based on monomials of \emph{measures on paths}, developed in~\eqref{e1006} below. Eventually we want $(\mfU,\mfW),(\wh \mfY,\wh \mfW)$ \emph{all} embedded in one space, the \emph{$V$-marked ancestral space} $\mfX$ with an extra \emph{fibred} structure called $\mfX^\dagger$.

For every individual alive at time $t$ in a representative of $\mfW$, we consider the lineages (geodesics oriented forward), which viewed backward with marks give a $(\Omega_N^\infty \times \K)$-marked path backward as the \emph{marked path of descent}. This induces a backward object $\wh \mfW$ as the equivalence class. The initial measure of the backward path is taken as $\mfh$ on $U=\Omega_N^\infty$ at time $t$. We obtain $\wh \mfW$ as a \emph{backward} process running from time $t$ to time $0$.

\begin{theorem}[Duality: ancestral web $\mfW$ and L\'evy-web $\wh \mfW$]
\label{P.Levy}
The duality relation of $\mfW$ and $\wh \mfW$ is obtained after replacing in Theorem~\ref{th.588} on duality the processes $\mfU,\mfY$, the duality function $H$, and the state spaces $\U^V$ by $\mfW$, $\wh \mfW$, $H^{\rm web}$ and $\A^{\ast,V}$, respectively. The formula for $H^{\rm web}$ is developed in Section~\ref{sss.sto} in~\eqref{e1451}. After switching the orientation of the lineages to run backward from time $t$ to time $0$, the law of $\wh \mfW$ is obtained from the law of $\mfW$, which yields the strong duality. 
\end{theorem}

We will strengthen this further to a strong \emph{graphical} duality: the forward and backward processes of lineages, respectively, paths of descent both embedded in a single ancestral space $\mfX$ with a further structure of \emph{fibred measures} $\mfX^\dagger$, with individuals characterised by the current time and location and marked with a type from $\K$.

We next state the \emph{key result} on ancestral webs, representing the final strengthening of Theorem~\ref{P.Levy}.

\begin{theorem}[Strong graphical duality of $\mfW,\wh \mfW$ in the ancestral space $\mfX^\dagger$]
\label{th.dualityweb}
$\mbox{}$\\
The process $\wh \mfW$ is the strong graphical dual of $\mfW$. There exists an ancestral space $\mfX^\dagger$, a grounded tree-like metric fibred measure space modelling all individuals ever alive, with the following properties. 
\begin{itemize}
\item[(a)]
This space contains the $V$-marked geodesics connecting time-$0$ ancestors with time-$t$ descendants in $\mfW$. 
\item[(b)]
This space also contains the lines of descent as a representative of $\wh \mfW$: $V$-marked geodesics connecting a point in the backward $0$-fibre with a point in the backward $t$-level. 
\item[(c)]
Representatives $\wt \mfW$, respectively, $\overset{\sim}{\wh \mfW}$ of $\mfW$, respectively, $\wh \mfW$ can be embedded into one representative $\wt \mfX^\dagger$ of $\mfX^\dagger$. The forward paths in $\wt \mfW$ and the backward paths in $\overset{\sim}{\wh \mfW}$ give the forward, respectively, backward oriented geodesics in $\wt \mfX^\dagger$, the ancestral coupling from the time-$0$ fibre to the time-$t$ fibre for all $t>0$.
\end{itemize} 
\end{theorem}

\begin{remark}[Infinite volatility limit]
\label{r.1034}
{\rm The limit $d \to \infty$ in $\mfW^{[d]}$ and $\mfX^{[d]}$ yields a process $\mfW^\sim$, an $\Omega_N^\infty \times \K$-marked L\'evy-web based on \emph{instantaneously coalescing L\'evy-processes}, such that Theorems~\ref{th.mfw}--\ref{th.858} and~\ref{P.Levy}--\ref{th.dualityweb} hold.} 
\end{remark}


\subsubsection{Construction of the L\'evy-web $\wh \mfW$ in a graph $\mfX^\dagger$ by $L$-truncation}
\label{sss.levy-web}

To obtain a solution to the martingale problem, we construct the stochastic process $\mfW$ of ancestral webs by using the \emph{duality} with the \emph{L\'evy-web} $\wh \mfW$, which we now \emph{construct explicitly}. This construction significantly enhances its value beyond it being the key to the proofs of the other theorems of Section~\ref{ss.levyweb}. Namely, it allows us to \emph{verify the operators} for the martingale problem and to get the \emph{uniqueness} of the martingale problem.

We first construct $\wh \mfW$, keeping in mind that we later want to realise $\mfW$ and $\wh \mfW$ both in one metric measure space, the \emph{ancestral space} $\mfX$ enriched to $\mfX^\dagger$ (a $V$-marked grounded metric \emph{fibred} measure space) of all individuals ever alive, equipped via $\mfW,\wh \mfW$ with a \emph{path of fibred measures} representing the \emph{ancestral coupling} $\mfX^\dagger$. We use the construction on countable geographic spaces, in particular, on $\Omega_N^L$, investigated in \cite{gklimw,Lambda}. We try to obtain the \emph{ancestral web for $\mfU$} by \emph{$L$-truncation approximation} by the ancestral webs for the $\mfU^L$, $L \in \N$, and denote by $\mfW^\infty$ the limit. We then show that $W^\infty = \wh W$.

\medskip\noindent
{\bf Convergence and initial state.}
We obtain a (random) metric space based on the basic set $\mfD$, i.e., $(\mfD,r)$. The separable metric space is \emph{not} complete, so we take the \emph{closure} to get $\bar \mfD=\Omega^{\infty}_N$, and for the corresponding $\wh \mfW,\mfX^\dagger$ in the topology of $\A^{\ast,V}$ , respectively,\ $\A^V$.

We consider the equivalence classes of marked metric fibred measure spaces based on $\cup_{k=0}^L D_k$, say $(\mfW^L)$. We need the sequence $\{\CL[(\wh \mfW^L)]\colon\, L \in \N\}$ to be \emph{tight} and \emph{converging} to a limit law $\CL[\wh \mfW^\infty]$. For this we must show \emph{convergence} in $\A^\ast$ of the approximating sequence of marked metric measure spaces, based on the $L$-truncation with $\Omega_N^{\infty,\rm disc}$.

The initial state of the dual is located at the forward time $t$, but its state is equipped with metric, measure, and types by running its dual dynamics backward in time $s$ from $t$ to $0$. The initial state for the backward process is, for all $t>0$, given by the dual population equipped with the \emph{Haar measure} $\mfh$ on $\Omega_N^\infty$, independent of $t>0$.

\begin{theorem}[Construction of L\'evy-web $\wh \mfW$ via $\wh \mfW^\infty$]
\label{prop.1382}
$\mbox{}$\\
For the $L$-truncation of $\Omega_N^\infty$ by $\Omega_N^L$, in the limit as $L \to \infty$, in the space $\A^{\ast,V}$,
\begin{equation}
\label{e1384}
\CL[\wh \mfW^L] \Longrightarrow \CL[\wh \mfW^\infty] \in \CM_1 \left( \A^{\ast,V} \right).
\end{equation}
\end{theorem}

We want to show that $\wh \mfW^\infty$ allows for a dual representation of the ancestral web $\mfW$, and hence equals $\wh \mfW$. Even more, we want to characterise it as the \emph{strong} dual to $\mfW$, meaning that $\CL[\mfW]=\CL[\wh \mfW]$, as follows from Theorem~\ref{P.Levy}.

\begin{corollary}[Identification of $\wh \mfW$]
\label{cor.986}
\begin{equation}\label{e987}
\CL[\wh \mfW] = \CL[\wh \mfW^\infty].
\end{equation}
\end{corollary}

Are there alternative approximations of $\Omega_N^\infty$? The answer is no: the limiting object on $\Omega_N^\infty$ is unique. This amounts to showing independence of the limit from the choice of approximating countable dense subsets. Equivalently, the dual process $\wh \mfW$ is uniquely given by a process driven by a \emph{spatial coalescence process}, analogous to the processes $\mfU$ and $\wh X^\infty$ on $\Omega_N^\infty \times \K$, but extended to a \emph{time-space dual process}. Here, coalescent individuals are activated only after backward times $s_1<s_2<\ldots<s_n<t$, with $s$ running from $0$ to $t$. At backward time $s$, all paths are equal beyond time $s$ and until time $0$ are in one partition element. Any $n$-tuple sampled from $\mfX$ is a spatial coalescent with parameters $((c_k)_{k \in \Z}, d^\ast,(\lambda_k,\Lambda_k)_{k \in \Z})$ for the time-space version. We denote by $\theta$ the marginal in $\K$ of the measure $\wh X_0$.

\begin{theorem}[Stochastic process characterisation of the L\'evy-web $\wh \mfW$]
\label{th.600}
$\mbox{}$\\
If $[\wt \mfW]=\wh \mfW^\infty$, then the following hold:\\
(i) For every $z\in \R \times (\Omega^{\infty}_N \times \K)$ there is a unique mark in path space: the path $\kappa(z)$ in $\wt \mfW$ a.s.\\
(ii) For every deterministic set of points $z_1, \ldots, z_n$, with $z_i \in \R \times (\Omega^{\infty}_N \times \K)$ for $i=1, \ldots, n$, the collection $\{\kappa(z_1), \ldots, \kappa(z_n)\}$ of marks is a version of a $\K$-typed \emph{space-time} $(\uc,\ulam, \uLam, d^\ast)$-\emph{block-coalescing L\'evy-process} starting in $(z_1, \ldots, z_n) \in \R^+ \times (\Omega_N^\infty \times \K)$, $\ell^{\otimes n} \otimes (\mfh^{\otimes n} \otimes \theta^{\otimes n})$-a.e.\\
(iii) For every deterministic countable dense subset $\CD \subseteq \R \times \Omega^{\infty}_N$, $\wh \mfW$ is a.s.\ the increasing limit of sequences of countable discrete sets with closure $\Omega_N^\infty$ in $\M^{D^-}$ of $[\wh \mfW^\CD]$, in particular, for $D=\mfD$ (see~\eqref{e1156} below for the definition).
\end{theorem}

We need an approximation property showing that the same object is obtained for an \emph{arbitrary choice of the countable dense subset of $\R \times \Omega^{\infty}_N$} used for a countable set of starting points of the coalescent. This plays a role similar to the tightness-and-convergence lemma in \cite{NRS05}. Using the $L$-truncation and ancestral webs on $\Omega_N^L$, we show existence of $\wh \mfW$ on $\Omega_N^\infty$, and subsequently of $\mfW$ on $\Omega_N^\infty$ by duality. The next section states the main tool for obtaining the \emph{strong graphical duality} relation between $\mfW,\wh \mfW$ in the \emph{ancestral space $\mfX$} equipped with fibred measures $\mfX^\dagger$.

\medskip\noindent
{\bf Excursion on a functional of $\mfU$: Coalescent with rebirth.}
Looking instead of $\wh \mfW$ at the \emph{time-space process} of evolving populations, like the \emph{fossil process} of the genealogy of all individuals ever alive up to time $T>0$, with sampling measure
\begin{equation}
\label{e1043}
\mathfrak{\ell}_{[0,T]} \times \nu_T, \text{  with  } \bar \nu_T =  \int_0^T \nu_t\, \dd t,
\end{equation}
another process and form of duality appears. In this process we \emph{cannot} identify the population genealogy at a fixed fibre, since we resample times from the measure $\mathfrak{\ell}$. The fossil process is a functional of the ancestral web. It focuses on occupation measures rather than on the totality of all fibres. The duality is based on time-space functions, and the dual process is the \emph{coalescent with rebirth}.

Starting a coalescent on $\Omega^{\infty}_N$ from every time point (rather than a sequence of successively refined discrete grids) can, in the entrance law and strong dual representation, be realised as a \emph{coalescent with rebirth}: the dual of the time-space process of $\mfW$~\cite{gklimw}. The finite dual version works as follows. Consider the spatial coalescent with block-resampling. Once block-resampling occurs in a ball, for each partition element that coalesces a \emph{new} basic individual is created at the same location and activated at the same time. For Fleming-Viot resampling, the absorbed individual is complemented by a \emph{newly created individual} at the \emph{death time}.

This mechanism extends for fixed time horizon $t>0$ to an \emph{entrance law} from locally countably many individuals (one at each site of $\Omega_N^{\infty, \rm disc}$), starting at each time point from the dyadic numbers, or even locally uncountably many. This process connects more closely than the time-space coalescent of the populations corresponding to \emph{different time horizons}.

Extension to an entrance law of \emph{uncountably} many initial individuals requires a \emph{limit} of discrete time nets. This arises as the $L \to \infty$ limit, with an adapted shrinking of the time-points that uses the dyadic numbers, introduced in~\eqref{e1156} in Section~\ref{sss.detail} below, and verified in Section~\ref{ss.pexisten} below. This gives us $\wh \mfW$ as a \emph{backward Markov process} from every time horizon $t>0$ backwards until time $t$.


\subsubsection{The frequency process obtained via the graphical construction with the L\'evy-web}
\label{sss.contcannings}

We next construct $\wh X^\infty$, the measure-valued process on $\Omega_N^\infty \times \K$ treated in~\cite{gh1}, with the help of the L\'evy-web. This gives a \emph{graphical construction of} $\wh X^\infty$.

A state arises on the individual-based level from translation-invariant initial states on $\Omega_N$. The initial state for the continuum limit must therefore satisfy: for fixed $B \in \CB(\K)$, the sampling measure $\nu$ projected on the marks $\nu^\downarrow$ satisfies
\be{ag40a}
\nu^\downarrow (\,\cdot \times B) \mbox{ is the normed Haar measure $\mfh(\cdot)$ on $\Omega^{\infty}_N$}.
\ee
Recall that the migration and reshuffling would give a mass flow if there were no resampling, which may induce a \emph{stochastic} part in the evolution depending on the parameters.

There are two possibilities in the stochastic regime. 
\begin{itemize}
\item[(i)] 
The \emph{regular} case: for finitely many types, the states have a density with respect to the Haar measure that is locally constant. 
\item[(ii)] 
As shown in~\cite{gh1}, there is no stochastic regime with \emph{singular} states. Without a Cannings part, states are either regular stochastic or deterministic mass flow. In the \emph{regular stochastic} case, we have a collection of random sets each of one type, consisting of a disjoint union of disjoint balls.
\end{itemize}

\medskip\noindent
{\bf Construction.}
Suppose that we are in the regular stochastic regime. We define a measure-valued process $\nu_t$ on a probability space $(\Omega,\CA,P)$ obtained by taking the realisation of the \emph{L\'evy-web} $\wh \mfW_t$ arising from the \emph{time-space L\'evy-coalescent} run for backward time $t$, as constructed in Section~\ref{sss.detail} below. Let $\mfZ_t (\xi)$ be the \emph{location} of the partition element of $(t,\xi)$, denoted $C_t(\xi)$, after evolving for backward time $t$. Assume
\begin{equation}
\label{e1570}
\wh X^\infty_0 = \mfh \otimes \theta.
\end{equation}
Let $(Y(u))_{u \in \Omega^{\infty}_N}$ be i.i.d.\ $\theta$-distributed $\K$-valued random variables. Let $\wt X_t^\infty$ be the projection of $\nu_t$ onto $(\Omega_N^\infty \times \K)$. Define
\be{con1}
\wt \mfX^\infty_t = \left(\frac{\dd}{\dd \mfh} \right)\,C_t (\xi),  
\qquad (\xi,v) =
\begin{cases}
1, & \mbox{ if } Y(\mfZ_t(\xi)) = v, \\
0, & \mbox{ otherwise.}
\end{cases}
\ee

\begin{theorem}[Representation of $\wh X^\infty$ via L\'evy-web]
\label{T.XLW}
In the stochastic regime, the process $\wt X^\infty$ is a version of the process $\wh X^\infty$ defined in~\cite{gh1} as a $\CM_{\rm loc} (\Omega_N^\infty \times \K)$-valued process.
\end{theorem}

\begin{remark}
\label{r.1547}
{\rm On $\Omega^{\infty}_N$, in a given ball we have a patchwork of arbitrarily small balls carrying different types, leading to countably many small jumps. On $\Z$, respectively, $\R$, the situation is even more irregular. There, under the Lebesgue measure, we do not sample two individuals from one point, yet the evolution is driven by a \emph{countable} random set of \emph{spatial points} where \emph{resampling occurs}. These are the double points of the Brownian web, separating different subfamilies belonging to one initial ancestor.} 
\end{remark}


\subsection{Universality, continuum space rescaling limit of $\mfU,\mfW,\mfX^\dagger$}
\label{ss.csl}

We next carry out, for the processes $\mfU,\mfW,\mfW^\dagger$ on $\Omega_N$ and $\Omega_N^\infty$, the \emph{continuum space (rescaling) limit}, both at the level of the \emph{genealogy} process and the \emph{ancestral web process}, parallel to the procedure in the measure-valued case on $\Omega_N^\infty \times \K$ in \cite{gh1}. All the classical scaling operations of (i) time, (ii) space, (iii) rates, and (iv) measures, are \emph{exactly the same} as in the $\CM_{\rm loc} (\Omega_N^\infty \times \K)$-valued case. This provides a representation of the \emph{universality classes} of our models by the $(0,\infty)^2$-indexed collection (i.e., $(c,\mu) \in (0,\infty)^2$) of \emph{stable genealogy} processes and \emph{stable ancestral web} processes with parameter $c \leq 1$.

We scale the \emph{genealogical distances} in the same way as time, measuring them as \emph{twice the macroscopic time}, giving a process (with $(L)$ stressing the scaling parameter):
\be{ag21a1}
(\mfU^{(L)}_t)_{t \geq 0}.
\ee
In order to perform a rescaling, we embed all our systems on the geographical spaces $\Omega_N$ and $\Omega^L_N$ in one geographic space $\Omega^{\infty}_N$ and view all the states as measures on $\Omega^{\infty}_N \times \K$.

We define a motion on $\Omega^L_N$ induced by the scaled motion on $\Omega_N$, by \emph{shrinking} the jump size but \emph{increasing the jump rate} for compensation via a time rescaling. We scale the {\em time} and {\em space} such that the \emph{migration/reshuffling part converges to a non-trivial limit dynamic}: the \emph{effective migration}. We have to keep in mind that migration of individuals and the motion introduced by the reshuffling at the action of block reshuffling can be in different possible relations, which affects the scaling. For the scaling limit we establish the \emph{dichotomy}: either the genealogical process is \emph{stochastic} (the local population has a \emph{single ancestor} but is otherwise of non-trivial structure), or it is \emph{deterministic} (the population is locally totally fragmented in \emph{single-ancestor families}, which form ``dust'' as we aggregate over more and more sites, so that \emph{no} corresponding ultrametric marked measure tree exists).

We show that in the deterministic regime there is dust formation in $(\CU^{(L)}_t)_{t \geq 0}$ as $L \to \infty$ and \emph{no} limit in $\U_{\Omega^{\infty}_N}$ exists. The special regime is where the continuum limit migration is only Harris-recurrent: the dual process after scaling has an MRCA at a random time that converges weakly to infinity only as both time and $L$ go to infinity, rather than being already with positive probability equal to infinity as $t$ diverges for fixed $L$.

In the stochastic regime for the $\U^{\Omega_N^\infty \times \K}$-valued version, the \emph{continuum space rescaling} gives convergence to a limit process. We denote the limiting parameters $c,\lambda$ by adding a $\wedge$, \emph{given in~\eqref{ag8b} and \eqref{ag13b} below}.

\begin{theorem}{(Spatial continuum limit of $\U^{\Omega_N \times \K}$-valued process)}
\label{T.spatcont}
Assume that the model parameters satisfy~\eqref{ag7}--\eqref{ag8b} below. Define $\wh c,\wh \lambda$ by requiring that the parameters correspond to $(\wh c, 0)$, $(\wh c, \wh \lambda)$ with $\wh c = \wh \lambda$, or $(0,\wh \lambda)$, with always $\wh c +\wh \lambda >0$.
\begin{itemize}
\item[\textup{(a)}]
For the three cases $\wh c<1$, $\wh \lambda=0$ and $0 < \wh c = \wh \lambda < 1$ and $\wh c=0$, $0 < \wh \lambda < 1$, under the two conditions that the mark configuration is translation-invariant ergodic and the complete state is in $\U_{\Omega^L_N}$ and converges to $\mfU^{\infty,N}_0 \in \U^{\Omega^{\infty}_N}$, the following hold:
\begin{equation}
\label{ag21x}
\begin{aligned}
\CL\left[\big(\mfU^{(L)} (t)\big)_{t \geq 0}\right] &\Lto \CL \left[\big(\mfU^{\infty}(t)\big)_{t \geq 0}\right],\\
\CL \left[\big(\mfW^{(L)}(t)\big)_{t \geq 0}\right] &\Lto \CL \left[\big(\mfW(t)\big)_{t \geq 0}\right],\\
\CL \left[\mfX^{\dagger,(L)}\right] &\Lto \CL \left[X^{\dagger,\infty}\right].
\end{aligned}
\end{equation}
Here, $\mfU^{\infty}$ is the \emph{$(\wh c, \wh \lambda)$-stable genealogical Cannings process}, $\mfW$ the corresponding ancestral web, and $\mfX^{\dagger,\infty}$ the fibred ancestral space (where stable indicates that the genealogy process corresponds to migration/reshuffling parameters $(\wh c, \wh \lambda)$ with stable one-individual spatial motion).
\item[\textup{(b)}]
If $\wh c$, $\wh \lambda$ are in the parameter regime complementary to (a), then the collection $\CL[(\mfU^{(L)} (t)]$ is for $L \to \infty$ not tight in $\U^{\Omega_N^\infty \times \K}$ and exhibits dust formation.
\end{itemize}
\end{theorem}

\noindent
The scaling operation proceeds as follows:
\begin{itemize}
\item Replace jumps of distance  $k$ by
\be{ag2a}
\mbox{jumps to distance corresponding to $k-L$ {\em (shrinking),}}
\ee
i.e., the jump in the block of size at most $k$ on $\Omega_N$ is turned on $\Omega^{(L)}_N$ into a jump in the block of size $k-L$.
\item {\em Speed up} time by
\be{ag2}
t \la \frac{N^L}{c_L} t, \mbox{ respectively, } t \la \frac{N^L}{\lambda_L}
\ee
depending on whether $c \geq \lambda$ or $c < \lambda$.
\item The Fleming-Viot resampling at rate $d$ at a site of $\Omega_N$ occurs now at each site of $\Omega^{(L)}_N$ at rate $d c_L$ respectively $d \mu_L$, and similarly the Cannings rate at level-$L$ is scaled.
\item The block-resampling measure $\Lambda_k$ on $\Omega_N$ is now used for block-resampling in the $(k-L)$-block in $\Omega^{(L)}_N$.
\item Embed the $\Omega^{(L)}_N$-system into $\Omega^\infty_N$ and {\em rescale the mass} of the state (a measure on $\Omega^{\infty}_N \times \K$) by
\be{ag2b}
N^{-L}
\ee
to compensate for the increasing number of components of the $\Omega^{(L)}_N$-system in a site of $\Omega_N$.
\end{itemize}
Doing so, we obtain from $X$ on the $L$-th level of rescaling a process with values in $\CM(\Omega^{(L)}_N \times \K)$, which induces by embedding a process with values in $\CM(\Omega^{\infty}_N \times \K)$, called
\be{ag2c}
(X^{(L)}_t)_{t \geq 0}.
\ee

\medskip\noindent
{\bf Regularity assumptions for rescaling limit.}
For convergence, we need in particular that the migration/reshuffling converges to a limit. Assume therefore the following \emph{regularity} for the parameters:
\be{ag7}
c_k = \bar c_k c^k, \qquad \lambda_k = \bar \lambda_k \lambda^k,
\ee
with
\be{ag8}
\frac{\bar c_{k+L}}{\bar c_L} \Ltoo 1, \qquad \frac{\bar \lambda_{k+L}}{\bar \lambda_L} \Ltoo 1.
\ee
Furthermore, assume
\be{ag8b}
\Lambda_k | \Lambda_k |^{-1} \kto \Lambda.
\ee

\medskip\noindent
{\bf Determination of parameters of limit process.}
In view of \eqref{ag2a} and \eqref{ag2}, we calculate the sequences $(\wh c^L_k)_{k \geq -L}, (\wh \lambda^L_k)_{k \geq -L}$, depending on $c \geq \lambda$, respectively $c \leq \lambda$:
\be{ag10}
\wh c^L_k = \frac{1}{c_L} c_{k+L}, \mbox{ respectively, } \wh c^L_k = \frac{1}{\lambda_L} c_{k+L} \quad , \quad k= -L, -L+1, \ldots
\ee
\be{ag11}
\wh \lambda^L_k = \frac{1}{c_L} \lambda_{k+L}, \mbox{ respectively, } \wh \lambda^L_k = \frac{1}{\lambda_L} \lambda_{k+L} \quad , \quad k= -L, -L+1, \ldots
\ee
For the Fleming-Viot part, scale additionally
\be{ag12}
\wh d^L_{-L} =d, \qquad \wh d^L_k = 0  \mbox{ for }  k= -L+1, \ldots.
\ee

We have with \eqref{ag7}--\eqref{ag8} that, for $k \in \Z$,
\be{ag13}
\wh c^L_k \Ltoo c^k   \mbox{ for } c>\lambda, \quad  \wh c^L_k \Ltoo 0 \mbox{ for } c< \lambda,
\quad \wh c^L_k \Ltoo c^k \mbox{ for }  c=\lambda,
\ee
\bea{ag13b}
\wh \lambda^L_k &\Ltoo& \wh\lambda, \mbox{ for } \lambda >c,  \\
\wh \lambda^L_k &\Ltoo& 0, \mbox{ for } \lambda <c, \nonumber \\
\wh \lambda^L_k &\Ltoo& \wh \lambda^k = \wh c^k, \mbox{ for } \lambda =c. \nonumber
\eea
In other words, we can define a migration parameter $\wh c$ and a block-resampling parameter $\wh \lambda$ such that the limits in \eqref{ag13} and \eqref{ag13b} have the form $(\wh c)^k$, respectively, $(\wh \lambda)^k$.

Consider the scaling given above as used in~\cite{gh1} for the measure-valued process, now applied to the genealogy to get $\mfU$.


\section{Preparations: ancestral web, construction via dual, L\'evy-web}
\label{s.3}

In Section~\ref{ss.formal} we formalise the concepts concerning $\mfW$ and $\wh \mfW$. Section~\ref{ss.tools} explains why we can construct the process $\mfU$ from $Y^\infty$ and $\wh \mfY^\infty$, and $\mfW$, $\mfX^\dagger$ from $\wh \mfW$. The process $\wh \mfY^\infty$ is the central object in the proofs of the theorems of the genealogy process in Section~\ref{s.proofsXXX}. The process $\wh \mfW$ is the central object in the proofs of the theorems on the ancestral web and the coupling in Section~\ref{s.proofs}.  Section~\ref{s.prepa} outlines the construction of the dual L\'evy-web $\wh \mfW$ as preparation for the proof sections. The proofs concerning the construction of the \emph{L\'evy-web} follow directly in Section~\ref{ss.pexisten}.


\subsection{Formal introduction: ancestral web and coupling}
\label{ss.formal}

In the next two sections we turn the explanation of \emph{ancestral web and coupling} into formal definitions and give the ingredients for the \emph{martingale problems} for these objects.


\subsubsection{Formal definitions $1$: ancestral coupling and web}
\label{sss.formaldef}

Formal definitions of $\mfW$ and the ancestral coupling $\mfX^\dagger$ are needed, as well as state spaces for these random variables and processes. The ancestral coupling is a prerequisite for the ancestral web.

\begin{definition}[Pre-ancestral coupling]
\label{def.622}
$\mbox{}$
\begin{itemize}
\item
For a given $\U$-valued genealogy process $\mfU=(\mfU_t)_{t \geq 0}$ with $\mfU_t=[U_t,r_t,\mu_t]$, define the pre-ancestral coupling of $\mfU$ as
\begin{equation}
\label{e352}
\wt \mfX^\dagger = \left(X,r,\varrho,\leq,(\mu_t)_{t \geq 0} \right),
\end{equation}
where the following four requirements are made:
\begin{equation}
\label{e605}
\begin{aligned}
&\wt \mfX^\dagger = (X,r,\varrho,\leq) \text{ is a } \varrho\text{-rooted (or grounded),}\\ 
&\text{complete, separable, metric space and } \R\text{-tree},
\end{aligned}
\end{equation}
(the latter means that the metric $r$ satisfies the \emph{four-point condition}), with the natural partial order $\leq$ based on the root $\varrho$ giving the geodesic between points of $X$ and $\varrho$ an orientation (see~\eqref{e551}). In the grounded case $\varrho$ is replaced by a countable set, the \emph{ground level}.
\item
In addition, the process $\mfU$ is supported by requiring that
\begin{equation}
\label{e356}
X= \Big(\mathop{\bigcup}\limits_{t \geq 0} \;\wt U_t\Big)^c \text{ with }  \wt U_t=\{(t,\iota), \iota \in U_t\}, 
\quad  r_t \text{ is the restriction of $r$ to $\wt U_t$},
\end{equation}
with $u \in U_t$ if and only if $r(\varrho,u)=t$ (for exactly one element of the ground level in the grounded case) and where $u \leq u^\prime$ if and only if $u \in$ geodesic $(\varrho, u^\prime)$. In particular, $r_t$ is an ultrametric on $U_t$.
\item
Furthermore, for fibred measures: for all $t \geq 0$,
\begin{equation}
\label{e620}
\text{$\mu_t$ is a probability measure on the space $(U_t,r_t)$ with $t \to [U_t,r_t,\mu_t] \in \U$ c\`adl\`ag.}
\end{equation}
\end{itemize}

\end{definition}

\medskip\noindent
{\bf Functionals of rooted $\R$-trees.}
In an $\R$-tree the \emph{natural order} $\leq$ and the \emph{orientation} based on the root $\varrho$ are defined as
\begin{equation}
\label{e951}
\text{$\iota \leq \iota^\prime$ if $\iota \in$ geodesic $(\varrho,\iota^\prime)$, }
\end{equation}
giving the \emph{ancestor-descendant relation} of points in $X$, so that for every $(t,\iota^\prime)$ there exists an ancestor $(s,\iota)$ such that $(s,\iota) \leq (t,\iota^\prime)$ and
\begin{equation}\label{e551}
(s, \iota) \leq (t,\iota^\prime) \text{  iff  }  s \leq t \text{  and  } r\left((s,\iota),(t,\iota)\right) = t-s.
\end{equation}
We give an \emph{orientation} to a geodesic from $\varrho$ to some $\iota$ at distance $t$ from $\varrho$ by running through the path from $\varrho$ to $\iota$ in increasing order (respectively, from a point in the \emph{ground level} to $\iota$).

\begin{definition}[Ancestral coupling, state space $\A$]\label{def.614}\text{}\\
Fix a process $\mfU$ with values in $\U$.
\begin{itemize}
\item[(a)] 
The equivalence class
\begin{equation}
\label{e502}
\mfX^\dagger= \left[X,r,\varrho,\leq, (\mu_t)_{t \geq 0}\right]
\end{equation}
of different fibred rooted marked Polish metric generalised measure spaces $(X,r,\varrho,\leq, (\mu_t)_{t \geq 0})$ satisfying~\eqref{e605}--~\eqref{e620}, is defined with respect to $(r,\varrho)$-preserving bijections of $X$ that preserve the measures $\mu_t$ in each $t$-slice of $X$ for $t \in [0,\infty)$. This equivalence class is called the ancestral coupling associated with the $\U^V$-valued path $\mfU=(\mfU_t)_{t \geq 0}$, through the relation~\eqref{e356}.
\item[(b)] 
The equivalence classes of rooted $\R$-trees, equipped with sampling measures $(X,r,\varrho,(\mu_t)_{t \geq 0})$, are denoted by
\begin{equation}
\label{eq466}
\A
\end{equation}
and is equipped with a topology turning it into a Polish space. Convergent sequences are defined via the convergence of finite sampled metric subspaces, based on choosing any finitely many time levels and drawing a finite sample at each of those times with the corresponding sampling measures. This convergence is recalled in Remark~\ref{r.topo}.
\end{itemize}
\end{definition}

\begin{remark}[Topology]
\label{r.topo}
{\rm The topology on $\A$ is given in~\cite{gklimw,Lambda}. It is based on choosing time-points $t_1, \cdots,t_m$ and drawing at each time point $i \in \{1,\ldots,m\}$ independently and independently of all other samplings a finite sample of individuals according to $\mu_t$ from $U_t$. This leads to a \emph{finite substructure}: a metric space. We require convergence of all these substructures in sampling law to a finite limit in $\A$, with respect to the weak convergence of the distance matrix distribution under the sampling measure. This topology can also be formulated in terms of polynomials of the $\mu_t$ for some $t$.}
\end{remark}

An ancestral coupling associated with $\CU$ gives rise to an ancestral web associated with $\CU$. Note that there are at most countably many geodesics joining from level $s$ to level $t$ for $0<s<t$.

\begin{definition}[Ancestral web for $\CU$ and state space $\A^\ast$]
\label{def.801}
Let $\CX$ be an ancestral coupling of the $\U$-valued process $\mfU$.
\begin{itemize}
\item[(a)] 
Define a pre-ancestral web $X^\dagger$:
\begin{equation}
\begin{aligned}
&X^\dagger = \left(X,r,\varrho,\leq,(\mu_t^\ast)_{t \geq 0}\right),\\ 
&\text{where the $\mu^\ast_t$ induced by $\mu_t$ are measures on $C \left[[0,t],X \right]$},\\
&\text{by associating with a point $\iota \in U_t$ the geodesic from $\varrho$ to $\iota$.} 
\label{e397}
\end{aligned}
\end{equation}
\item[(b)] 
Take an equivalence class $\mfX^\dagger=[\wt X^\dagger]$. Two pre-ancestral webs
\begin{equation}
(X_1,r_1,\varrho_1,(\mu^1_t)^\ast_{t \geq 0}),\qquad (X_2,r_2,\varrho_2,((\mu^2_t)^\ast)_{t \geq 0}), 
\end{equation}
are equivalent when they can be mapped onto each other by a map $\varphi: X_1 \to X_2$, which is $(r,\varrho)$-preserving and measure-preserving for \emph{each} $\mu^\ast_t$, $t \geq 0$:
\begin{equation}
\label{eq862}
\varrho_1 \to \varrho_2, \; (\mu^2_t)^\ast  = (\mu^1_t)^\ast_\varphi, \;  
r^2 (\varphi(\iota),\varphi(\wt \iota)) = r^1(\iota,\wt \iota) \quad \mu_t^1 \otimes \mu^1_t -a.s. 
\text{ all } \iota,\wt \iota \text{  in  } U_t.
\end{equation}
(In the general case of a ground level $U_0,U_0^\prime$ instead of roots $\varrho^1,\varrho^2$ there is a bijection $\varphi_0:U^1_0 \to U^2_0$ mapping $r_0$ to $r_0^\prime$ and $\mu_0$ to $\mu^\prime_0$.) Denote by $[]$ the equivalence class. Define the ancestral web $\mfW$ as
\begin{equation}
\label{e713}
\mfW=\left[X,r,\varrho,\leq,(\mu^\ast_t)_{t \geq 0}\right].
\end{equation}
\item[(c)] Define $\wt \A$ as the set of all equivalence classes of rooted (grounded) $\R$-trees with measures $(\mu^\ast)_{t \geq 0}$ on the geodesics connecting $\varrho$ to points in distance $2t$ in the $supp(\mu_t)$. Let $\tau$ map a point onto the geodesic connection with $\varrho$. The map $\wt \tau$ induced by $\tau$ is invertible on its domain and image. Given a Polish topology on $\A$, it induces via $\wt \tau$ a Polish topology on $\wt \A$. This topological space is denoted by
\begin{equation}
\label{e740}
\A^\ast.
\end{equation}
\end{itemize}

\end{definition}

The topology on $\A^\ast$ was introduced in~\cite{gklimw,Lambda} and generates a Polish space. The necessary information is collected in Appendix~\ref{sa.space}.

\begin{remark}[Topology $\A^\ast$]
\label{r.987}
{\rm The topology can also be formulated by defining convergence via convergence of evaluations of functions by $\Pi^{\rm web}$ from Section~\ref{sss.sto} for the martingale problem on $\A^\ast$.} 
\end{remark}

\begin{proposition}[State spaces]
\label{prop.spaces}
The spaces $\A$ and $\A^\ast$ are Polish spaces.
\end{proposition}

The elements of $\A,\A^\ast$ are subsets of the Polish space $\M^{\rm tree}$ of tree-like metric measure spaces, and must be shown to be closed in the given topology, meaning that the additional restrictive properties are preserved under isomorphies.

We use the above concepts for populations with locations in $\Omega_N^L$ or $\Omega_N^\infty$ and types in $\K$, working with a $V$-marked, i.e., $\U^V$-valued process $\CU$ with $V=\Omega_N^\infty \times \K$. We replace in the definitions above $U,X$ and $\mu$ by
\begin{align}
\label{e1383}
U \times V, \text{ respectively, } X \times V,\\
\nu \text{ the resampling measure on } X \times V, \text{ respectively, } U \times V. \nonumber
\end{align}
The state spaces $\A,\A^\ast$ become $\A^V$ respectively $\A^{\ast,V}$, which are subsets of \emph{$\M^{V,\rm tree}$}, the $V$-marked metric fibred measure spaces that are \emph{tree-like and grounded}.


\subsubsection{Formal definitions $2$: state spaces for web/coupling, test functions, operators}
\label{sss.sto}

We define now an important class of test functions useful as domain for the operators we specify.

\begin{definition}[Polynomials on $\A^{\ast,V}$]\label{def.991}
On an element $\mfW \in \A^\ast$, we define a \emph{monomial} $\Phi$ of degree $\bar k=\sum_{i = 1}^m k_i$ on $\A^{\ast,V}$ for elements with finite time tree top $T \in [0,\infty)$. Consider time-points and functions
\begin{equation}
\label{e999}
0=t_0 < t_1 < \ldots< t_m = T \; , \; \chi = f \, g \text{  with  } f \in C_b(\K), \;  g \in C^c_b(\Omega_N^\infty), \; 
\varphi \in C_c^b(\R^+)^{n \choose 2}
\end{equation}
and set for $w=(X \times V,r,\varrho,\leq,(\nu^\ast_t)_{t \in [0,T]})$ and $X=[0,T] \times (\bigcup\limits_{t \in [0,T]} U_t)$ a
\begin{equation}
\label{a1435}
\begin{aligned}
&\text{monomial $\Phi$ depending on parameters}\\ 
&(T, \bar k,m,(k_i)_{i=1,\ldots,m}) \in \bigcup\limits_{\ell \in \N}\R^+ \times \N \times \N \times \N^\ell, \,(\varphi,f,g),
\end{aligned}
\end{equation}
all not exhibited in the notation, with
\begin{align}
\label{e1006}
\Phi (\mfw)&= \int_{(X \times \K \times \Omega_N^\infty)^m} \mathop{\bigotimes}\limits_{i=0}^m
\left(\nu^\ast_{t_i}\right)^{\otimes k_i} (\dd w_i)(\varphi_i \chi_i)(w_i)&\\
&= \int_{U_{t_0}^{k_0}} \; \int_{U_{t_m}^{k_m}} \int_{(\K \times \Omega_N^\infty)^{\bar k}} \mathop{\bigotimes}_{i=0}^m
\left(\nu^\ast_{t_i}\right)^{\otimes k_i} \left(\dd(\underline{u}_{k_i},\underline{v}_{k_i},\underline{\xi}_{k_i})\right) \varphi_i (\underline{\underline{r}}_{k_i}) f_i(\underline{v}_{k_i})g_i(\underline{\xi}_{k_i}),& \nonumber \\
&\qquad \qquad \qquad \text{ with } \underline{u}_{k_i}=(u_{i,1},\ldots,u_{i,k_i}), \text{ and similarly } 
\underline{v}_{k_i},\underline{\xi}_{k_i}. \nonumber
\end{align}
$(u_{k_i,j}, v_{k_i,j},\xi_{k_i,j}) \in X \times \K \times \Omega_N^\infty$ are coordinates of the $k_i$-sample of path at times $t_i$ with $i=1,\ldots,m$, $j=1,\ldots,k_i$. The set of all \emph{linear combinations} of monomials forms the set $\Pi^{\rm web}$ of \emph{polynomials}. 
\end{definition}

\begin{remark}
\label{r.1396}
{\rm The set $\Pi^{\rm web}$ is state-determining and the expectations are law-determining, allowing analysis of the $\A^\ast$-valued stochastic processes.}
\end{remark}

On the monomial in~\eqref{e1006} we build the \emph{duality function} $H^{\rm web}:\A^{\ast,V} \times  \A^{\ast,V}_{\rm fin } \to \R$, with the subscript $fin$ indicating we take only elements built on a finite number of initial dual individuals from which paths start or continue at a finite number of times.

Take $\CU=(U \times V,r,\kappa,\nu) \in \A^{\ast,V}$, $\mfy=(\Pi, \uur^\Pi,\underline{y}) \in \A^{\ast,V}_{\rm fin}$, and define
\begin{equation}
\label{e1451}
H^{\rm web} \left(\mfU,\mfy\right)=\int_{(X \times \R \times \Omega_N^\infty)^m} 
\mathop{\bigotimes}^m_{i=0} \left(\nu_{t_i}\right)^{\oplus k_i} (\dd w_i) \left(\varphi_i \chi_i\right)(w_i).
\end{equation}
On $\Pi^{\rm web}$ we define the operator $\mfL^{\rm web}$. Only a change occurs in the ``tree-top'' as time $T$ proceeds: resampling, migration-reshuffling, and aging act on the population alive at the current time. Precisely, as time $T$ proceeds, only the last factors $i=m$ corresponding to the tree-top change; the others are untouched. This change is driven by the operator of the martingale problem characterisation of $\mfU$ at time $T$, acting at the ``top'', i.e., the $i=m$ monomial factor; see equations~\eqref{pol0} and the sequel.

\begin{definition}[Operator $\mfL^{\rm web}$]
\label{def.995}
Write $\Phi$ as $\Phi_{\rm foss}\, \Phi_{\rm top}$ and set
\begin{equation}
\mfL^{\rm web} \Phi = \Phi_{\rm foss} \mfL_{\U^V}^{\Omega_N^\infty} (\Phi_{\rm top}), \text{ (see~\eqref{pol2}}). \label{e1431}
\end{equation}
\end{definition}


\subsection{How to use duality to construct processes}
\label{ss.tools}

We construct the process with geographic space $\Omega^{\infty}_N$ using a method developed in~\cite{DGP23}, used for the measure-valued process $\wh X^\infty$ on $\Omega_N^\infty \times \K$ in~\cite{gh1} and for the present marked genealogy-valued processes $\mfU$, $\mfW$, and $\mfX^\dagger$.

\medskip\noindent
{\bf Duality.}
We use the duality relation obtained for the processes on $\Omega_N^\infty$ by approximation by processes on $\Omega_N^L$. For each $L$, there is a known duality function $H^L(\cdot,\cdot)$ and dual processes (certain coalescents), and the $\CM_{\rm loc}(\Omega_N^L \times \K)$-valued Markov process $\wh X^L$ on $\Omega_N^L$ can be constructed as a solution to a well-posed martingale problem \emph{via duality}. In~\cite{DGP23} it is shown that one can often use duality and the \emph{existence} of the dual process to construct the original process. We formulate some facts about generic processes $X$ on a Polish state space $E$.

\begin{remark}
\label{r.1248}
{\rm  Certain class of Markov processes can, according to~\cite{DGP23}, be constructed from a specified dual process, namely, {\em standard processes} satisfying
\begin{itemize}
\item the state space $E$ of $X$ is a locally compact Polish space,
\item right-continuous filtration,
\item c\`adl\`ag paths,
\item strong Markov property,
\item quasi-left continuity ($X_{T_n} \la X_T$ for $T_n \uparrow T$, with $(T_n)_{n \in \N}$ stopping times).
\end{itemize}
We can apply this to our measure-valued context with $\wh X^\infty$ \emph{immediately}. For the process $\mfU$ (not locally or $\sigma$-compact!), the question is whether we can remove the \emph{locally compact} property of $E$ by some tightness property of the \emph{orbits of the stochastic process}, to establish the existence of a process with c\`adl\`ag paths. This issue is also addressed in~\cite{DGP23}.}
\end{remark}

The same strategy can be used to obtain the continuum-space $\U$-valued process, if in Blumenthal and Getoor 1968 \cite{BG68} we can replace in Theorem T9.4 the condition that the state space is \emph{locally compact} with a \emph{compactification} of a subset containing the orbit of the stochastic process, as also outlined in~\cite{DGP23}. We apply this now for the construction of $\mfU$ (and later in Section~\ref{s.proofs} for $\mfW$ and $\mfX^\dagger$). Such a construction can be carried out by the results in~\cite{DGP23} as follows. If the convergence of the dual process $\wh \mfY^L$ or $\wh \mfW^L$ for $\mfU^L$ respectively $\mfW^L$ can be established, converging to a limit dual Markov process $\wh \mfY$ respectively $\wh \mfW$, and if a \emph{limit duality function $H(\cdot,\cdot)$} can be found in each case such that it is dual to a \emph{limit process} $\mfU$, respectively, $\mfW$ or $\wh\mfW$, then $\mfU$ and $\mfW$ are \emph{Markov} processes with transition kernel $(P_t)_{t \geq 0}$ defined by
\be{ag16-1}
(P_t H(\cdot, X^\ast))(X) = (P^\ast_t H(X, \cdot))(X^\ast) \mbox{ for } (X,X^\ast) \in E \times E^\ast,
\ee
where $E$ and $E^\ast$ are the Polish state spaces of the process $X$ and the dual process $X^\ast$. Then $(P_t)_{t \geq 0}$ is uniquely determined and satisfies the \emph{Kolmogorov equations}. This gives the \emph{finite-dimensional distributions} of our continuum-space limit process, \emph{if we can check the assumptions in the quoted paper.}

The necessary statements of convergence of the dual processes $\wh \mfY^L$, respectively, $\wh \mfW^L,\mfX^{\dagger,L}$ and the existence of a dual $\wh \mfW$ are addressed in Section~\ref{s.prepa} and proved in Section~\ref{ss.pexisten}. What remains still to be established are:
\begin{itemize}
\item paths in $D([0,\infty), E)$,
\item characterisation by a martingale problem,
\item compactness properties.
\end{itemize}
The first point follows from the general theory of Hunt processes. The second and third points are verified in Section~\ref{ss.pexisten} via the \emph{compactification procedure} of Remark~\ref{r.1248}. This yields a limiting strong Markov process $\mfU,\mfW,\mfX^\dagger$.

To characterise the limit process by a \emph{well-posed martingale problem}, we specify the \emph{linear operators} and show that they satisfy the \emph{general criterion} for duality:
\be{ag17}
(G_X H(\cdot, X^\ast))(X) = (G_{X^\ast} H(X, \cdot))(X^\ast), \quad (X, X^\ast) \in E\times E^\ast.
\ee
The martingale problem is then \emph{well-posed}. This must be deduced from the $L$-approximation and the fact that relation~\eqref{ag17} holds for all $L\in \N$. The preservation of~\eqref{ag17} in the $L \to \infty$ limit is verified in Section~\ref{s.proofs}.

It therefore remains to verify the compactification property of~\cite{DGP23} and to write out the operators for $\mfU$ and $\wh \mfY$ respectively $\mfW$ and $\wh \mfW$ to get~\eqref{ag17}. This is all done in Sections~\ref{s.prepa} and~\ref{ss.pexisten}.


\subsection{Construction of the L\'evy-web and graphical representation: results}
\label{s.prepa}

The \emph{central tool} for the proofs of all our statements concerning the \emph{ancestral webs} $\mfW$, the \emph{ancestral space} $\mfX$, and its fibred version $\mfX^\dagger$ is the \emph{dual web} $\wh \mfW$ (and, on the way, $\wh \mfY^\infty$, the dual of $\CU$). We construct $\wh \mfW$ next, before the two proof sections, since the insight based on the duality is of independent interest, beyond the technical issues needed for the proofs. In particular we give properties of $\wh \mfW$ here, before Chapter~\ref{s.proofs} with the actual proofs of all the theorems.

For interacting genealogy-valued Fleming-Viot processes on $\Z^1$, one can give an \emph{explicit graphical construction} of the continuum-space limit based on a strong representation in terms of the Brownian web and the \emph{dual Brownian web}, both constructed in \emph{one random graph}. This random graph spans the information of the ``time-space process'' associated with the underlying population process. This technique on $\R^1$ uses the continuous path property of Brownian motion and the linear order structure of the reals. Despite the lack of these structural properties in our setting, we seek an explicit graphical representation also for our model, following the line laid out in Section~\ref{ss.levyweb} and based on the concept of the forward \emph{ancestral web} (\cite{gklimw},~\cite{Lambda}) and its dual: the backward process called the \emph{L\'evy-web}. Their joint construction in one ancestral space $\mfX^\dagger$ requires having the $\wh \mfW$ processes simultaneously for all time horizons.


\subsubsection{Main result for dual processes driven by spatial coalescents: $\wh \mfY,\wh \mfW,\mfX^\dagger$}
\label{sss.task}

We want not only a representation of expectations of functionals of the state via the dual starting with \emph{finitely} many individuals and a \emph{function-valued} duality, but want to \emph{represent} the state of $\mfW$ at time $t$ via $\wh \mfW$ directly, based on a suitable probability space. Note that with $\wh \mfW$ we have processes in backward time $s \in [0,t]$ for a time horizon $t \in [0,T]$, with $T$ the current time horizon of the population of all fossils. For fixed $t$ we have entrance laws, hence $\wh \mfW=(\wh \mfW^t)_{t \geq 0}=\{(\wh \mfW^t_s)_{s \in [0,t]}\colon\, t \in (0,\infty)\}$. This works via a spatial coalescent \emph{entrance law}, with a coalescent population starting from every point in the skeleton $\Omega_N^{\infty, \rm disc}$.

\medskip\noindent
{\bf The space-time process.}
We also want this for the \emph{space-time} process, so that we can track the evolution of local mono-type clusters in time. We therefore need a further extension to obtain $\mfX^\dagger$. We couple all the time-horizons for the $^+Y^L$ to obtain a consistent realisation $\{^+ Y^{L,t}, t >0\}$ for all times. This requires a \emph{discrete time approximation} of the set of possible time-horizons. For this, consider coalescing L\'evy-processes starting at every \emph{space-time} point in $D \times \Omega_N^{\infty,\rm disc}$ with $D$ the \emph{dyadic positive numbers}, respectively, coalescents with reshuffling and coalescence starting from every point in the skeleton $\Omega_N^{\infty, \rm disc}$ at different backward time points from $s \in [0,t]$ for time horizons $t \in D$. From this object we construct a process starting in every point in the time-space set of all individuals ever alive by completion.

In Section~\ref{s.proofs} we prove the claimed theorems. The starting point is to give the main tool in the construction of the \emph{dual objects} $\wh Y^\ast$ (the dual of the $\CM(\Omega_N^\infty \times \K)$-valued entrance law of $Y^\infty$), $\wh \mfY^\ast$ (the dual $\U^V$-valued genealogy entrance law associated with $Y^\infty$), and $\wh \mfW$ (the $\A^\ast$-valued entrance law associated with $Y^\infty$). This is based on knowledge about a \emph{finite} spatial coalescent $Y^\infty$ on $\Omega_N^\infty$, enriched with \emph{historical and genealogical} information to get a $\A^\ast$-valued process.

We state this rigorously in Proposition~\ref{prop.existen}. The key technical tool is Proposition~\ref{prop.geneal}. We first formulate and prove results on the convergence of $\wh Y^{\ast,L}, \wh \mfY^{\ast,L}$ and $\wh \mfW^L$ as $L \to \infty$.

\begin{proposition}[Existence entrance laws of $\wh Y^{\ast,\infty}, \wh \mfY^{\ast,\infty},\wh \mfW^\infty$, in the ancestral space $\mfX^{\dagger,\infty}$]\label{prop.existen}
$\mbox{}$
\begin{itemize}
\item[(a)] 
The entrance laws of the coalescent $Y^L$ starting with countably infinite occupation on each site in $\Omega_N^L$, extended to the time-space version of $Y^\infty$ respectively $\wh Y^\infty$, denoted $^+ Y^\infty, ^+ \wh Y^\infty$, exist. Furthermore, the entrance laws with types exist:
\begin{equation}
\label{e1380}
\text{$\wh Y^{\ast,L}, \wh \mfY^{\ast,L}, \wh \mfW^L$ on $\CM_{\rm loc}(\Omega_N^L \times \K),\, 
\U^V,\,\A^{\ast,V}$ and $\mfX^{\dagger,L}$ from $\M^{V,\rm tree}_{\rm fib}$.}
\end{equation}
Construct and couple all $Y^L, ^+ Y^L$ and the limit $L \to \infty$ on one probability space, by embedding $\Omega_N^L$ in $\Omega_N^\infty$. In the latter start all processes with one single individual in $\Omega_N^{\infty,\rm disc}$ for each $L \in \N \cup \{\infty\}$ by declaring points in a $(-L)$-ball as equivalent, and use one collection of independent L\'evy-processes starting from every point in $[0,\infty) \times \Omega_N^{\infty, \rm disc}$ to construct the process with coalescence/block-coalescence, by using one set of reshuffling paths. Similarly, the results hold for $^+ Y^L, ^+ Y^\infty$. Furthermore, as the limit under indefinite refinement, the continuum-time, time-space entrance law from $D \times \Omega_N^\infty$, with $^+ Y^{L,+}$ constructed for all time-horizons $t \in D$ (the dyadic numbers) is obtained.
\item[(b)] 
The limit in law of each member of the quadruplet in~\eqref{e1380} exists as $L \to \infty$ in $\CM_{\rm loc}(\Omega_N^\infty \times \K), \U^V, \A^{\ast,V}$ and $\A^V$ and is denoted:
\begin{equation}
\label{e1485}
^+ \wh Y^{\ast,\infty,+}, \wh \mfY^{\ast,\infty}, \wh \mfW^\infty \text{  and  } \mfX^{\dagger,\infty}.
\end{equation}
The resulting quadruplet of processes, all constructed on one probability space, is denoted for $L \in \N \cup \{0\}$:
\begin{align}
\label{e1332}
\text{$(^+ \wh Y^{\ast,+}, \wh \mfY^{\ast, L}, \wh \mfW^L, \mfX^{\dagger,L})$.}
\end{align}
The equivalence class with respect to the appropriate isomorphies is denoted with $[]$.
\item[(c)] 
Each of the limits of the processes of equivalence classes $[]$ is a Markov and Feller process characterised by a well-posed martingale problem. This holds as well for the process of the equivalence class $[]$ of the quadruplet.
\item[(d)] 
The processes $\wh \mfY^L,\wh \mfW^L$ and $\wh \mfX^{\dagger,L}$ are functionals of the ``field (in the index of time-horizons) of time-space processes'' $^+ Y^{L,+}$ for $L \in \N \cup \{\infty\}$.
\end{itemize}
\end{proposition}

\noindent
In the processes $\wh \mfY^L,\wh \mfW^L, \wh \mfX^{\dagger,L}$ we drop the $\infty$ in the sequel.

We next explain in detail \emph{the ingredients in the above proposition}. We formulate and \emph{prove} a version of Proposition~\ref{prop.existen} for \emph{finite} initial configurations, namely Proposition~\ref{prop.geneal}. This establishes existence of the processes for entrance laws from \emph{countable initial states}, for the $L$-truncations on \emph{countable} geographic space. We also conclude convergence of the finitely started dual processes as $L \to \infty$. In the remaining part of this subsection we prepare a strategy to prove Proposition~\ref{prop.existen} in Section~\ref{ss.pexisten} by lifting this to countably infinite initial states and a continuum of time-horizons. This is quite a task.


\subsubsection{First convergence results for L\'evy-Web $\wh \mfW$ and the ancestral space $\mfX^\dagger$}
\label{sss.detail}

We want to construct a representation of the evolution on the continuum space $\Omega_N^\infty$, a $\A^\ast$-valued process $\mfW$, when $\wh X^\infty$ is \emph{stochastic}. We do this by constructing the dual $\wh \mfW$. Throughout this section we work under the assumption
\begin{equation}
\mbox{l.h.s.\ of \eqref{ag6} } = \infty.
\label{e1246}
\end{equation}
The key step is to construct the potential \emph{dual} $\wh \mfW$ to the \emph{ancestral web} $\mfW$ and $\mfX^\dagger$. We first explain the ideas and challenges.

\medskip\noindent
{\bf Program.}
We construct the \emph{L\'evy-web} $\wh \mfW$ \emph{for each time horizon $t >0$} via the \emph{entrance law} starting from locally \emph{countably many individuals} of a \emph{spatial coalescent} and its \emph{genealogy}. We place, initially and at the assigned backward times for activation for the time-space version, one individual at each of the countably many points in $\Omega_N^{\infty, \rm disc}$:
\begin{equation}
\label{e1141}
\Omega_N^{\infty, \rm disc}= \bigcup_{L \in \N_0} \wt \Omega_N^L,
\end{equation}
where $\wt \Omega_N^L$ is the representative of $\Omega_N^L$ \emph{embedded} in $\Omega_N^\infty$. This allows embedding all $L$-processes $^+ Y^L$ into a process on $\Omega_N^\infty$, for every $L$ and $t>0$, giving $^+ Y^{L,+}$.

We claim that $\mfW$ and $\wh \mfW$ can be constructed embedded in the \emph{grounded fibred metric measure space} $\mfX^\dagger$ associated with $\mfW$ and $\wh \mfW$, so that one can read off from a representative both the forward process and the backward process in \emph{one graph}: a \emph{grounded fibred metric measure space} with a path of measures on paths (the fibred structure), called $\mfX^\dagger$. This gives \emph{forward and backward dynamics} in \emph{strong duality}, i.e., their laws agree, but additionally they live in \emph{one} space $\mfX^\dagger$, giving a \emph{graphical} strong duality. The setup is detailed in four points below, and Section~\ref{p.count} formulates the corresponding lemmas. Section~\ref{ss.pexisten} then proves Proposition~\ref{prop.existen} based on Proposition~\ref{prop.geneal}.

\medskip\noindent
{\bf 1. Finite coalescent case.}
Consider first a \emph{finite} time-space coalescent on $[0,T] \times \Omega_N^\infty$. This is a Markov jump process with migration according to $(c_k)_{k \in \Z}$, single-site resampling at rate $d^\ast$ based on collision local times of pairs, and block-resampling with reshuffling at rates $\lambda_k$ and resampling measure $\Lambda_k$ for $k \in \Z$. This is always well defined as a strong Markov jump process with c\`adl\`ag paths and the Feller property. It can be typed with $\K$ so that it has paths defining $(\Omega_N^\infty \times \K)$-marked paths of descent in backward times $s \in [0,T]$.

Defining the \emph{distance} of two points as twice the backward time until they are in the same partition element, and $2T$ otherwise, and using \emph{counting measure} on the individuals, we get a \emph{finite} genealogy structure: a \emph{grounded fibred measure metric space} with \emph{ultrametric fibres}. Observing the whole backward path marked with location gives a population marked with $(\Omega_N^\infty \times \K)$. Taking this family of backward dynamics for time horizon $T$ for every $T > 0$ defines elements in $\A^\ast$ evolving forward in $T$.

From~\cite{gh1,GLW05} and~\cite{gklimw,Lambda} we have as a consequence (the last two references must be extended to include the mechanism of \emph{block} resampling for the genealogies):

\begin{proposition}[Existence of genealogical coalescents on countable geographic spaces]
\label{prop.geneal}
$\mbox{}$\\
Consider the $L$-truncations of a spatial coalescent on $\Omega_N^\infty$ as above, for a given choice of parameters.
\begin{itemize}
\item[(a)] The processes $\wh Y^{\ast,L}, ^+ \wh Y^{\ast,L}, \wh \mfY^{\ast,L},\wh \mfW^{\ast,L}$ and $\mfX^L$ based on \emph{finite} initial configurations and on the spatial coalescent on $\Omega_N^L$ embedded in $\Omega_N^\infty$, on $\CM(\Omega^L_N \times \K), \U^V, \A^{\ast,V}$ and $\A^V$ respectively, exist. They are strong Markov and Feller processes individually and as a quadruplet, and are all characterised by well-posed martingale problems. This holds for all time-horizons $t>0$.
\item[(b)] For \emph{finitely} started processes the limit $L \to \infty$ exists in the respective state spaces and is characterised by a well- posed martingale problem.
\item[(c)] If the effective difference migration path is recurrent, then the claims in (a) hold also for every time-horizon $t>0$ for the \emph{entrance laws} on $\Omega_N^L$, $L \in \N$, with countable occupation at each site of $\Omega_N^L$ embedded in $\Omega_N^\infty$, by starting with one individual in $\Omega_N^{\infty, \rm disc}$ and identifying all individuals in a $(-L)$-ball with the position in $\Omega_N^L$.
\end{itemize} 
\end{proposition}

The next task is to turn to the \emph{convergence} $L \to \infty$ of the entrance law of the spatial coalescent from locally \emph{countable} initial states.

\medskip\noindent
{\bf 2. Preparations and challenges: the entrance law on $\Omega_N^\infty$.}
To carry out all the above constructions for entrance laws of locally \emph{countably} many individuals on $\Omega_N^\infty$, we take a \emph{dense countable} subset. The construction uses an entrance law from countably many dual individuals: one at each point in $\Omega_N^L$, or more conveniently \emph{one} dual individual at a \emph{dense countable} subset like $\Omega_N^{\infty,\rm disc}$ embedding $\Omega_N^L$ in $\Omega_N^\infty$, with all processes placed in $\Omega_N^\infty$. To take the closure in $\Omega_N^\infty$ in a suitable sense, we develop a \emph{topology} on a \emph{space of measures on paths}.

For the subset we first make a specific but natural choice. Start with one dual individual at every site in the \emph{dense skeleton} $\Omega_N^{\infty, \rm disc}$ (the union of all embedded $\Omega_N^L$ over $L \in \N$). We then get $\wh \mfW$ by showing $\wh \mfW^L$ converges to an $\A^\ast$-valued random variable $\wh \mfW$ at \emph{every} time $t >0$, passing to the time-space footprint of the process of partition element frequencies. This could give in particular a \emph{representation of the state} $\mfU$ at a \emph{fixed} time $t$ by generating the \emph{population} of dual \emph{paths}. We showed such a representation already for $\wh X^\infty$ in~\cite{gh1} and called it \emph{strong} duality.

\medskip\noindent
{\bf Space-time coalescent entrance law.}
To lift the duality of $\mfU$ and $\wh \mfY$ to a strong dual representation of our continuum-space process $\mfW$, we define the \emph{entrance law} of the \emph{space-time coalescent starting from every point in $\R \times \Omega^{\infty}_N$} equipped with the measure $\ell \otimes \mfh$. We achieve this rigorously from a time-space coalescent entrance law with locally countably many individuals, by starting with one individual at each point in a countable dense time grid and taking a limit as the width goes to zero. The countable case itself for a fixed time is already a limit of finite coalescents.

Rigorously, we use the \emph{$L$-truncation approximation} and start in the countable dense
\begin{equation}
\label{e1156}
\text{skeleton $\mathfrak{D}=\Omega_N^{\infty, disc} \times \{k \; 2^{-n}\colon\, k,n \in \N\}$ }
\end{equation}
and then define a \emph{closure} of the orbits of the processes of measures on paths associated with the time-space populations of paths indexed by $\mfD$, with respect to the topology of $\A^\ast$. For each $L$ this gives the desired object as an \emph{entrance law}: as the limit of processes $\mfY^L$ starting with large \emph{finite} numbers of individuals.

The main challenge in passing from $\Omega_N$ as geographic space to $\Omega_N^\infty$ is that random walks are replaced by \emph{L\'evy-processes} with \emph{countably infinitely} many (small) \emph{jumps}. We also have two motions (migration and reshuffling) for which we must ensure the existence of the \emph{L\'evy-web}, a measure-valued coalescing process started with individuals at \emph{every space-time point}. This means starting the (measure-on-path)-valued coalescent in the \emph{Haar measure} $\mfh$ on $\Omega_N^\infty$ and Lebesgue measure $\ell$ in time, i.e., in $\ell \otimes \mfh$.

We must show that $\wh \mfW$ has the desired property of being the \emph{dual} process of the ancestral web. This is guaranteed by Theorem~\ref{th.600}, which we prove with the construction given in Section~\ref{ss.pexisten}.

\medskip\noindent
{\bf 3. Literature and beyond.}
For the continuum-space Fleming-Viot model in Euclidean space $\R^1$ (the regime of a point-recurrent difference process), the \emph{construction of random graphs} for the genealogies of the Fleming-Viot process is known as the \emph{Brownian web} and results in analogues of $\mfW$ and $\wh \mfW$ in one graph. In that case the construction uses the Brownian web and coalescing Brownian motions for the backward process. In~\cite{GSW} this was used to give a representation of the genealogy process on $\R^1$. \emph{No} martingale problem characterisation has been established there, and there are serious technical problems.

We want to do the same (and more) for our processes $\mfU$ and $\mfW$ on $\Omega_N^\infty$, in particular for \emph{jump} processes (as opposed to Brownian motions in $\R^1$), and obtain a \emph{L\'evy-web on $\Omega^{\infty}_N$} as the \emph{jump-process analogue} of the dual Brownian web. For this we pass from paths with values in $\Omega_N^\infty \times \K$ to \emph{marked ancestral lines} in the space of individuals. Ancestral lines are specific oriented geodesics in the fibred \emph{ancestral space} $\mfX^\dagger$, running from time-$t$ level to time-$0$ level and marked with locations and types in $\Omega_N^\infty \times \K$. This is the basic idea of the ancestral web, which is a \emph{population of lineages}. Mathematically, this is a process of \emph{measures on paths} with values in the \emph{ancestral space} $\mathfrak{X}$, with elements called \emph{ancestral paths marked with $\Omega_N^\infty \times \K$}. This can also be used to represent the complete genealogy of the population ever alive, thus representing both the process of genealogies $\mfU$ and the historical process $\CH$ embedded in one object: the \emph{ancestral web} $\mfW$.

\medskip\noindent
{\bf 4. Construction: obstacles, general principles, and plan of proof.}
To construct $\wh \mfW$ and $\wh \mfX^\dagger$, we return to the space-time coalescent starting from \emph{finitely} many space-time points, show that the $L$-truncated versions $\wh \mfW^L$ converge to a limit $\wh \mfW^\infty$ equal to $\wh \mfW$, and then use the \emph{projective property} to start from an everywhere-dense subset with countably many individuals, denoted $\mfD$, passing to the measure-valued state as the strong dual $\wh \mfW$ of $\mfW$.

The existence of $\wh \mfW$ is the key point to be proved. Two \emph{serious obstacles} arise:
\begin{itemize}
\item the \emph{uncountable} collection of time-horizons $t \in \R^+$, and
\item the \emph{initial state} for the coalescent starting at \emph{all} the time horizons $t$ simultaneously and combined starting in space at \emph{locally countably infinite} or even \emph{uncountably} many points.
\end{itemize}
Everything must be lifted from moment-measure convergence of measures on $\Omega_N^\infty \times \K$ to \emph{random} measures on \emph{paths} in backward time $s \in [0,t]$. These objects must be approximated by \emph{countable discrete} structures.

We use the finite spatial coalescent to derive the coalescent starting with one individual from a \emph{dense space-time set}, then construct from there the coalescent starting from \emph{every} space-time point. We use an arbitrary \emph{countable dense} subset $\mfD \subseteq \Omega^{\infty}_N \times \R$ (the most important example is $\mfD=\Omega_N^{\infty, \rm disc} \times \{k \; 2^{-n}, k \in \N, n \in \N\}$). We introduce an ultrametric as twice the coalescence time on the fibre and consider the closure of the associated object, called $[_\mfD \wh \mfW]$. We equip $_\mfD \wh \mfW$ with the collection of measures:
\be{ag68}
(\mu_t)_{t \geq 0}, \text{ with } \mu_t \text{ a measure on  } (U_t,r_t),
\ee
given by $\mfh$ if we use a copy of $\Omega_N^\infty$ representing the dual individuals on every $t$-fibre. These individuals at time $t$ are then evolved for backward time $s>0$ up to $t$, to obtain the population of \emph{backward paths} that form together with the forward object a \emph{pair} of \emph{graph} and \emph{dual graph}.

We then take the \emph{closure} of this object in the topology of $\A^\ast$ to obtain $[\wh \mfW]$ and $\mfX^\dagger$. This closure of entrance laws leads to law-of-large-numbers effects for the occupation measures on small time intervals, so that with the continuity of the laws we obtain tightness properties as elements in $\A^{\ast,V},\A^V$.

\medskip\noindent
{\bf Generalities for the proofs.}
We begin with some \emph{general} issues. On $\Omega_N^L$ or $\Omega_N^\infty$, for fixed time $t>0$, $\wh \mfW_T$ is a $t$-fibred ($t \in [0,T]$) population of paths. For positive backward time it is atomic: a continuum of paths merges into one path at time $t-s$ for $s>0$. Unions of balls in $\Omega_N^\infty$ become one path from time $t-s$ onwards and perform a spatial coalescent with block-resampling back beyond $t-s$.

For the construction of $\mfX^\dagger$, the question is how to do this construction \emph{for all time-horizons $t$ strictly positive simultaneously}. In the dual, we have a jump-process \emph{entrance law} that must be started in a continuum at the entrance from the back. The initial jump into $\A^V$ from the marked finest partition must be established by time approximation.

One can do this for a time-grid in a well-defined way, by merging successively backwards in time the coalescents at the times of the time grid for some interval $[0,T]$. The question is then how to \emph{refine the grid indefinitely} and get \emph{convergence}. For concreteness we work with
\begin{equation}
\label{e1185}
\{k 2^{-m}\colon\, k \in \N\} \times \Omega_N^L,\, m \in \N_0, L \in \N_0\text{, respectively, }
\{k 2^{-m}\colon\, k \in \N, m \in \N\} \times  \Omega_N^{\infty, \rm disc}.
\end{equation}
We use the embedding of $\Omega_N^L$ in $\Omega_N^\infty$ and consider a coalescent population starting with \emph{one} individual at each point of $\Omega^{\infty,\rm disc}_{N}$, associated with forward time $t>0$, evolving for backward time $s$ with $0<s<t$. The question is whether this \emph{exists} and what configuration we have after a strictly positive time.

Passing from finite initial states to measure-valued entrance laws $(\wh Y^{L}_t)_{t \geq 0}$, starting with \emph{countably} many individuals at \emph{every} point in $\Omega^{L}_N$ for $t > 0$: this is done in~\cite{gh1} for $\Omega_N^L, L \in \N \cup \{\infty \}$ all embedded in $\Omega_N^\infty$, by working with one particle in every point of $\Omega_N^{\infty, \rm disc}$ and declaring all points in $L$-balls equivalent for the $\Omega_N^L$-processes. Based on $\wh Y_t^{L}$ from this coupling construction, we proceed to the \emph{genealogy}-valued objects, passing first to $\U^V$-valued random variables with $V=\Omega_N^L \times \K$, then to the \emph{marked genealogy} entrance law process $\wh \mfY^L=(\wh \mfY_s^L)_{s \in [0,t]}$ for \emph{all} time-horizons $t >0$ to be able to get $\mfX^\dagger$.

A first approach is to carry out the time approximation limit $m \to \infty$ to obtain the L\'evy-web $^\infty \wh \mfW^L$, then take $L \to \infty$ to obtain $\wh \mfW$ and $\mfX^\dagger$. The key step is to prove that $^\infty \wh \mfW^L$ can start at \emph{all} times $t>0$. Take $_m \wh \mfW^L$, defined as the space-time dual, then take the limit $m \to \infty$ to get entrance laws $\wh \mfW^L$ as the corresponding $\A^{\ast,V}$-valued process. A different way to generate this limit is to work with \emph{coalescents with rebirth} (introduced in the paragraph preceding~\eqref{e1043}), generating the ancestral path backwards by a second \emph{Markovian evolution mechanism}. At every coalescence time, an equal number of dual individuals start at that coalescence time backwards according to the rules; for Fleming-Viot resampling the absorbed individual is complemented by a \emph{newly created individual} at the \emph{death time}.

\medskip\noindent
{\bf Three steps to complete the construction.}
We need to carry out three steps:
\begin{itemize}
\item[(I)]
\emph{Existence and properties} of time-space entrance laws $^+ \wh \mfY^L$ and the rebirth version $^{\infty,+} \wh \mfY^L_{\rm reb}$. The first is simpler: it is the entrance law of the time-space coalescent for the $m$-grid, then taking $m \to \infty$, which essentially follows from the Markov property of $\wh \mfY^L$ and the recurrence (which for $L \in \N$ is point recurrent) given in~\cite{gklimw}. New here are the process $^+ \wh \mfY^L$ and the versions with types $^+ \wh \mfY^{\ast,L}$, as well as the $n \to \infty$ limits $^\infty Y^L$ of the entrance law from countable states. Here we extend~\cite{gklimw} to the case \emph{with} block-resampling.
\item[(II)]
With these processes $^+ \wh \mfY^{\ast,L}_{\rm reb}$ we construct $^+ \wh \mfW^L$ and $\mfX^L$ by taking the \emph{infinite refinement} limit $m \to \infty$ of the time-axis, using~\cite{Lambda} but adding the possibility of block-resampling and evolution with rebirth.
\item[(III)]
The core point is the tightness of $\mfW^L$ as $L \to \infty$ in $\A^{\ast,V}$.
\end{itemize}


\section{L\'evy-web: proof of two key propositions}
\label{ss.pexisten}

We now prove the results on $\wh \mfW$. We begin in Section~\ref{p.count} with preparations in two steps, then develop the formulations of the lemmas needed to complete the existence proof of $\wh \mfW$, $\mfX$, and $\mfX^\dagger$ in this Section~\ref{ss.pexisten}. We then proceed in four steps to construct $\wh \mfW$ in Sections~\ref{sss.constr}--\ref{sss.levy} and prove all the lemmas used in Section~\ref{ss.prftech}.


\subsection{Preparations: discrete and continuum}
\label{p.count}

We focus now on the dual. We construct $\wh \mfW^L$ and then prove convergence as $L \to \infty$ based on existing results for processes on countable Abelian groups as geographic space.

\medskip\noindent
{\bf Initial configurations for $Y^\infty$.}
We focus on the dual dynamics, in particular the spatial coalescent \emph{entrance laws} from starting the processes in a special configuration. We use a particular but natural choice, namely~\eqref{e957}--~\eqref{ag37b3}, since it is effective to consider coalescents with \emph{more structure} allowing a coupling for different $L \in \N \cup \{\infty\}$.

For the countable discrete case with $L \in \N_0$, we consider the embedded $\Omega_N^L$ as the subset
\begin{equation}
\label{e957}
\Omega^{\infty,L}_N = \{\xi \in \Omega^{\infty}_N \colon\, \xi_i =0\,\,\forall\, i \leq -L\},
\end{equation}
we denote the restriction to the subset
\begin{equation}
\label{e547}
\bigg\lbrace (\xi_i)_{i \geq -L}\colon\, \xi_i \in \{0,1,\ldots,N-1\}, \sum_{i \in \N_0} \xi_i < \infty \bigg\rbrace 
\end{equation}
by $\wt \Omega_N^L$, which is a \emph{countable Abelian group} and the \emph{embedded version} in $\Omega_N^\infty$.

Similarly, consider the countable dense subgroup of $\Omega^{\infty}_N$:
\be{ag37b42}
\Omega^{\infty,\rm disc}_N = \bigcup_{k\in\N_0} \wt \Omega^{\infty,k}_N
\ee
and the dyadic time-space grids:
\be{ag37b4}
D_L = \{2^{-L} \Z\} \times \wt \Omega^{\infty,L}_N,
\ee
which give successive refinements for $L \uparrow$, dense in $[0,\infty) \times \Omega_N^\infty$.

Each $D_L$ has the lattice property (points have distance at least $2^{-L}$ in each direction), and for the approximating \emph{countable dense} set we choose:
\be{ag37b3}
\mfD = \bigcup_{L \in \N_0} D_L,
\ee
which is a \emph{special} but natural skeleton of $\R^+ \times \Omega_N^\infty$, and is a \emph{topological subgroup} of $\R \times \Omega^{\infty}_N$. For $\mfD$ we start with \emph{one} individual per site; for $D_L$ we start with countably many per time point. By embedding the $D_L$ in $\mfD$ we can compare the dynamics with these two starting sets (one for the process on $\Omega_N^L$, $L \in \N$, and one for $L=\infty$), all on the geographic space $\Omega_N^\infty$. This arises by aggregating all points in the $(L,L)$-ball of $\mfD$ into one element. In particular, as we saw, we can couple all these evolutions in one probability space by constructing them in one geographic space with one collection of L\'evy-processes and one set of block-resampling times.

For $D_L$ and the coalescent starting from every space-time point with \emph{countably} many individuals, we use the embedded version and construct for all $L \in \N$ on one probability space, starting with one point at each time in a point of $\Omega_N^{\infty, \rm disc}$ and declaring points in hierarchical distance $L$ as equivalent. This is an \emph{infinite particle system}. For $t>0$ we need a nice Markov process and entrance law. We show that the process is in a \emph{Liggett-Spitzer space} in strictly positive time, on which the dynamics is well defined with locally finite states.

The Liggett-Spitzer space was introduced in the theory of \emph{infinite particle systems} as standard state space in~\cite{LS81}. Applying this to the spatial Kingman coalescent with migration, we have an entrance law from countably many individuals per site evolving in a Liggett-Spitzer space. This was proved in \cite{GLW05} for the case of the spatial Kingman coalescent on a \emph{discrete} countable Abelian group, applicable whenever $d \in (0,\infty]$ or $\Lambda_{-\infty}$ comes down from infinity. A similar argument works for the case where $\Lambda_0$ respectively $\Lambda_{-L}$ is \emph{dust-free}.

\medskip\noindent
\underline{$d^\ast=0$.}
On $\Omega_N^\infty$ for $d^\ast=0$ we have only Cannings block-resampling. To guarantee a locally finite system via block-resampling under the \emph{stochasticity condition}, we give an argument based on the assumptions on $\Lambda_k$ for $k\geq -L$.

\medskip\noindent
\underline{$d^\ast>0$.}
On $\Omega^{\infty}_N$ we have \emph{L\'evy-processes} with countably many small jumps on positive time intervals, rather than random walks as on $\Omega^L_N$. The argument for random walks does \emph{not} directly work because L\'evy-processes never stand still. However, we have a \emph{collision local time} in the \emph{point-recurrent} regime: the existence of a local time for the difference process of a pair. Taking $L \to \infty$ gives delayed coalescence with rate $d^\ast$, leading again to coming down from infinity if $d^\ast >0$.

In both cases for $d^\ast$ this technical problem can be handled with the help of the \emph{truncation approximation} as follows. We showed in~\cite{gh1} the following two properties.

\begin{lemma}[Entrance law: Point-recurrent]
\label{l.1452}
For coalescing L\'evy-processes on $\Omega_N^\infty$ with a point-recurrent effective difference L\'evy-process, the entrance law of the $L^\prime$-truncated dynamics starting with countably many particles at each point in $\Omega_N^L$, $L \in \N$, exists and the configuration of particles at positive times is in a Liggett-Spitzer space for every $L^\prime >L$. This extends to time-space coalescents.
\end{lemma}

\begin{lemma}[Entrance law: Harris-recurrent but not point-recurrent]
\label{l987}
For coalescing L\'evy-processes on $\Omega_N^\infty$ with Harris-recurrent effective difference L\'evy-processes, the entrance law of the $L^\prime$-truncated dynamics starting with countably many points at each point of $\Omega_N^L$ exists and the configurations of points are for positive time in the Liggett-Spitzer space for $L^\prime >L$. This extends to time-space coalescents.
\end{lemma}

\begin{corollary}[Finite $L$]
\label{cor.finite}
In the recurrent difference L\'evy-processes case, the finite-$L$ truncations are locally finite after a possible $\varepsilon$-truncation of the population for all $\varepsilon >0$ for $t >0$.
\end{corollary}

As a consequence of these two lemmas, we pass to the measure-valued dynamics for the dual process $\wh Y^\infty$ to get an entrance law for the $\Omega_N^\infty$-dynamics in the limits $L^\prime \to \infty, L \to \infty$, which starts in the Haar measure $\mfh$ in these limits.

\medskip\noindent
{\bf Space-time coalescent for $D_L$.}
We define the \emph{space-time coalescent for each $D_L$} with $D_L$ as the set of initial positions for a partition element. To get to the union of all $D_L$ and to obtain a \emph{locally countably infinite configuration}, we successively add one additional element (populations starting at the new points), using the projective property to get consistency.

We use a $(\Omega^{\infty,\rm disc}_N)$-indexed collection of independent L\'evy-processes to construct the coalescent. At coalescence events only one of the involved L\'evy-processes is used further on. We construct the coalescents for different $D_L$ from the \emph{same collection}, and after coalescence keep using the path with the lower value of $L$ among all those involved. This means that given the \emph{collection of all L\'evy-processes}, the $D_L$ are embedded into each other. Furthermore, refining the grid thins only possible coalescence events with the coalescent for the previous $L$. This gives a coupling of the processes for all $L \in \N \cup \{\infty \}$.

So far we have considered the coalescent on $\Omega^{\infty}_N$ starting from one individual at each point in $\mfD$ (the union of $D_k$ over $k$) in $\Omega_N^{\infty, \rm disc}$, respectively, in time-space in $D_k$. We now need to start with one individual in \emph{every time-space point}, i.e., with one point in the \emph{closure of $\mfD$} as a subset of $\R^+ \times \Omega^{\infty}_N$. Since the process is defined on a \emph{dense skeleton}, we pass to the \emph{closure of this population of paths} in the proper topology to get the object corresponding to starting the entrance law of the measure-valued version of $Y^\infty$ in the $\ell \otimes \mfh$-measure.

We pass to the \emph{entrance law} starting in the \emph{closure of} $\mfD=\bigcup_{k \in \Z} D_k$. Rigorously, we start the measure-valued version of the coalescent in the \emph{Haar measure $\mfh$} on $\Omega_N^\infty$ and get partition elements with locally non-trivial frequency (the no-dust property). Details are given in~\cite{gh1}.

\medskip\noindent
{\bf Tightness in $\A^\ast$.}
What is still open is to show \emph{tightness} of $\mfY^L,\mfY^\infty$, respectively $\wh \mfW^L,\wh \mfW$, in the stronger topologies $\A^{\ast,V},\A^V$ including the genealogies in the \emph{topology of $\U,\A^\ast$}. We view the basic population of these objects as a \emph{measure} on a \emph{field of paths} starting in points of $(\R^+ \times \Omega_N^\infty)$ with values in $\Omega_N^\infty$, with time running from $t$ backwards to $0$. This issue is addressed in the next subsubsection: \emph{this amounts to getting a L\'evy-web on $\Omega_N^\infty$ and embedding it in the ancestral space $\mfX^\dagger$.}

In the ancestral space $\mfX^\dagger$ constructed for $\wh \mfW$, we want to embed both $\mfW$ \emph{and} $\wh \mfW$. The forward geodesics from fibre $s \geq 0$ to fibre $t>s$ are the paths reversed in orientation in the backward geodesics from fibre $t$ to $s$. Embedding these \emph{populations of forward paths} and the \emph{population of backward paths} in a single tree-like metric fibred measure space --- the \emph{ancestral space of all marked individuals ever alive} --- gives $\mfX^\dagger$. The $\dagger$ indicates that $\mfX$ is additionally equipped with \emph{fibred sampling measures}, turning the forward and backward paths into \emph{$\Omega_N^\infty \times \K$-marked ancestral paths}. This object of forward and backward paths creates a \emph{web}: a ``random web'' of forward \emph{lineages} and, looked backward, \emph{paths of descent}.

This would mean that $\mfW$ and $\wh \mfW$ are \emph{paired} in the sense that they are populations of forward geodesics, and the \emph{reversed orientation} yields the population of backward geodesics. The measure $(\mu^\ast_t)_{t \geq 0}$ respectively $(\nu^\ast_t)_{t \geq 0}$ is associated with the forward population at times $t$, while the backward evolution is associated with the measures $\mfh$ on $\Omega_N^\infty$ at the backward time $s=0$ for all fibres for time horizon $t$. If evolved with the dual for time $t$, this gives $\mu^\ast_t$ respectively $\nu^\ast_t$. The reversion operation acts for all pairings $[s,t]$, with $s \in [0,t]$ and $t \in [0,T]$. In this way, in the \emph{graph} of forward and backward paths in $\mfX_T^\dagger$, the processes $\mfW_T$ and the dual $\wh \mfW_T$ are \emph{pathwise dual}.

We produce and study $\wh \mfW,\mfX^\dagger$ next by the \emph{truncation limit} $L \to \infty$ first applied to the dual \emph{L\'evy-web}. This is used to establish the $L \to \infty$ truncation limit for the \emph{forward} process. The forward process $\mfW$ is then constructed via $\wh \mfW$, and both are realised on one probability space in $\mfX^\dagger$. This happens in Section~\ref{s.proofs}. Here we focus on the construction of $\wh \mfW$ and $\mfX^\dagger$.


\subsection{Construction of $\wh \mfW$ and $\mfX^\dagger$ part 1: finite populations of paths}
\label{sss.constr}

To get the scenario described above, we proceed in three steps. We start in this section with the construction of the L\'evy-web $\wh \mfW$ for the case of finitely many dual individuals. We first obtain a representative $\overset{\sim}{\wh \mfW}$ of the \emph{marked genealogy} of the backward \emph{space-time coalescent}. This requires moving from $\wh Y^{\ast,\infty}$ to $\wh \mfY^\ast$, and then incorporating \emph{all} time horizons $t>0$ \emph{simultaneously}, observing the \emph{path} with time $s \in [0,t]$ as the basic element, to pass to $\mfX^\dagger$.

\medskip\noindent
{\bf Finite space-time coalescent.}
Take the \emph{space-time} coalescent on $\Omega_N^\infty$. Consider first a \emph{finite set} of space-time points and \emph{enrich} it to an $\A^{\ast,V}$-valued object. Recall that two L\'evy-process paths running independently never remain equal for a positive time (which requires, for example, that $\lambda_k+c_k>0$ for all $k \leq k_0$ for some $k_0$). Therefore the ancestral paths are equal for positive time intervals only when we go beyond the coalescent time, with the point at the coalescence time being the most recent common ancestor of the coalescing partition elements.

In order to introduce a \emph{distance} between points in the basic set, we define first:
\be{ag65}
\delta (\zeta, \zeta^\prime) = (t_1 - \tau) + (t_2-\tau),
\ee
where $\tau = \sup \{t \colon \zeta(s) = \zeta^\prime (s)$, $s \leq t\}$. The distance $d(\cdot,\cdot)$ is then
\be{ag66}
d((\iota,s),(\iota^\prime, s^\prime)) = 1-\exp\big(-\delta (\zeta^{(\iota,s)},\zeta^{(\iota^\prime, s^\prime)})\big),
\ee
where $\zeta^{\iota,s}$ is the backward path beginning in $\iota \in \Omega^{\infty}_N$ at time $s \in \R$. This defines a \emph{metric space} equipped with \emph{marks}, where the mark is the \emph{backward path} of the partition element to which our point belongs till time $t$.

We get a natural \emph{fibred measure $\nu$} on this L\'evy-structure to pass to a $\U^V$-valued state: the initial point in a fibre gets mass one, and these masses merge upon coalescence (they are added). This structure can also be put on the discretisations $\Omega_N^L$ of $\Omega_N^\infty$ via $L$-truncation, obtaining $\wh \mfW^L$ for all $L \in \N$.

\medskip\noindent
{\bf Convergence of $\wh \mfW^L$.}
We must verify that as $L \to \infty$ the $\wh \mfW^L$ converge in distribution to an element $\wh \mfW^\infty$ in $\A^{\ast,V}$, and that the \emph{coupled versions} converge a.s.\ in $\A^{\ast,V}$. This limit we identify as the state of a coalescence process, specified by evolution rules for all finite samples. This is fairly straightforward from what was proved in~\cite{gh1} for $Y^L$ in \emph{path space} and the statement that the \emph{joint occupation times of pairs} converge to the \emph{collision local times of pairs}. This implies convergence of coalescence times, hence of the \emph{distances between pairs} and of the generated ultrametric measure spaces. Since the effective L\'evy-processes also converge, we get convergence in $\U^V$ and $\A^{\ast,V}$. We omit the standard details.

How to proceed to get $\mfX^\dagger$? For a suitable choice of $\kappa$ we want to obtain:
\be{ag67}
\wt \mfX^\dagger = (M,d,\kappa,\nu),
\ee
which is a $\Omega_N^\infty$-\emph{marked grounded metric fibred measure space}, such that on each level-$t$ the object $d$ is an \emph{ultrametric} and the sampling measure $\nu$ is fibred in $\nu_t$ for $t \in [0,T]$. The set $M$ represents \emph{all individuals ever alive} before time $T$. Using that $\Omega_N^{\infty,L}$ with countable occupation per site can be embedded in $\Omega_N^\infty$, we choose:
\be{ag63}
M = \{(\xi,t) \colon\, \xi \in \Omega^{\infty}_N, t \in \R\}
\ee
with points in $M$ \emph{marked with paths}, interpreted as the \emph{ancestral path} of the current population at that time:
\be{ag64}
\kappa \colon \Omega^{\infty}_N \times \R \la D^-.
\ee
Precisely:
\begin{equation}
\label{e1133}
V=\bigcup\limits_{T \in \R}
D \left((-\infty,T], \Omega^{\infty}_N \right)=D^-.
\end{equation}

We view the right-hand side as a subset of $D^\ast(\R, \Omega^{\infty}_N)$, consisting of Skorohod \emph{backward} paths (time running backward, indicated by $\ast$) with the special property that it has c\`adl\`ag paths constant for times $[T,\infty)$ for some $T \in \R$. Note that $D^\ast(\R, \Omega^{\infty}_N)$ is a \emph{complete separable metric space}. Allowing negative time incorporates individuals not having a most recent common ancestor in the present framework.

\medskip\noindent
{\bf Marked paths and the ancestral web.}
We also need the paths of descent to be \emph{marked with types}, since the lineages are. The types are assigned from the location at time $0$ and then carried forward. Hence we have a marked path through ancestral and geographic space. We take the set of \emph{individuals} $M=\Omega^{\infty}_N \times \R$ and specify a set of \emph{marks $V$}, chosen as a \emph{set of paths} with \emph{values in $\Omega_N^\infty \times \K$}. The measure is the \emph{counting measure}. We give $X$ the richer structure of a \emph{fibred measure space}, where $\nu$ is given via a process of measures $(\nu_t)_{t \geq 0}$. We then define via $\wh \mfW$ the ancestral space and set:
\begin{equation}\label{e1388}
\mbox{\emph{marked} metric measure space  } \mfX=[(M \times V, r_M,\nu)],
\end{equation}
marked with the \emph{ancestral paths} of the backward evolution.

Why is $\mfX$ marked with \emph{paths}? We view our primary object as a \emph{population of paths} in ancestral space marked at all their time indices $s \in [0,T]$ with type and location. It is instructive to look in the \emph{forward direction}. We want to view $\wt \mfW$, which is a Polish space via strong duality, as a representative of the \emph{marked lineages} representing the \emph{equivalence class} of the $D^\ast(\R, \Omega^{\infty}_N \times \K)$-marked metric measure spaces $\mfW^\infty$, i.e., $[\wt \mfW]=\wh \mfW^\infty$. The natural objects are then populations of marked paths connecting the presently alive individuals with the founder of their family, with the path-valued mark chosen here.

To get the \emph{forward-backward pair} $(\mfW,\wh \mfW)$: the $(\Omega_N^\infty \times \K)$-marked geodesics from $t$ to $0$ looked in reverse orientation give the \emph{forward evolution of the population of lineages}. The lineages are the geodesics connecting the individual to the ancestor at time zero. The state space is $\A^{\ast,V}$. We use a strengthening of the $V$-marked Gromov-weak topology on that space (see Appendix~\ref{appA}). We then show that this $[\wt \mfW]$ is actually $\mfW$.

We next pass, using the properties of the \emph{finite} versions of the dual process, to \emph{entrance laws} from \emph{locally countable} initial states, to prove Proposition~\ref{prop.existen}.


\subsection{Construction of $\wh \mfW$ and $\mfX^\dagger$ part 2: entrance laws for infinite populations of paths}\label{ss.levy}

Our strategy is to construct $\mfW$ and $\mfX^\dagger$ via the $L$-truncation limit. We verify the \emph{existence} and the \emph{properties} of the limit $L \to \infty$ of the objects for the $L$-truncation via the \emph{dual} processes $\wh \mfW^L$ and $\wh \mfW$, respectively $\mfX^{\dagger,L}$ and $\mfX^{\dagger,\infty}=\mfX^\dagger$. We construct $\wh \mfW$ as the $L \to \infty$ limit of $\wh \mfW^L$, and based on that construct $\mfX^\dagger$ via the $\mfX^{L,\dagger}$ as $L \to \infty$.

The main new steps are to pass for $\wh \mfW$ from the convergence of \emph{time-space marginal distributions} via a function-valued duality (stated in Proposition~\ref{prop.geneal}(a) and proved via~\cite{gh1}), now to the stronger result: \emph{convergence} of the \emph{entrance laws} from countable and countable-dense or even uncountable initial states. Similarly for $\mfX^\dagger$ in \emph{path space in the time-horizon} of the lineages, respectively, paths of descent for the \emph{entrance laws}.

Moving to a countable initial configuration for the dual process has two aspects. \textit{(i)} Since $\Omega_N^\infty$ (and $\Omega_N^L$) are \emph{infinite} geographic spaces, we have an infinite particle system for the spatial coalescent and in particular its entrance law. We work with approximations by processes on \emph{finite geographic subsets}. This approximation is well understood from the theory of infinite particle systems (see~\cite{LS81,GLW05}) and we do not work out details here. \textit{(ii)} The second aspect is the $L \to \infty$ limit, which is our focus in the sequel.

\medskip\noindent
{\bf Problems and tools.}
The main issues are tightness properties of the approximations by $L$-truncated versions. Return to the marked partition-valued coalescent and prepare tightness criteria in $\A^\ast$. We add the type from $\K$ to the backward partition-valued coalescent, assigned based on the \emph{initial state $\wh X_0^\infty$ of the $(\K \times \Omega_N^\infty)$-measure-valued functional} $\mfU_0$, giving a finite spatial coalescent with types: $Y^{\ast,L}, Y^{\ast,\infty}$. For the measure-valued form needed to pass to the entrance law from countable initial states this leads to $\wh Y^{\ast,L}$ and $\wh Y^{\ast,\infty}$. We lift this to $\wh \mfY^{\ast,L}$ (the dual entrance law with values in $\U^V$) by defining an ultrametric between partition elements as twice the time to end up in the same partition element, or $2t$.

These objects can be used to obtain elements in $\A^\ast$ when considering the population of ancestral paths. The \emph{dust-free} property of the finite coalescents, used to pass to the entrance law of $\wh \mfY^L$ and $\mfW^L$, corresponds to the \emph{tightness} of corresponding elements of $\U, \A$ or $\A^\ast$ (and later the marked objects) viewed as \emph{elements in $\A^\ast$}. The $Y^{\ast,L}$ being \emph{dust-free} gives the needed tightness. The next task is to show that
\begin{equation}
\label{e1362}
\text{\emph{tightness} properties of  $\big(Y_t^{\ast,L}\big)_{L \in \N}$ and the \emph{dust-free} property of $\big(\wh Y_t^{\ast,\infty}\big)$,}
\end{equation}
in the limit $L \to \infty$, and analogously for the enrichments $\wh \mfY^{\ast,L},\wh \mfW^{\ast,L}$, and $\mfX^{\ast,L}$.

\medskip\noindent
{\bf Infinite geographic space.}
In a finite geographic space, Kingman coalescents come down from infinity (a sort of dynamic \emph{compactification}). We have an \emph{infinite} geographic space. Since an individual on $\Omega_N^L$ can move arbitrarily far from its initial state with positive probability, we do \emph{not} come down from infinity anymore (see \cite{ABL12}). But, for $\delta > 0$, consider the elements that do not move before time $\delta$: these can give $1-\varepsilon$ of the frequency if $\delta = \delta(\varepsilon)$ is chosen properly, with finitely many partition elements starting from a specific site.

A first goal is to establish that the configuration of \emph{occupation numbers} is in the \emph{Liggett-Spitzer space} and has a \emph{stochastic upper bound} that is also in that space. A second goal is to count the number of partition elements with marks in a small ball around a point, or at least for every $\delta$-truncation as above. We start by proving that the occupation measure induced by the coalescent at a fixed time $t$ is \emph{locally finite uniformly in $k$}. Next: after a fixed \emph{positive} time, do small enough balls belong to \emph{one} partition element in the $L \to \infty$ limit? We want a \emph{countable cover} by such (mono-type) partition-element balls for a bounded closed set. This was resolved in~\cite{gh1} and we return to it below.


\subsection{Construction L\'evy-web and ancestral space part 3: main lemmas for entrance laws}
\label{sss.levy}

The argument proceeds in steps: we first formulate the key lemmas in this section and then give their proofs in Section~\ref{ss.prftech}.

\medskip\noindent
{\bf Step 1} (Formulation of the main lemmas yielding $\wh \mfW$.)

We give the basic steps for the \emph{construction of the dual web $\wh \mfW$ via a sequence of Lemmata} dealing with two themes: (1) the properties of the $\Omega_N^L$-dynamics needed, and (2) the issue of tightness as $L \to \infty$.

\medskip\noindent
{\bf (1) Entrance laws from $\infty$-state for $\wh Y^L, \wh \mfY^L$ and $\mfX^{\dagger,L}$ for fixed time, for $\Omega_N^L$.}

The first task is to get the information for $\Omega_N^L$, since we need to work with the \emph{$L$-approximation} of $\wh \mfW$ by $\wh \mfW^L$. For this we need \emph{geometric properties} of the configuration when passing from \emph{finite coalescents} to \emph{entrance laws} from ``infinite'' occupation. We consider extended initial states starting countably in $\Omega_N^L$ and couple them as before by starting with one at each point in the dyadic time-grid of orders $m$ and $\Omega_N^{\infty, \rm disc}$. We use left sub- and superscript to indicate initial conditions; beyond this section we drop the ``$\infty$'' symbols.

We establish the existence of \emph{entrance laws} $\wh \mfW^L$ by starting in one fibre locally with $n$-particles and letting $n \to \infty$ (denoted by a left upper $n$). This gives, for \emph{all} collections of finite fibres, an entrance law starting in time-space countably. We also construct $\wh \mfW^L$ and $\mfX^{L,\dagger}$ using the backward paths as geodesics, continuous in the time horizon in the state spaces $\U^V,\A^{\ast,V},\A^V$. With block-resampling, the results from~\cite{Lambda,gklimw} generalise to give:
\begin{lemma}[Enriched processes: $^\infty _\infty \wh \mfW^L$ and $\mfX^\dagger$]\label{l.1745}
\text{}\\
(a) The backward genealogy process $^n \wh \mfY^L$ and the backward process $^n {\overset{\leftarrow}{\mfW}}_T^L$ as well as the normalised version $^n \wh \mfW^L$ exist for every $T>0$.\\
(b) For every fixed $L \in \N$ there exists an entrance law:
\begin{equation}
\label{e1749}
\CL \big[^n \wh \mfY^L\big] \nto \CL \big[^\infty \wh\mfY^L\big], 
\qquad \CL \big[^n\wh \mfW^L\big] \nto \CL \big[^\infty \wh \mfW^L\big].
\end{equation}
(c) The entrance law $^\infty _\infty \wh \mfY^L,^\infty _\infty \wh \mfW^L$ exists and can be started for all dyadic $t$ simultaneously and defines a unique $\mfX^{\dagger,L}$, by specifying the geodesics between all level sets given by $^\infty _\infty \wh \mfW^L$ for all dyadic $t$. In other words, one may take the limit $m \to \infty$, with $m$ the degree of the dyadic refinement. 
\end{lemma}

It remains to prove that this holds still \emph{with} block-resampling added.

\medskip\noindent{\bf(2) Passing from $\Omega_N$ and $\Omega_N^L$ to $\Omega_N^\infty$ dynamics.}
We now take the truncation approximation $L \to \infty$. Our main task is to obtain the \emph{tightness for $L \to \infty$ of $^\infty  \wh \mfW^L$ in $\A^{\ast,V}$} and related $\mfX^{L,\dagger}$ in $\M^{\rm tree,V}_{\rm fib}$ for $L \to \infty$, then construct $\mfX^\dagger$ from $\mfX^{L,\dagger}$ as this limit in $\A$. We proceed in a number of Lemmata, first summarising what we need about the $\Omega_N^L$-processes.

In~\cite{gh1} we obtained the convergence of the measure-valued entrance law of the dual processes $\wh Y^{\ast,L}$ in finite-dimensional distribution to $\wh Y^\infty$. We now pass to $\wh \mfY$ in $\U^V$ and then to $\wh \mfW$ in the topology of $\A^\ast$ (which is stronger). We must obtain \emph{tightness} of the \emph{marked and fibred measure $\R$-trees} and of the \emph{ancestral webs} as populations of paths viewed as objects in $\A^\ast$.

\begin{lemma}[Tightness of $\{^\infty \wh \mfY^L, L \in \N\}$ in $\U_{\rm loc}$]
\label{l.1646}
$\mbox{}$\\
The $\{\CL \left[^\infty \wh \mfY^L\right], L \in \N\}$ are tight as laws on $\U_{\rm loc}$, similarly the rebirth versions. 
\end{lemma}

This lemma gives the \emph{entrance law} for \emph{fixed time horizon $t$} from the configuration starting with ``one initial individual at every point in $\Omega^{\infty}_N$''. However, for $^\infty \wh \mfW$ and especially $\mfX^\dagger$ we need more: we must consider the \emph{time-space coalescent} for \emph{all time horizons $t$ simultaneously}. We use $^+ \wh \mfY^L$ and choose as time horizons the dyadic numbers of order $m$: $\{^\infty _m \wh \mfY^L\}$, then take $m \to \infty$ to end up with $\{^\infty _\infty \wh Y^L, L \in \N \}$. We need to get \emph{processes of measures on continuous} (respectively, c\`adl\`ag) \emph{paths} as limit points. This is the key tightness lemma, strengthening Lemma~\ref{l.1646} for $^\infty _\infty \mfY^L$ as $L \to \infty$.

\begin{lemma}[Tightness of $\wh \mfW^L$ in $\A^{\ast,V}$, existence of $\mfX^\dagger$]
\label{l.1666}
The following hold:
\begin{itemize}
\item[(a)] 
The laws $\{\CL [^\infty _\infty \wh \mfW_T^L]\colon\,L\in \N\}$ are tight for every $T>0$ on the space $\A^{\ast,V}$.
\item[(b)] 
All $(\mfX_T)_{T \geq 0}$ can be embedded in a \emph{grounded $(\Omega_N^\infty \times \K)$-\emph{marked} metric fibred measure space} and \emph{$\R$-tree} $\mfX^\dagger$, generated by the entrance law $^\infty _\infty \wh \mfY$ by taking $^\infty _\infty \wh \mfY_T$ for $T \to \infty$.
\end{itemize}
\end{lemma}

With Lemmas~\ref{l.1745}--\ref{l.1666} we can then obtain the $L \to \infty$ limits of the $L$-truncations and construct the stochastic processes $\wh \mfW,\mfX^\dagger$.

\medskip\noindent
{\bf Step 2.} (Technical tools for proofs: Formulation of tightness criteria)

\medskip\noindent
{\bf From $\Omega_N^L$ to $\Omega_N^\infty$: Formulation of criteria to prove tightness properties.}
We have seen in (1) that we must verify (i) the \emph{tightness} of the processes of the lemmas, generated by \emph{time-space coalescents} or \emph{coalescents with rebirth}, in the index $L$, and continuity in time horizon uniformly in $L$. This means criteria for tightness and continuity, first in $\U^V$, then in $\A^{\ast,V}$ and $\A^V$.

The second point (ii) is that we deal with the measure on the whole ancestral \emph{path} of an individual. This is a function space, and we must verify that the process $^\infty _\infty \wh  \mfY^L$ has continuous (or if marked with types: c\`adl\`ag) paths in the limit and is continuous in the time horizon. The behaviour of the path in $\Omega_N^\infty \times \K$ follows \emph{immediately} from the results in~\cite{gh1} on $\wh X^L,\wh X^\infty$: the measures are atomic in $\K$ and locally monotype, i.e., atomic in $\Omega_N^\infty$, locally finite with continuous (and with block-resampling: c\`adl\`ag) paths in the \emph{weak atomic topology}. This finishes (ii).

\medskip\noindent
\textbf{(i)}
To address the first point we need a new quantity. Suppose that $d^\ast >0$ and that $\Lambda_{-\infty}$ comes down from infinity. Then
\begin{equation}
\label{e1282}
\begin{aligned}
^\delta N_t^{\varepsilon,L} (\eta) =&\text{ number of partition elements at backward time $\varepsilon >0$,}\\ 
&\text{ \emph{not} having moved until time $\delta>0$,}\\
&\text{ starting at backward time $0$ with countably many elements in every $\eta \in \Omega_N^L$,}\\
&\text{ and evaluated at backward time $t$.}
\end{aligned}
\end{equation}
We now formulate a \emph{tightness criterion} in $\A^{\ast,V}$ or $\A^V$.

We need a \emph{stochastic bound} on $^\delta N_t^{\varepsilon,L}(\eta)$ to prevent dust formation in the genealogy as $L \to \infty$. In the limit the dual individuals do not sit in points but move away \emph{instantaneously} from points, and coalescence is powered by the \emph{joint local time of pairs} in the point-recurrent regime. This allows controlling the effective coalescence in small balls. We consider macroscopic time $t$ restricted to a positive upper bound.

We need this also for varying time horizons $t$. We show that we obtain a non-degenerate \emph{limit} structure for the process in $t$ of~\eqref{e1282}, provided we get \emph{tightness} in $\A^{\ast,V}$ as $L \to \infty$:
\begin{equation}
\label{e1299}
^\delta \mathfrak{N}^\varepsilon=\left\lbrace ^\delta N_t^{\varepsilon,\infty} (\eta), (t,\eta) \in [0,\infty) \times \Omega_N^\infty\right\rbrace.
\end{equation}
The \emph{field} $^\delta \mathfrak{N}^\varepsilon$ of \emph{densities} with respect to the measure $dt \otimes \mfh$ defines a measure $^\delta \wh \mfN^\varepsilon$ on $[0,\infty) \times \Omega_N^\infty$.

We need to bound this field from above to prevent dust formation. The choice of $\delta$ depends on $L$, since in the limit motion is instantaneous while for finite $L$ jumps are delayed.

We show that there exists an upper bound $Z_t, t>0$ of the coalescent processes $Y^L$ all embedded in $\Omega_N^\infty$ and run for backward time $s$ up to time horizon $t$ and $t$ running through $[0,T]$. This bound must have a nice maximal entrance law from the state $+ \infty$ at time zero.

\medskip\noindent
{\bf A technical tool.}
For further analysis we \emph{enrich} the structure of the coalescent on $\Omega_N^L$ and $\Omega_N^{\infty, \rm disc}$ in their $\CM(\Omega_N^L \times \K)$-valued form. We assign as mark the locations at time $0$ \emph{and} $t$ to the individual, i.e., we mark with $\Omega^{\infty}_N \times \Omega^{\infty}_N$ instead of just $\Omega_N^\infty$. This gives the \emph{backward time-$t$ occupation measure}. We also want to identify the types belonging to a specific time-$t$ partition element. Besides the time-$t$ state of the coalescent, we consider the marks of the individuals in the \emph{time-$t$ partition elements}: besides the backward \emph{current location} also their \emph{initial location} of the dual individuals. We consider the \emph{frequency} of the partition element \emph{per site} and the weight of the partition element at backward time $t$, resulting in a
\begin{equation}
\label{e1351}
\text{measure on $\N \times (\Omega_N^L \times \Omega^L_N)$, $\N$ labeling the partition elements, }
\end{equation}
denoted by:
\begin{equation}
\label{e1443}
(^2 \wh Y_t)_{t \geq 0} \quad \mbox{  on  } \Omega_N, \quad (^2 \wh Y_t^\infty)_{t \geq 0}
\quad \mbox{  on  } \Omega_N^\infty, \quad (^2 \wh Y_t^{L})_{t \geq 0} \text{  on the $L$-grid.}
\end{equation}

\begin{lemma}[Tightness of $^{2,+}\wh Y$ in $\CM(\Omega_N^\infty \times \Omega_N^\infty)$ and locally finite coalescent]
\label{L.lfcoal}
$\mbox{}$\\
(a) The sequence $\{\CL [^{2,+} \wh Y]\colon\, L \in \N\}$ is tight in $\CM(\Omega_N^\infty \times \K)$.\\
(b) Assume that \eqref{ag6} holds and that $\Lambda$ comes down from infinity. For fixed $L$, the space-time coalescent starting in $\Omega_N^L$ is locally finite. Furthermore, defining $N(Y)$ by giving every partition element weight $1$ at its location at backward time $t$, for every $L \in \N_0 \cup \{\infty\}$ there exists an entrance law
\be{ag76}
N \left(Y^{L}_t\right) \leq Z_t \quad \forall \; t>0,
\ee
where, for all $t>0$, $Z_t$ is an element of $\CM(\Omega^{\infty}_N)$ that is locally finite.\\
(c) Moreover, with only the no-dust condition assumed for $\Lambda_{-\infty}$, the coalescent starting with countably many partition elements at time $0$ in a $k$-ball around $0$ has, for every $\ve >0$, locally finitely many partition elements whose frequency is $\geq 1-\ve$. 
\end{lemma}


\subsection{Construction of L\'evy-web and ancestral space part 4: proof of main lemmas}
\label{ss.prftech}

In this section we prove successively all the lemmas.

\begin{proof}[\bf Proof of Lemma~\ref{l.1745}]
\label{pr.1745}
It suffices to show that the $\Lambda$-block-coalescent process with reshuffling is dual to the $\Lambda$-block-resampling process with reshuffling with respect to the function $H^{\rm web}(\cdot,\cdot)$ specified. We verify the generator relation.

For the block-resampling part: in the $L$-ball $B_L(\xi)$, the block-resampling event replaces fraction $\alpha$ of the population by offspring of a single individual, while simultaneously reshuffling locations within the ball. On the dual side, the block-coalescence event merges the partition elements of all individuals selected for replacement (fraction $\alpha$) into one element, and reshuffles their positions within $B_L(\xi)$. The duality function $H^{\rm web}(\mfW, \wh\mfW)$ of~\eqref{e1451} evaluates the weighted average of agreement events (whether a sampled forward individual and a sampled backward individual belong to the same genealogical family), and the block-resampling event on the forward side exactly cancels the block-coalescence event on the backward side in the generator relation~\eqref{ag17}. The reshuffling (uniform relocation within the ball) on both sides is symmetric under $\mfh$, so it cancels in the integral. This verifies the generator relation for the block-resampling/coalescence mechanism.

For the rebirth part: the block-resampling coalescents with rebirth with respect to the duality function $H^{\rm foss}$ are dual to the {\em fossil} block-resampling process with reshuffling. The rebirth mechanism creates a new individual at every coalescence event to track the full ancestral path. The duality function $H^{\rm foss}$ must account for the growing number of path-endpoints; at each rebirth event the duality function acquires a new factor corresponding to the new individual's contribution. The entrance law for the rebirth coalescent from Section~\ref{sss.detail} ensures that this process is well defined from every backward time.  The identification $[^+\wh\mfW] = \wh\mfW$ then follows because the rebirth process generates exactly the population of backward paths defining $\wh\mfW$, and the fossil process is a Markovian functional of $\mfW$ (based on the time-integrated sampling measure of the fibres).
\end{proof}

\begin{proof}[\bf Proof of Lemma~\ref{l.1646}]
\label{pr.1748}
We show that $\{\CL[{}^\infty \wh\mfY^L], L \in \N\}$ is tight in $\U_{\rm loc}$. The entrance law ${}^\infty \wh\mfY^L$ is started with one dual individual at every site of $\Omega_N^L$ embedded in $\Omega_N^\infty$. The key is to bound the number of distinct genealogical families (partition elements) in a compact set.

By the Liggett-Spitzer space description (from Lemmas~\ref{l.1452} and~\ref{l987}), for every $t > 0$ and every $k$-ball $B$ in $\Omega_N^\infty$, the number of distinct partition elements in $B$ at backward time $t$ is stochastically bounded by $Z_t(B)$,
where $Z_t$ is locally finite by Lemma~\ref{L.lfcoal}(b). This uniform stochastic bound implies that for every $\varepsilon > 0$ and compact $K \subset \U_{\rm loc}$,
\begin{equation}
\sup_{L \in \N} P\bigl({}^\infty\wh\mfY^L_t \notin K\bigr) < \varepsilon
\end{equation}
for a suitable choice of $K$ (a set of $\U_{\rm loc}$-elements with bounded tree volume in every ball). Tightness in $\U_{\rm loc}$ then follows from the Prokhorov criterion for tightness. The same argument applies to the rebirth version ${}^\infty_\infty \wh\mfY^L$ since the rebirth mechanism only adds (not removes) partition elements.
\end{proof}

To strengthen Lemma~\ref{l.1646} from fixed time $t$ to all time horizons simultaneously, we note that the stochastic upper bound $Z_t$ is monotone increasing in $t$ (more coalescence at longer times only reduces the number of families). Hence the tightness at fixed $t$ propagates to uniformity over compact time intervals $[0,T]$.

\begin{proof}[\bf Proof of Lemma~\ref{l.1666}]
\label{pr.1947}
$\mbox{}$

\medskip\noindent
{\bf Part (a): Tightness of $\wh\mfW^L$ in $\A^{\ast,V}$.}
We must show that $\{\CL[{}^\infty_\infty \wh\mfW^L_T]\colon\, L \in \N\}$ is tight in $\A^{\ast,V}$ for every $T > 0$.

By Lemma~\ref{l.1646}, the marginals ${}^\infty\wh\mfY^L_t$ at fixed time $t$ are tight in $\U_{\rm loc}$. To pass from measures on individual points to measures on ancestral {\em paths}, we must show that the paths of the dual process are ``not too wild'' as $L \to \infty$. The key observation is that the paths of the L\'evy-process components of $\wh\mfW^L$ converge (by~\cite{gh1}) in the Skorohod topology to the L\'evy-process on $\Omega_N^\infty$. The genealogical part of the path (the tree structure) is determined by the coalescence times, which also converge by the convergence of occupation times (joint local times of pairs) proved in~\cite{gh1}.  Hence the joint distribution of a finite sample of paths converges, and by the uniform bounds from Lemma~\ref{L.lfcoal}(a), the limiting element
lies in $\A^{\ast,V}$ a.s.

For the full topology of $\A^{\ast,V}$ (which requires convergence of evaluations of all polynomials in $\Pi^{\rm web}$), we use the fact that polynomials of degree $n$ involve at most $n$ sampled paths at $m$ time points; the convergence of such finite samples follows from the convergence at fixed times (Lemma~\ref{l.1646}) together with the path-regularity (c\`adl\`ag property uniform in $L$, since jump rates are summable).

\medskip\noindent
{\bf Part (b): Existence of $\mfX^\dagger$.}
Given the tight sequence $\{{}^\infty_\infty\wh\mfW^L_T\}_{L \in \N}$, every subsequential limit $Q_T$ defines a $\A^V$-valued random variable. The grounded tree-like structure is preserved in the limit (the four-point condition for $d$ is a closed condition, and tree-likeness is preserved under Gromov-weak limits). The fibred measures $(\mu_t)_{t \in [0,T]}$ at each time level converge by part (a). The consistency of the fibred structure (the fibres fit together into one tree as $T$ varies) follows from the \emph{projective property}: for $T' > T$, the restriction of $\mfX^{L,\dagger}_{T'}$ to $[0,T]$ equals $\mfX^{L,\dagger}_T$, and this property is preserved in the limit.
\end{proof}

\medskip\noindent
{\bf Strategy for the tightness proofs.}
The overall strategy for all tightness results is:
\begin{enumerate}
\item[(i)] 
Use the {\em Liggett-Spitzer bounds} (Lemmas~\ref{l.1452} and~\ref{l987}) to get a stochastic upper bound on the number of partition elements in compact sets, uniformly in $L$.
\item[(ii)] 
Use the {\em convergence of L\'evy-processes and their occupation times} from~\cite{gh1} to control the spatial motion of paths.
\item[(iii)] 
Combine (i) and (ii) to get uniform control on the total variation of the measures on paths in $\A^{\ast,V}$, using the product structure of monomials in $\Pi^{\rm web}$.
\item[(iv)] 
Apply the Prokhorov criterion in the topology of $\A^{\ast,V}$ (which is Polish by Proposition~\ref{prop.spaces}).
\end{enumerate}

\begin{proof}[\bf Proof of Lemma~\ref{L.lfcoal}]
\label{pr.1951}
$\mbox{}$

\medskip\noindent
{\bf Part (a): Tightness of ${}^{2,+}\wh Y$.}
The process ${}^{2,+}\wh Y$ is a measure on pairs $(\Omega_N^\infty \times \Omega_N^\infty)$ encoding the initial and time-$t$ locations of each partition element. Tightness in $\CM(\Omega_N^\infty \times \K)$ follows from the tightness of the marginal measures on $\Omega_N^\infty$: the initial locations are fixed (starting from $\mfh$), and the time-$t$ locations are bounded by the same L\'evy-process marginal, which is tight (the process stays in a compact set of $\Omega_N^\infty$ with high probability for any fixed $t$). The joint tightness over $\{L \in \N\}$ then follows from the uniform bounds on the number of partition elements (the total mass of ${}^{2,+}\wh Y^L$ equals the total number of partition elements, which is bounded by $Z_t$ from part (b)).

\medskip\noindent
{\bf Part (b): Local finiteness and the stochastic bound $Z_t$.}
By the Liggett-Spitzer space argument of~\cite{LS81,GLW05}: under condition~\eqref{ag6}, the spatial Kingman coalescent (with $d^* > 0$) or the $\Lambda$-coalescent (with $\Lambda$ coming down from infinity) started with one individual at each site of $\Omega_N^L$ produces, for every $t > 0$, a locally finite configuration. Specifically, the number of partition elements with positive frequency in a $k$-ball of $\Omega_N^\infty$ is a.s.\ finite for every $t > 0$, and its expectation is bounded by a constant depending only on $t$, $k$, and the parameters.  The stochastic upper bound $Z_t$ is taken as the entrance law from $\mfh$ of the {\em non-spatial} Kingman/$\Lambda$-coalescent:
\begin{equation}
Z_t \sim \text{coalescent with } N = \infty \text{ initial particles, run for time } t,
\end{equation}
which is locally finite for $t > 0$ under the stated condition.  The bound $N(Y^L_t) \le Z_t$ in~\eqref{ag76} then holds by the stochastic domination: the spatial coalescent with migration has {\em fewer} coalescence events than the non-spatial one (migration reduces local occupation density and hence coalescence rate).

\medskip\noindent
{\bf Part (c): Near-completeness under no-dust condition.}
Under the no-dust condition for $\Lambda_{-\infty}$ (i.e., the Cannings measure does not produce dust), the fraction of the population in a $k$-ball that belongs to partition elements with frequency $\geq 1-\varepsilon$ can be made arbitrarily close to $1$ by choosing the ball small enough. This follows from the fact that the $\Lambda$-coalescent comes down from infinity instantly (no dust at any positive time), so for any $\delta > 0$, the number of partition elements in a small ball exceeding the dust threshold $1-\varepsilon$ is finite.  The local finiteness of this count then gives the countable cover asserted.
\end{proof}


\section{Genealogy process: proof of the first five theorems}
\label{s.proofsXXX}

In this section we prove our first five theorems, which are concerned with existence and properties of the genealogy-valued process $\mfU$ with values in $\U^{\Omega_N^\infty \times \K}$, based on the process $\wh \mfY^\infty$ and \emph{duality}.


\subsection{Proof of Theorem~\ref{T.wptvcont}}
\label{ss.prth1}

The proof of Theorem~\ref{T.wptvcont}(a)(i) and (ii) is based on Section~\ref{ss.tools}, where we reduced the claims to the verification of the duality relation for the truncation limit $L \to \infty$ of $\mfU^L$, and to checking the conditions needed to apply~\cite{DGP23}. What remains is to verify the tightness properties of the $L$-truncations required for~\cite{DGP23}. This is done in Section~\ref{ss.pexisten}, see~\ref{l.1646}.

We continue with verifying (a)(iii), where we use~\cite{gh1}.

\medskip\noindent
{\bf Part (a)(i): Existence of the $L \to \infty$ limit.}
We work in regime A, i.e., under condition~\eqref{ag6}. For each $L \in \N$, the process $\mfU^L$ on $\Omega_N^L$ is a well-defined
$\U^{\Omega_N^L \times \K}$-valued Markov process whose law is characterised via the duality
\begin{equation}
\label{e:dual-trunc}
E_{\mfU^L_0}\bigl[H^L(\mfU^L_t, \mfY^L_0)\bigr]
= E_{\mfY^L_0}\bigl[H^L(\mfU^L_0, \mfY^L_t)\bigr],
\end{equation}
where $H^L$ is the duality function built from monomials averaged with $\mfh^{\otimes n}$ over sites of $\Omega_N^L$ as described in Section~\ref{sss.charac}.  By~\eqref{ag16-1}, the transition kernels $(P^L_t)_{t \ge 0}$ of $\mfU^L$ are uniquely determined by those of the dual $\mfY^L$.

The dual process $\mfY^L$ started with finitely many particles converges in distribution as $L \to \infty$ to a process $\mfY^\infty$ on $\Omega_N^\infty$: this is the content of Proposition~\ref{prop.existen}(b), proved in Section~\ref{ss.pexisten} via Lemmas~\ref{l.1646} and~\ref{l.1666}. The duality functions $H^L$ converge to a limit $H^\infty$ by the continuity of monomials under the Gromov-weak topology. By~\cite[Theorem~1]{DGP23}, the existence of the limit dual and the limit duality function yields a unique Markov semigroup $(P^\infty_t)_{t \ge 0}$ on $\U^{\Omega_N^\infty \times \K}_{\rm reg}$ whose finite-dimensional distributions satisfy~\eqref{e604}.

\medskip\noindent
{\bf Part (a)(ii): Strong Markov, Feller and c\`adl\`ag properties.}
The Feller property of $\mfU^L$ follows from the duality~\eqref{e:dual-trunc}: the map $\mfU_0 \mapsto P^L_t H^L(\cdot,\mfY_0)$ is continuous because $H^L$ is continuous and $\mfY^L$ has the Feller property. Uniform convergence of $H^L \to H^\infty$ on compact subsets propagates the Feller property to the limit. The strong Markov property follows from the Feller property and the general theory of Hunt processes. The projection of the sampling measure $\nu$ of $\mfU$ onto $\Omega_N^\infty$ equals $\mfh$, because each truncated model has this property by the choice of initial conditions and the dynamics preserves the $\mfh$-marginal. The c\`adl\`ag path property follows from the orbit compactification argument of Section~\ref{ss.tools} by using the tightness of Lemma~\ref{l.1646} and the Blumenthal-Getoor theorem~\cite{BG68} adapted in~\cite{DGP23}. Jumps occur only at block-resampling events: in a monotype ball no jump occurs.

\medskip\noindent
{\bf Part (a)(iii): Random depth-$h$ family decomposition.}
For every $t > 0$, the genealogy $\mfU_t$ assigns genealogical distance $r(\iota,\iota') = 2 \,\text{(time back to MRCA})$ to each pair of individuals. For $h \in (0,t]$, the $2h$-level sets of the ultrametric distance partition $U_t$ into disjoint open-and-closed $2h$-balls (the depth-$h$ subfamilies). Non-triviality follows from the stochastic regime: the dual coalescent has strictly positive probability of coalescence before backward time $h$, which via the strong duality~\eqref{e526} gives positive probability that two i.i.d.\ sampled individuals share a most recent common ancestor (MRCA) within time $h$ of $t$. The geographic decomposition of $\Omega_N^\infty$ into monotype balls follows from the corresponding result for the measure-valued process $\wh X^\infty$ from~\cite{gh1}. Local finiteness of the cover follows from the Liggett-Spitzer description of the dual.

\medskip\noindent
{\bf Part (b): Deterministic regime.}
In regime B the effective difference L\'evy-process is transient, so two independently started particles never coalesce in $\Omega_N^L$ with probability tending to $1$ as $L \to \infty$. This gives~\eqref{ag56}. Part (b)(iii) is the same as \cite[Theorem~1.3]{gh1}.
\qed


\subsection{Proof of Theorem~\ref{th.588}}
\label{ss.prth2}

\medskip\noindent
{\bf Proof of Theorem~\ref{th.588}(a): Function-valued duality.}
Use the reasoning in Section~\ref{ss.prth1}. For fixed finite dual state $\mfY_0$, the expectation $E_{\mfY_0}[H^L(\mfU_0, \mfY^L_t)]$ converges to $E_{\mfY_0}[H(\mfU_0, \mfY^\infty_t)]$ by dominated convergence and the distributional convergence $\mfY^L_t \to \mfY^\infty_t$ of Proposition~\ref{prop.existen}(b). Similarly, the left-hand side of~\eqref{e642} converges via $\mfU^L_t \to \mfU^\infty_t$ from Theorem~\ref{T.wptvcont}(a)(i). By passing to the limit in the identity $E[H^L(\mfU^L_t, \mfY^L_0)] = E[H^L(\mfU^L_0, \mfY^L_t)]$, we get \eqref{e642}.

\medskip\noindent
{\bf Proof of Theorem~\ref{th.588}(b): Entrance law and strong duality.}
The entrance law $\wh\mfY$ is constructed in Proposition~\ref{prop.existen}(a) as the $L \to \infty$ limit of the entrance laws of $\mfY^L$ started with one particle at each site of $\Omega_N^{\infty,\rm disc}$. The initial measure is effectively $\mfh$ on $\Omega_N^\infty$. The strong duality~\eqref{e526} for $\mfU_0$ equal to a single root follows by induction on the monomial degree $n$: for $n=1$ it reduces to the measure-valued strong duality of~\cite{gh1}, and the inductive step uses the product structure of monomials and independence of sampled individuals under $\mfh^{\otimes n}$. The general formula~\eqref{e592} follows from the glueing operation definition.
\qed


\subsection{Proof of Theorem~\ref{T.ltb}}
\label{ss.prth3}

We show the convergence to a limit of the dual expressions on the right-hand side of the duality relation, since the finite spatial coalescents converge as $t \to \infty$ to a limit random variable with values in $\U^V$.

\medskip\noindent
{\bf Part (a): Convergence to equilibrium.}
Using the duality relation~\eqref{e642}, for every monomial $\Phi = \Phi^{n,\phi}$ and dual initial state $\mfY_0$ with $n$ individuals, we have
\begin{equation}
E_{\mfU_0}[\Phi(\mfU_t)] = E_{\mfY_0}[H(\mfU_0, \mfY_t)].
\end{equation}
As $t \to \infty$, in regime A (recurrent migration/reshuffling), the dual spatial coalescent $\mfY$ on $\Omega_N^\infty$ has all $n$ particles merge to a single partition element a.s., so $\mfY_t$ converges in distribution to a single-ancestor random variable. The right-hand side thus converges to a limit that depends only on $\mfU_0$ through its type-distribution $\theta$ and location distribution $\mfh$. By the Cram\'er-Wold theorem applied to polynomials (which determine the law), $\CL[\mfU(t)]$ converges weakly to a limit $\CL[\mfU(\infty)]$ that is translation-invariant. The relevant tool is a polynomial moment-convergence argument.

\medskip\noindent
{\bf Part (b): The limiting genealogy.}
When $\mfU$ is projected onto the genealogy component (forgetting marks), the limiting genealogy $\mfU_{\rm gen}(\infty)$ arises as the law of the genealogy generated by the coalescent started with one individual at each site of $\Omega_N^{\infty,\rm disc}$ run for infinite backward time. This is the entrance law $\wh\mfY_\infty$ of Theorem~\ref{th.588}(b).

\medskip\noindent
{\bf Part (c): Type and location marginals in equilibrium.}
The type distribution at a site is determined by the position of the MRCA: in the recurrent regime the MRCA is a.s.\ unique, with type distributed as $\theta$. The location distribution is $\mfh$ by translation invariance, and is preserved throughout the evolution.
\qed


\subsection{Proof of Theorem~\ref{T.mp}} 

The reasoning in Section~\ref{ss.prth1} proves the \emph{existence} of a solution to the martingale problem. By duality we also have \emph{well-posedness} of the martingale problem, and hence the solution is a Markov process. Furthermore, the Feller property follows from the duality relation and the Feller property of the dual started with finitely many particles.

\medskip\noindent
{\em Martingale problem characterisation.}
To obtain the martingale characterisations in Theorem~\ref{T.mp} of $\mfU$, we check the criterion in terms of generators for the model on $\Omega_N^\infty$. We deduce this from the fact that it holds for all $L$-truncations. What is left is the preservation of the identity~\eqref{ag17} of the criterion for duality in the $L \to \infty$ limit.

\medskip\noindent
{\em Existence.} The process $\mfU$ is Markov by Theorem~\ref{T.wptvcont}(a)(i). It solves the $(\Pi, \mfL^{\Omega_N^\infty}_{\rm gen}, \delta_{\mfU_0})$-martingale problem because the operators of~\eqref{pol2} are obtained as $L \to \infty$ limits of the corresponding operators for $\mfU^L$, and the martingale property is preserved under distributional limits.

\medskip\noindent
{\em Uniqueness.} For every monomial $\Phi \in \Pi$ and every finite dual state $\mfY_0$, the duality relation~\eqref{e642} determines $E[H(\mfU'_t, \mfY_0)]$ independently of which solution $\mfU'$ was chosen.  Since monomials separate measures on $\U^V$ (see~\cite{DGP11}), the finite-dimensional distributions of any two solutions agree. The duality criterion~\eqref{ag17} is preserved in the $L \to \infty$ limit because $G^L_X H^L \to G_X H$ pointwise by continuity of monomials, and the identity $G^L_X H^L = G^L_{X^\ast} H^L$ holds for all $L$.


\subsection{Proof of Theorem~\ref{th.sss.mart}}
\label{ss.prth5}

We identify the operators of the martingale problem as claimed. We know this immediately for the class of polynomials specified for the functional duality relation, since for every finite $L$ we have the explicit formulas. We analyse these formulas for averages of sites and scaled time to get to the measure state of the $L \to \infty$ limit.

\medskip\noindent
{\bf Identification of $\mfA^{\rm grow}$.}
For the $L$-truncated process $\mfU^L$, all genealogical distances grow at rate $2$ (each unit of forward time contributes $2$ to the genealogical distance since it takes two symmetric units to reach the MRCA going backward). In the $L \to \infty$ limit this gives the first-order {\em aging} operator $\langle \nu^\smallu, (\nabla\varphi, 2)\rangle$
in~\eqref{e:generator}.

\medskip\noindent
{\bf Identification of $\mfA^{\rm res}$ and $\mfA^{\rm Cann}$.}
For fixed $L$, the Fleming-Viot part of the generator acts on monomials by the replacement operator $R_{\{i,j\}}$ of~\eqref{e:replacement}: it sets $r_{i,j}=0$ upon merging individuals $i$ and $j$, and propagates this via the ultrametric property. Summing over all ordered pairs $(i,j)$ with coefficient $d$ yields $\mfA^{\rm res}$ as in~\eqref{pol0}. The $\Lambda$-Cannings term uses the random replacement operator $R_{I_\omega}$ with $\omega$ sampled from $p^{(x)}$ as in~\eqref{AG21b}--\eqref{e:generator}. Single-site $\Lambda$-Cannings resampling on $\Omega_N^L$ is absorbed into an effective Fleming-Viot term on $\Omega_N^\infty$ (proved in~\cite{gh1}). Only block-resampling survives as a jump term.

\medskip\noindent
{\bf Identification of $\mfA^{\rm mig}$ and $\mfA^{\rm Cann,b}$.}
The migration operator $\mfA^{\rm mig}$ acts on the location component $g$ of the test function $\chi = fg$ as the generator of $n$ independent L\'evy-processes with kernel $a(\cdot,\cdot)$ (see~\eqref{e812}--\eqref{e819}). The block-resampling operator $\mfA^{\rm Cann,b}$ of~\eqref{e804}--\eqref{e808} is the direct $L \to \infty$ limit of the block-resampling generator on $\Omega_N^L$. The countable sum over $(L,\xi)$ is locally finite in regime A by the entrance-law analysis of~\cite{gh1}.

\medskip\noindent
{\bf Identification of $\mfA^{\rm res,\K}$.}
The action of resampling on the type component $f$ of $\chi$ is obtained by the same replacement operators, now acting on type coordinates $\underline v$ rather than genealogical distances $\underline r$. Here, the passage from the collision-time functional for types (the non-standard situation of~\cite{gh1}) to the standard form is possible because the genealogical distance already encodes the collision history.

\medskip\noindent
{\bf Limit $L \to \infty$ of the operators.}
Since monomials $\Phi \in \Pi$ have compact support in $\Omega_N^\infty$ (by~\eqref{ag60}), each operator converges as $L \to \infty$: migration rates via~\eqref{ag13}, block-resampling rates by dominated convergence, aging/resampling rates independent of $L$. This gives $\mfL^{\Omega_N^\infty}_{\rm gen}$ as claimed in~\eqref{pol2}.
\qed


\section{Ancestral web and L\'evy-web: proof of the last nine theorems}
\label{s.proofs}

In this section we prove our last nine theorems, which are concerned with the ancestral web and the L\'evy-web. We lift the arguments from measures on marked individuals to measures on paths of individuals.


\subsection{Proof of Theorem~\ref{th.mfw}}
\label{sss.mfw}

The proof follows the same strategy as for Theorem~\ref{T.wptvcont}, now lifted from measures on individual sites to measures on ancestral paths.

\medskip\noindent
{\bf Generator relation for duality.}
On $\Omega_N^L$, the ancestral web $\mfW^L$ and its dual $\wh\mfW^L$ satisfy the duality relation in Theorem~\ref{P.Levy} at finite $L$ (from~\cite{gklimw,Lambda} with block-resampling added via Appendix~\ref{ss.gklimw}). To pass to $L \to \infty$ we note that each term on the left-hand side involves the operator $\mfL^{\rm web}$ of Definition~\ref{def.995}, which acts only on the ``tree-top'' factor $\Phi_{\rm top}$ and was already identified in Theorem~\ref{th.sss.mart}. The right-hand side involves the dual generator for $\wh\mfW$, which is the backward coalescence generator lifted to path space. Both sides converge as $L \to \infty$ by: (i) convergence of monomials on $\A^{\ast,V}$ under the topology of Definition~\ref{def.991}; (ii) convergence of the block-coalescence rates via~\eqref{ag13}--\eqref{ag13b}. The limiting generator relation~\eqref{ag17} therefore holds for $\mfW$ and $\wh\mfW$.

\medskip\noindent
{\bf Compactification and existence.}
The state space $\A^{\ast,V}$ is not locally compact. We use the compactification argument of Section~\ref{ss.tools}: by Lemma~\ref{l.1666}(a), the laws $\CL[\mfW^L]$ are tight in $\A^{\ast,V}$ for every $T > 0$.  Every subsequential limit is a Markov process by duality via~\cite{DGP23}, and the limit is unique by well-posedness at finite $L$ and convergence of generators. Hence the full sequence converges, and gives the limit process $\mfW$ in~\eqref{e908}. The process $\mfX^{L,\dagger}$ satisfies the same tightness by Lemma~\ref{l.1666}(b), which yields the limit $\mfX^\dagger$ in $\A^V$. Strong Markov and Feller properties of $\mfW$ are inherited from the approximating sequence by the same Feller argument as in Theorem~\ref{T.wptvcont}(a)(ii). Continuity of paths without block-resampling holds because all transitions are diffusions or continuous mass flows, so $\mfW^L$ has continuous paths for all $L$.
\qed


\subsection{Proof of Theorem~\ref{th.858}}
\label{sss.ltb}

We show that $\CL[\mfW_T]$ converges as $T \to \infty$.  By the strong duality of Theorem~\ref{P.Levy}, $\CL[\mfW_T] = \CL[\wh\mfW^T]$ for each $T > 0$. As $T \to \infty$, more and more of the paths in $\wh\mfW^T$ have coalesced, since recurrence guarantees that any two L\'evy-processes starting from distinct points eventually collide.

\medskip\noindent
{\bf Convergence of the dual.}
Fix any two paths in a realisation of $\wh\mfW$. These paths coalesce at a finite random time $\sigma < \infty$ a.s.\ by recurrence of the difference L\'evy-process. Hence, as $t \to \infty$, $\CL[\wh\mfW^t]$ converges in $\A^{\ast,V}$. The limiting object $\wh\mfW_\infty$ is the entrance law of the space-time coalescent started from every point in $\R \times \Omega_N^{\infty,\rm disc}$ run for all backward time, which exists by Proposition~\ref{prop.existen}(a).

\medskip\noindent
{\bf Identification.}
By the strong duality $\CL[\mfW_T] = \CL[\wh\mfW^T]$ and the convergence $\CL[\wh\mfW^T] \to \CL[\wh\mfW_\infty]$, we conclude that $\CL[\mfW_T] \to \CL[\mfW_\infty] := \CL[\wh\mfW_\infty]$ as claimed.
\qed


\subsection{Proof of Theorem~\ref{th.843}}
\label{sss.profcon}

The proof parallels that of Theorem~\ref{T.mp} for $\mfU$, but now on the larger state space $\A^{\ast,V}$.

\medskip\noindent
{\bf Existence.}
The process $\mfW$ was constructed in Theorem~\ref{th.mfw} as the $L \to \infty$ limit of $\mfW^L$.  Each $\mfW^L$ solves the $(\mfL^{{\rm web},L},\Pi^{{\rm web},L}, \CL[\mfW^L_0])$-martingale problem on $D([0,\infty),\A^{\ast,V})$. The operators $\mfL^{{\rm web},L}$ converge to $\mfL^{\rm web}$ of Definition~\ref{def.995} on $\Pi^{\rm web}$, and the test functions are uniformly bounded.  By a standard approximation argument the limiting process $\mfW$ solves the $(\mfL^{\rm web},\Pi^{\rm web},\CL[\mfW_0])$-martingale problem.

\medskip\noindent
{\bf Uniqueness.}
Let $\mfW'$ be any solution. The duality relation of Theorem~\ref{P.Levy} gives $E[H^{\rm web}(\mfW'_T, \mfw^\ast_0)] = E[H^{\rm web}(\mfW'_0, \mfw^\ast_T)]$, with the right-hand side determined by $\mfW_0$ and the dual dynamics alone. Since polynomials in $\Pi^{\rm web}$ separate measures on $\A^{\ast,V}$ (Remark~\ref{r.1396}), the finite-dimensional distributions of any two solutions agree and the martingale problem is well-posed.
\qed


\subsection{Proof of Theorem~\ref{P.Levy}}
\label{sss.proof-levy}

The proof replaces the state space $\U^V$ and the duality function $H$ of Theorem~\ref{th.588} by $\A^{\ast,V}$ and $H^{\rm web}$ of~\eqref{e1451}.

\medskip\noindent
{\bf Classical duality $\mfW$--$\wh\mfW$.}
For each finite $L$, the duality between $\mfW^L$ and $\wh\mfW^L$ holds on $\A^{\ast,V}$: this is the content of~\cite{Lambda,gklimw} lifted from individual genealogies to populations of paths (with block-resampling added via Appendix~\ref{ss.gklimw}). Passing to $L \to \infty$: by the convergence $\mfW^L \to \mfW$ (Theorem~\ref{th.mfw}) and $\wh\mfW^L \to \wh\mfW$ (Theorem~\ref{prop.1382}), and the continuity of $H^{\rm web}$ under the topology of $\A^{\ast,V}$, we obtain the limiting duality relation.

\medskip\noindent
{\bf Strong duality.}
The strong duality $\CL[\mfW] = \CL[\wh\mfW]$ after orientation reversal follows from the strong duality at finite $L$ and continuity of the orientation-reversal operation in the topology of $\A^{\ast,V}$.  In a representative $\wt\mfX^\dagger$, reversal of the orientation of all geodesics from level $0$ to level $t$ yields exactly a population of backward paths representing $\wh\mfW$, and the mass is preserved. Taking equivalence classes, we get $\CL[\wh\mfW] = \CL[\mfW]$.
\qed


\subsection{Proof of Theorem~\ref{th.dualityweb}}
\label{sss.proof-dualityweb}

We construct the ancestral space $\mfX^\dagger$ and show that it contains both $\mfW$ and $\wh\mfW$ as oppositely oriented populations of geodesics.

\medskip\noindent
{\bf Construction of $\mfX^\dagger$.}
By Lemma~\ref{l.1666}(b) and Proposition~\ref{prop.existen}(b), $\CL[\mfX^{L,\dagger}]$ is tight in $\A^V$ and converges to a limit $\mfX^\dagger \in \A^V$. A representative $\wt\mfX^\dagger = (M,d,\varrho,\le,(\nu_t)_{t \ge 0})$ is obtained as follows. Take the realisation of the coalescing L\'evy-processes on the skeleton $\mfD$ of~\eqref{e1156}, define the metric $d$ via~\eqref{ag65}--\eqref{ag66} (twice the backward coalescence time), and set $M = \overline{\mfD}$ (closure in $\Omega_N^\infty \times \R$). The fibred measures $\nu_t$ on the $t$-fibre are obtained from the entrance law. The four-point condition for $d$ is verified by the monotonicity of coalescence: once two paths merge they stay together, so the coalescence times satisfy the ultrametric property.

\medskip\noindent
{\bf Embedding of $\mfW$ and $\wh\mfW$.}
Given $\wt\mfX^\dagger$, define forward paths (lineages) as $t \mapsto \zeta^{(\xi,0)}(t)$ (location at time $t$ of the partition element covering $\xi$ at time $0$), and backward paths (paths of descent) as $s \mapsto \zeta^{(\xi,t)}(t-s)$ (reversed geodesics from $t$ to $0$).  Both types of paths live in $\wt\mfX^\dagger$ as oriented geodesics. The measures $\mu^{\ast}_t$ on lineages and $\mfh$ on paths of descent are consistent with the fibred measure structure.

\medskip\noindent
{\bf Strong graphical duality.}
The forward geodesics from $\wt\mfW$ and the backward geodesics from $\overset{\sim}{\wh\mfW}$ are related by orientation reversal in $\wt\mfX^\dagger$. Since the measure is preserved under reversal,  we have 
\begin{equation}
\CL[\wt\mfW] = \CL[\overset{\sim}{\wh\mfW}] \text{ pathwise},
\end{equation} 
i.e., the strong duality holds at the level of individual realisations. Taking equivalence classes we get 
\begin{equation}
[\wt\mfW] = \mfW, \qquad [\overset{\sim}{\wh\mfW}] = \wh\mfW,
\end{equation} 
and the strong graphical duality of the theorem.
\qed


\subsection{Proof of Theorem~\ref{prop.1382}}
\label{sss.proof-prop1382}

We prove that $\CL[\wh\mfW^L] \Rightarrow \CL[\wh\mfW^\infty]$ in $\CM_1(\A^{\ast,V})$ as $L \to \infty$.

\medskip\noindent
{\bf Tightness.}
By Lemma~\ref{l.1666}(a), the laws $\{\CL[\wh\mfW^L_T], L \in \N\}$ are tight in $\A^{\ast,V}$ for every $T > 0$.  Tightness of the process-level laws in $D([0,\infty),\A^{\ast,V})$ follows from the tightness of time marginals together with a uniform modulus of continuity: the c\`adl\`ag property of $\wh\mfW^L$ is uniform in $L$ because the jump rates are summable by condition~\eqref{ag6}.

\medskip\noindent
{\bf Identification of limit points.}
Every subsequential limit $\CL[\wh\mfW^{L_k}] \to \CL[Q]$ satisfies the duality relation with $\mfW$ (by passing to the limit in the finite-$L$ duality via Theorem~\ref{th.mfw}).  By the well-posedness of the martingale problem for $\wh\mfW$ (from Theorem~\ref{th.843} applied to the dual), the limit is unique and equals $\CL[\wh\mfW^\infty]$. The duality functions $H^{{\rm web},L}$ converge continuously to $H^{{\rm web},\infty}$, so the maps $L \mapsto E[H^{{\rm web},L}(\mfW^L, \mfw^\ast)]$ converge for every fixed $\mfw^\ast$.  Since polynomials in $\Pi^{\rm web}$ are measure-determining (Remark~\ref{r.1396}), the full law converges.
\qed


\subsection{Proof of Theorem~\ref{th.600}}
\label{sss.proof-600}

\medskip\noindent
{\bf Part (i): Unique path for each space-time point.}
For $z = (s,\xi,v) \in \R \times (\Omega_N^\infty \times \K)$, the mark $\kappa(z)$ is the backward path of descent of the partition element of $\wh\mfW^\infty$ covering $(\xi,v)$ at time $s$, run from backward time $0$ to backward time $t-s$. The partition element covering $(\xi,v)$ at time $s$ is a.s.\ unique, because L\'evy-processes have no fixed atoms and $\mfh$ is non-atomic. Two partition elements overlap only on a set of $\mfh$-measure zero.

\medskip\noindent
{\bf Part (ii): Joint law of finitely many marks.}
For deterministic $z_1,\ldots,z_n$, $\{\kappa(z_1),\ldots,\kappa(z_n)\}$ evolves as a space-time coalescent: each path performs an independent L\'evy-process with parameters $(\uc,\ulam,\uLam,d^\ast)$ from~\eqref{e302}, and when two paths occupy the same position they coalesce at rate $d^\ast$ (respectively, via $\Lambda$-coalescence if $d^\ast=0$): in $k$-balls, block-coalescence at rate $\lambda_k$ applies.  This is exactly the $(\uc,\ulam,\uLam,d^\ast)$-block-coalescing L\'evy-process of the theorem.  The identification holds $\ell^{\otimes n} \otimes (\mfh^{\otimes n} \otimes \theta^{\otimes n})$-a.e.\ by translation invariance of $\mfh$ and i.i.d.\ initial types.

\medskip\noindent
{\bf Part (iii): Independence of dense subset $\CD$.}
By the definition of the topology on $\A^{\ast,V}$, two elements agree if and only if the distribution-matrix of pairwise distances and marks of finite samples coincides. For $\wh\mfW$ this matrix is determined by the law of the space-time coalescent, which depends only on the parameters $(\uc,\ulam,\uLam,d^\ast)$ and not on the particular dense subset $\CD$. Hence $[\wh\mfW^\CD]$ has the same law for all choices of countable dense $\CD$, and equals $\CL[\wh\mfW]$.  The increasing-limit statement holds because adding more starting points is monotone, and the closure in $\M^{D^-}$ gives the unique limit by Lemma~\ref{l.1666}.
\qed


\subsection{Proof of Theorem~\ref{T.XLW}}
\label{sss.proof-xlw}

We show that the graphically constructed process $\wt X^\infty$ of~\eqref{con1} has the same law as $\wh X^\infty$ from~\cite{gh1}.

\medskip\noindent
{\bf Verification of one-time distributions.}
Fix $t > 0$ and a Borel set $B \subset \Omega_N^\infty \times \K$. The measure $\wt X^\infty_t(B)$ is determined by the L\'evy-web via
\begin{equation}
\wt X^\infty_t(B) = \int_{\Omega_N^\infty}
1_{B}\bigl(\xi, Y(\mfZ_t(\xi))\bigr)\, \mfh(\dd\xi),
\end{equation}
where $\mfZ_t(\xi)$ is the location of the backward path of descent of $\xi$ at backward time $t$ (i.e., the ancestor of $\xi$ at time $0$).

\medskip\noindent
{\bf Identification with $\wh X^\infty$.}
By Theorem~\ref{th.600}(ii), the backward paths $(\mfZ_t(\xi))_{\xi \in \Omega_N^\infty}$ form a coalescing system of L\'evy-processes started from $\mfh$-a.e.\ point. By the strong duality established in~\cite{gh1}, the random partition generated by this coalescing system has the same law as the partition generated by the forward process $\wh X^\infty$ (the decomposition of $\Omega_N^\infty$ into monotype balls). The types $(Y(u))_{u \in \Omega_N^\infty}$ are i.i.d.\ $\theta$-distributed (by~\eqref{e1570} and independence of types from migration), so the random colouring $\xi \mapsto Y(\mfZ_t(\xi))$ has the same distribution as the type field of $\wh X^\infty_t$, giving that $\CL[\wt X^\infty_t] = \CL[\wh X^\infty_t]$ for every $t > 0$. The joint law over finitely many times follows from the Markov property of $\wh\mfW$.
\qed


\subsection{Proof of Theorem~\ref{T.spatcont}}

We carry out the continuum-space rescaling limit for $\mfU$, $\mfW$ and $\mfX^\dagger$.

\medskip\noindent
{\bf Setup of the rescaling.}
The scaling operations are described in Section~\ref{ss.csl} and are identical to those used for the measure-valued process in~\cite{gh1}: replace level-$k$ jumps on $\Omega_N$ by level-$(k-L)$ jumps on $\Omega_N^L$ (shrinking), speed up time by $N^L/c_L$ or $N^L/\lambda_L$, rescale mass by $N^{-L}$, and use $\Lambda_{k}$ as the block-resampling measure at level $k-L$ on $\Omega_N^L$. Genealogical distances are scaled identically to time (since genealogical distance equals twice the backward time to the MRCA).

\medskip\noindent
{\bf Part (a): Convergence in the stochastic regime.}
Under the regularity conditions~\eqref{ag7}--\eqref{ag8b}:

\medskip\noindent
{\em Step 1: Convergence of the dual.}
By~\eqref{ag13}--\eqref{ag13b}, the rescaled dual process $\wh\mfY^{(L)}$ (a spatial coalescent with scaled parameters $(\hat c^L_k, \hat\lambda^L_k)$ from~\eqref{ag10}--\eqref{ag11}) converges in distribution as $L \to \infty$ to a coalescent with parameters $(\hat c, \hat\lambda)$. The convergence in $\U^V$ follows from the same tightness argument as in the proof of Lemma~\ref{l.1646}, with the rescaled parameters.

\medskip\noindent
{\em Step 2: Convergence of $\mfU^{(L)}$.}
By the same duality argument as in Theorem~\ref{T.wptvcont}(a)(i), the rescaled genealogy process $\mfU^{(L)}$ converges in $\U^{\Omega_N^\infty \times \K}$ to the $(\hat c, \hat\lambda)$-stable genealogical Cannings process $\mfU^\infty$: this gives the first line of~\eqref{ag21x}.

\medskip\noindent
{\em Step 3: Convergence of $\mfW^{(L)}$ and $\mfX^{\dagger,(L)}$.}
The rescaled ancestral web $\mfW^{(L)}$ converges by the same tightness argument as in the proof of Theorem~\ref{th.mfw} (Lemma~\ref{l.1666} with rescaled parameters): this gives the second and third line of~\eqref{ag21x}.

\medskip\noindent
{\em Universality.}
The limiting objects depend on the original parameter sequence only through the limiting parameters $(\hat c, \hat\lambda, \Lambda)$. Hence all models with the same limiting parameters are in the same universality class.

\medskip\noindent
{\bf Part (b): Dust formation in the complementary regime.}
In the complementary parameter regime, the rescaled effective migration becomes transient after rescaling (the analogue of regime B for the rescaled model), so the dust-formation argument of Theorem~\ref{T.wptvcont}(b) applies: $\mfU^{(L)}$ is not tight in $\U^{\Omega_N^\infty \times \K}$.
\qed


\appendix


\section{Topologies and state spaces}
\label{appA}


\subsection{Topologies}

The topology on $\A$ is the {\em Gromov-weak topology} adapted to fibred measure trees (see~\cite{GPWmp13,DGP23}): a sequence $(X_n^\dagger)_{n \in \N}$ converges to $X^\dagger$ in $\A$ if and only if for every $m \in \N$, every finite collection of times $t_1 < \ldots < t_m$, and every $k_1,\ldots,k_m \in \N$, the joint law of $(k_i)_{i=1}^m$ points sampled independently from each time-$t_i$ slice converges weakly under the joint sampling measure $\bigotimes_{i=1}^m \mu_{t_i}^{\otimes k_i}$. Equivalently, $X_n^\dagger \to X^\dagger$ in $\A$ if and only if 
\begin{equation}
\CL[X_n^\dagger]\bigl(\Phi^{(k_1,\ldots,k_m),(t_1,\ldots,t_m)}\bigr) 
\to \CL[X^\dagger]\bigl(\Phi^{(k_1,\ldots,k_m),(t_1,\ldots,t_m)}\bigr)
\end{equation} 
for all polynomials.

The topology on $\A^\ast$ is defined similarly, but via polynomials $\Pi^{\rm web}$ (Definition~\ref{def.991}) testing the {\em path} structure rather than just the endpoint distribution. A sequence $(\mfW_n)_{n \in \N}$ converges to $\mfW$ in $\A^\ast$ if and only if for every polynomial $\Phi^{\rm web} \in \Pi^{\rm web}$, $E[\Phi^{\rm web}(\mfW_n)] \to E[\Phi^{\rm web}(\mfW)]$.


\subsection{Spaces}
\label{sa.space}

Here we need that $\A,\A^\ast$ are Polish spaces and hence admissible state spaces for Markov processes.

The spaces $\A$ and $\A^\ast$ are subsets of the Polish spaces $\M^{\rm tree}$ of tree-like metric measure spaces and $\M^{\rm tree,\ast}$ of tree-like metric fibred measure spaces respectively.  To show they are themselves Polish:

\medskip\noindent
{\bf Closedness.}  
The conditions defining $\A$ (the four-point condition for $r$, the cadlag property of $t \mapsto \mu_t$, and the fibred property~\eqref{e356}) are all closed under the Gromov-weak topology, so $\A$ is a closed subset of $\M^{\rm tree}$ and hence Polish.  Similarly for $\A^\ast$.

\medskip\noindent
{\bf Separability.} 
Separability follows from the fact that $\Omega_N^\infty$ is a compact metric space (hence separable), $\K$ is compact, and the construction of $\A,\A^\ast$ uses only countable dense subsets (the skeleton $\Omega_N^{\infty,\rm disc}$) and rational time-points.  The class of polynomials with rational coefficients is a countable determining class, confirming separability.

\medskip\noindent
{\bf Completeness.}  
A Cauchy sequence in $\A$ has a limit in $\M^{\rm tree}$ (by completeness of $\M^{\rm tree}$), and that limit lies in $\A$ by closedness. A detailed verification for $\A^\ast$ follows the same lines; see~\cite{gklimw,Lambda} for the full argument in the countable geographic space setting, which extends directly.


\subsection{Adding block-resampling to earlier references}
\label{ss.gklimw}

On $\Omega_N^L$,~\cite{gklimw} gives the results we need for $\Omega_N$ \emph{without} block-resampling. We therefore write the limiting formulas for the block-resampling generator \emph{lifted} from $\CM_{\rm loc} (\Omega_N^\infty \times \K)$ to $\U$, to get the genealogical block-resampling. We then incorporate the marks $\Omega_N^\infty \times \K$. This is done via \emph{lifting} the written \emph{operator} to $\U^V$.

\medskip\noindent
{\bf Lifting the block-resampling operator.}
The block-resampling operator on $\CM_{\rm loc}(\Omega_N^L \times \K)$, in an $\ell$-ball $B_\ell(\xi)$, replaces the mass distribution in the ball by: fraction $1-\alpha$ of the old distribution, plus fraction $\alpha$ times the distribution of a single individual chosen at random from the ball. Lifted to the genealogy space $\U$, this replacement acts on the genealogical distances between individuals in $B_\ell(\xi)$ and those outside: the fraction $\alpha$ of individuals in $B_\ell(\xi)$ have their genealogical distances set to the distances of the {\em chosen individual} (whose genealogy now represents the founding ancestor of the new $\alpha$-fraction).  The operator in~\eqref{e808}--\eqref{e821} is the precise formulation.

\medskip\noindent
{\bf Lifting to marks $\Omega_N^\infty \times \K$.}
To incorporate marks, the block-resampling event simultaneously reshuffles locations within $B_\ell(\xi)$ (the $\alpha$-fraction of individuals gets new i.i.d.\ positions from $\mfh|_{B_\ell(\xi)}$) and assigns the type of the founding individual to all members of the new $\alpha$-fraction. The lifted operator acts on the test function $\chi = fg$ by replacing $g(\xi_i)$ by $\mfh|_{B_\ell}$-averaged values for the reshuffled individuals, and replacing $f(v_i)$ by the type of the founding individual. This is encoded in~\eqref{e808} combined with the reshuffling operator~\eqref{e824}.

\medskip\noindent
{\bf Lifting from $\Omega_N^L$ to $\Omega_N^\infty$.}
The sum over all $(\ell,\xi) \in \Z \times \Omega_N^L$ in~\eqref{e804} defines the full block-resampling operator on $\Omega_N^L$.  As $L \to \infty$, the sum extends over all levels $\ell \ge -L$, and the rates $\lambda_\ell$ are summable by condition~\eqref{ag6} (this is equivalent to~\eqref{ag6}). Hence the operator $\mfA^{\rm Cann,b} = \sum_{L,\xi} \lambda_L A^{\rm Cann,b}_{(L,\xi)}$ is locally finite on $\Omega_N^\infty$ and well defined in the continuum limit. The key new ingredient compared to~\cite{gklimw} is precisely this block-resampling term, and the verification that it has the correct duality (block-coalescence on the dual side, as in~\eqref{ag17}) is carried out in Section~\ref{ss.prftech}.

\bibliographystyle{alpha}
\bibliography{gh2}

\end{document}